\documentclass[final,12pt,times]{elsarticle}
\usepackage{lineno}
\modulolinenumbers[200]
\usepackage{epstopdf}
\usepackage{color}
\usepackage{colortbl}
\usepackage{transparent}
\usepackage{graphicx}
\usepackage{import}
\usepackage{diagbox}
\usepackage{algorithmic}
\usepackage{algorithm}
\usepackage{appendix}
\usepackage{amsmath}
\usepackage{amscd}
\usepackage{mathrsfs}
\usepackage{graphicx,color}
\usepackage{arydshln}
\usepackage{multirow}
\usepackage{booktabs}
\usepackage{rotating}
\usepackage{setspace}
\usepackage{indentfirst}
\usepackage{cases}
\usepackage{epstopdf}
\usepackage{xcolor}

\usepackage{amssymb,color}
\usepackage{amsthm}
\usepackage{epsfig,subfigure}

\newcommand{\edge}[1]{\ar@{-}[#1]}
\numberwithin{equation}{section}

\biboptions{numbers,sort&compress}%

\begin{document}

\begin{frontmatter}

    \title{Bayesian Inversion with Neural Operator (BINO) for Modeling Subdiffusion: Forward and Inverse Problems}

	\author[1,4]{Xiong-bin Yan}
	\ead{yanxb2015@163.com}

	\author[1,2,3]{Zhi-Qin John Xu\corref{cor}}
	\ead{xuzhiqin@sjtu.edu.cn}

	\author[1,2,3,4]{Zheng Ma\corref{cor}}
	\ead{zhengma@sjtu.edu.cn}

	\cortext[cor]{Corresponding author}
    \address[1]{School of Mathematical Sciences, Shanghai Jiao Tong University}
    \address[2]{Institute of Natural Sciences, MOE-LSC, Shanghai Jiao Tong University}
    \address[3]{Qing Yuan Research Institute, Shanghai Jiao Tong University}
    \address[4]{CMA-Shanghai, Shanghai Jiao Tong University}

\begin{abstract}
Fractional diffusion equations have been an effective tool for modeling anomalous diffusion in complicated systems. However, traditional numerical methods require expensive computation cost and storage resources because of the memory effect brought by the convolution integral of time fractional derivative. We propose a Bayesian Inversion with Neural Operator (BINO) to overcome the difficulty in traditional methods as follows. We employ a deep operator network to learn the solution operators for the fractional diffusion equations, allowing us to swiftly and precisely solve a forward problem for given inputs (including fractional order, diffusion coefficient, source terms, etc.). In addition, we integrate the deep operator network with a Bayesian inversion method for modelling a problem by subdiffusion process and solving inverse subdiffusion problems, which reduces the time costs (without suffering from overwhelm storage resources) significantly. A large number of numerical experiments demonstrate that the operator learning method proposed in this work can efficiently solve the forward problems and Bayesian inverse problems of the subdiffusion equation.
\end{abstract}
\end{frontmatter}

\linenumbers

\section{Introduction}\label{Sec1}
Reactions and diffusion phenomena, across a wide range of applications including animal coat patterns and nerve cell signals, are modeled extensively by standard reaction-diffusion equations, which complies with a basic assumption that the diffusion obeys
the standard Brownian motion. The most distinctive characteristic of the standard Brownian motion is that the mean squared displacement $\langle x^2(t)\rangle$ of diffusing species follows a linear
dependence in time, i.e., $\langle x^2(t)\rangle\sim K_{1}t$. However, over the past few decades, it has been found that anomalous diffusion~\cite{Klafter_Sokolov_2005}, in which diffusive motion cannot be modeled as the standard Brownian motion, has occurred in many phenomena.
Unlike the standard diffusion, a hallmark of the anomalous diffusion is the non-linear power law growth of mean square displacement with time, i.e., $\langle x^2(t)\rangle\sim K_{\alpha}t^{\alpha}$, where
$\alpha>1$ is referred as superdiffusion and $0<\alpha<1$ is referred as subdiffusion. To derive the fractional diffusion equation, at a mesoscopic level, one of the most efficient methods is continuous time random walks (CTRWs). In the CRTWs, subdiffusion arises when the waiting time distribution is heavy-tailed (e.g., $\psi(t) \sim t^{-\alpha-1}$ with $\alpha\in (0, 1)$), where $\psi(t)$ denotes a waiting time of the particle jumping between two successive steps. And the associated probability density function of the particle appearing at time $t>0$ and spatial location $x\in \mathbb{R}^{d}$ satisfies the subdiffusion equation. Recently, the subdiffusion model has been found in many applications in physics, finance, and biology, including solute transport in heterogeneous media \cite{Berkowitz_Cortis_Dentz_Scher_2006, Dentz_Cortis_Scher_Berkowitz_2004}, thermal diffusion on fractal domains
\cite{Nigmatullin_1986}, protein transport within membranes
\cite{Kou_2008, Ritchie_Shan_Kondo_Iwasawa_Fujiwara_Kusumi_2005},
flow in highly heterogeneous aquifers \cite{Berkowitz_Klafter_Metzler_Scher_2002}.

In this paper, let $\Omega=(0,1)^2$, we consider the following subdiffusion problem:
\begin{align}\label{1.1}
	\left\{
	\begin{aligned}
		 & \partial_{0+}^{\alpha}u(x,t)= \nabla\cdot(a(x)\nabla u(x,t)) + c(x)u(x,t) + f(x),
		~ x\in\Omega,~0<t\leq T,                                                             \\
		 & u(x,0)=u_0(x),~x\in \Omega,                                                       \\
		 & u(x,t)=0,~x\in\partial \Omega,~ 0<t\leq T,
	\end{aligned}
	\right.
\end{align}
where $0<\alpha<1$ and $\partial_{0+}^\alpha u(x,t)$ is the  Caputo
left-sided fractional derivative defined by
\begin{equation*}
	\partial_{0+}^\alpha u(x,t)=\frac{1}{\Gamma(1-\alpha)}\int_0^t(t-s)^{-\alpha}\frac{\partial u}{\partial s}(x,s)ds,\quad t>0,
\end{equation*}
$a(x)$, $c(x)$, $f(x)$ and $u_{0}(x)$ represent the diffusion coefficient,
reaction coefficient, source and initial value, respectively.

It is widely known that it is very challenging to find the time-fractional diffusion equation's analytical solution, thus, it is very important to design numerical methods for solving it. There are currently many numerical methods for discretizing the above fractional equation. Compared with solving standard diffusion equations, much effort is always used to construct appropriate methods to discretize the time-fractional derivative. 
Existing schemes include: the $L_1$-type approximation \cite{Langlands_Henry_2005, Li_Xu_2009,Lin_Li_Xu_2011}, the Gr${\rm\ddot{u}}$nwald-Letnikov approximation
\cite{Yuste_Acedo_2005, Zeng_Li_Liu_Turner_2013, Chen_Liu_Turner_Anh_2007} and
Convolution quadrature \cite{Jin_Li_Zhou_2017, Jin_Lazarov_Zhou_2019}, and so on.
These methods of solving the problem (\ref{1.1}) require using and storing the result at all previous time steps, which means that obtaining its numerical solution requires a lot of time and computational resources. Additionally, in order to solve the inverse problems for the time-fractional diffusion equation, it is necessary to repeatedly solve the forward problems, which brings more difficulties and computational complexity of the inverse problem solutions. At the same time, with numerous applications of machine learning in solving PDEs, it is increasingly possible to solve the above problems with machine learning approach. In this paper, we seek a machine learning-based method to deal with these difficulties of the subdiffusion problem (\ref{1.1}).

Using the machine learning methods, particularly deep learning methods, to solve partial differential equations has grown rapidly. Here, we will discuss two major neural network-based approaches for solving partial differential equations. The first method is to use a deep neural network (DNN) to parameterize the differential equation solutions, so that the neural network's outputs satisfy the partial differential equations, and then use the strong or weak form of the differential equations as the loss function. Several examples are listed: Physics-informed neural network methods \cite{dissanayake1994neural,Raissi_Perdikaris_Karniadakis_2019,liu2020multi}, variational methods
\cite{E_Yu_2018,Liao_Ming_2021}, weak adversarial network methods
\cite{Zang_Bao_Ye_Zhou_2020, Bao_Ye_Zang_Zhou_2020}. Several works solve fractional partial differential equations by deep learning methods. For example, Pang et al. and Guo et al. in \cite{Pang_Lu_Karniadakis_2019,Guo_Wu_Yu_Zhou_2022} apply fPINN and Monte Carlo fPINN to solve the forward and inverse problems in fractional partial differential equations, respectively. However, the methods
mentioned above suffer from retraining a new neural network for every new input parameter. 
The second approach is to parameterize the solution
operator by a deep neural network, for example, deep convolutional encoder-decoder networks \cite{zhu_Zabaras_2018,Afshar_Bhatnagar_Pan_Duraisamy_Kaushik_2019}, deep operator networks (DeepONet) \cite{Lu_Jin_Karniadakis_2019,Jin_Meng_Lu_2022, Wang_Wang_Perdikaris_2021,Mao_lu_Marxen_Zaki_2021}, neural operator \cite{Li_Kovachki_Azizzadenesheli_liu_Bhattacharya_Stuart_Anandkumar_2021,Li_Kovachki_Azizzadenesheli_Liu_Bhattacharya_Anandkumar_2020,Li_Kovachki_Azizzadenesheli_Liu_Bhattacharya_Stuart_Anandkumar_2020_1}, and MOD-Net \cite{zhang_luo_zhang_e_xu_ma_2022}. These operator networks only need to be trained once. To obtain a prediction solution for a new parameter, it only needs a forward pass of the network, which saves a lot of computational time. Therefore, learning operator can be an appropriate method for efficiently solving the complex subdiffusion problem (\ref{1.1}).


In the paper, we propose a Bayesian Inversion with Neural Operator (BINO) approach for modeling subdiffusion problem and solving inverse subdiffusion equation (\ref{1.1}). We use data-driven method to train a DeepONet for solving forward problems of the subdiffusion equation (\ref{1.1}).
The operator network only needs to be trained once, and then it can repeatedly solve the forward problems for different parameters fast. Then, we incorporate the deep operator learning into the Bayesian inversion approach for inverse problems of the subdiffusion model, which allows efficiently performing the Bayesian inversion.
In summary, the main contributions of this paper are listed as follows:
\begin{itemize}
	\item We propose a BINO method for modelling a subdiffusion problem by identifying the fractional order of the equation (\ref{1.1}).
	\item We use BINO to recover the diffusion coefficient of subdiffusion equation (\ref{1.1}), which can be used to rapidly and efficiently implement the Markov chain Monte Carlo (MCMC) sampling in Bayesian inference.
\end{itemize}

The rest of the paper is structured as follows: In section
\ref{sec2}, we introduce the deep operator learning and use it to study the solution operator learning problems in subdiffusion (\ref{1.1}). In section \ref{sec4}, we combine the operator learning method with the Bayesian inference
method to solve several inverse problems in model (\ref{1.1}). Finally, we conclude the paper in section \ref{sec5}.

\section{Solving forward problems of subdiffusion by the deep operator networks.}\label{sec2}

Let $a\in C^1(\bar{\Omega})$, $a(x)>0$, $c\in C(\bar{\Omega})$, $c(x)\leq 0$, $x\in \Omega$ and $u_{0},~f\in L^2(\Omega)$. Then, by \cite{Sakamoto_Yamamoto_2011}, the solution of problem (\ref{1.1}) is unique and satisfies the following form
\begin{equation}\label{3.0}
	u(x,t) = \sum_{n=1}^{\infty}(u_0,\varphi_n)E_{\alpha,1}(-\lambda_nt^{\alpha})\varphi_n(x)+\sum_{n=1}^{\infty}\frac{1-E_{\alpha,1}(-\lambda_n t^{\alpha})}{\lambda_n}(f,\varphi_n)\varphi_n(x),
\end{equation}
where the notation $(\cdot,\cdot)$ denotes the scalar product of $L^2(\Omega)$ space. The $\lambda_n,~\varphi_n,~n=1,\cdots,\infty$ in (\ref{3.0}) are the eigenvalues and eigenfunctions of elliptic operator, which means
$-\nabla\cdot(a(x)\nabla \varphi_n(x)) - c(x)\varphi_n(x)=\lambda_n \varphi_n(x)$ and $\varphi_n(x)|_{\partial\Omega}=0$. The funciton $E_{\alpha,1}(-\lambda_n t^{\alpha})$ in (\ref{3.0}) denotes the Mittag-Leffler function \cite{Podlubny_1999}. By simplifying the formula (\ref{3.0}), it can be seen that
\begin{equation}\label{3.1}
	u(x,t)=\int_{\Omega}k_{1}(x,\xi,t)u_{0}(\xi)d\xi+
	\int_{\Omega}k_2(x,\xi,t)f(\xi)d\xi,
\end{equation}
where $k_{1}(x,\xi,t)=\sum_{n=1}^{\infty}E_{\alpha,1}(-\lambda_n t^{\alpha})\varphi_n(x)\varphi_n(\xi)$, $k_{2}(x,\xi,t)=\sum_{n=1}^{\infty}\frac{1-E_{\alpha,1}(-\lambda_n t^{\alpha})}{\lambda_n}\varphi_n(x)\varphi_n(\xi)$. From (\ref{3.1}), we can easily see that the solution to the problem (\ref{1.1}) depends on the fractional order $\alpha$, coefficients $a(x),~ c(x)$ as well as the initial value $u_0(x)$ and source term $f(x)$. Therefore, based on (\ref{3.1}), we will use the deep operator learning method to learn the solution of subdiffusion problem (\ref{1.1}).

\begin{figure}[H]
	\centering
	\subfigure[Unstacked DeepONet]{
		\begin{minipage}[t]{0.5\linewidth}
			\centering
			\includegraphics[width=3in]{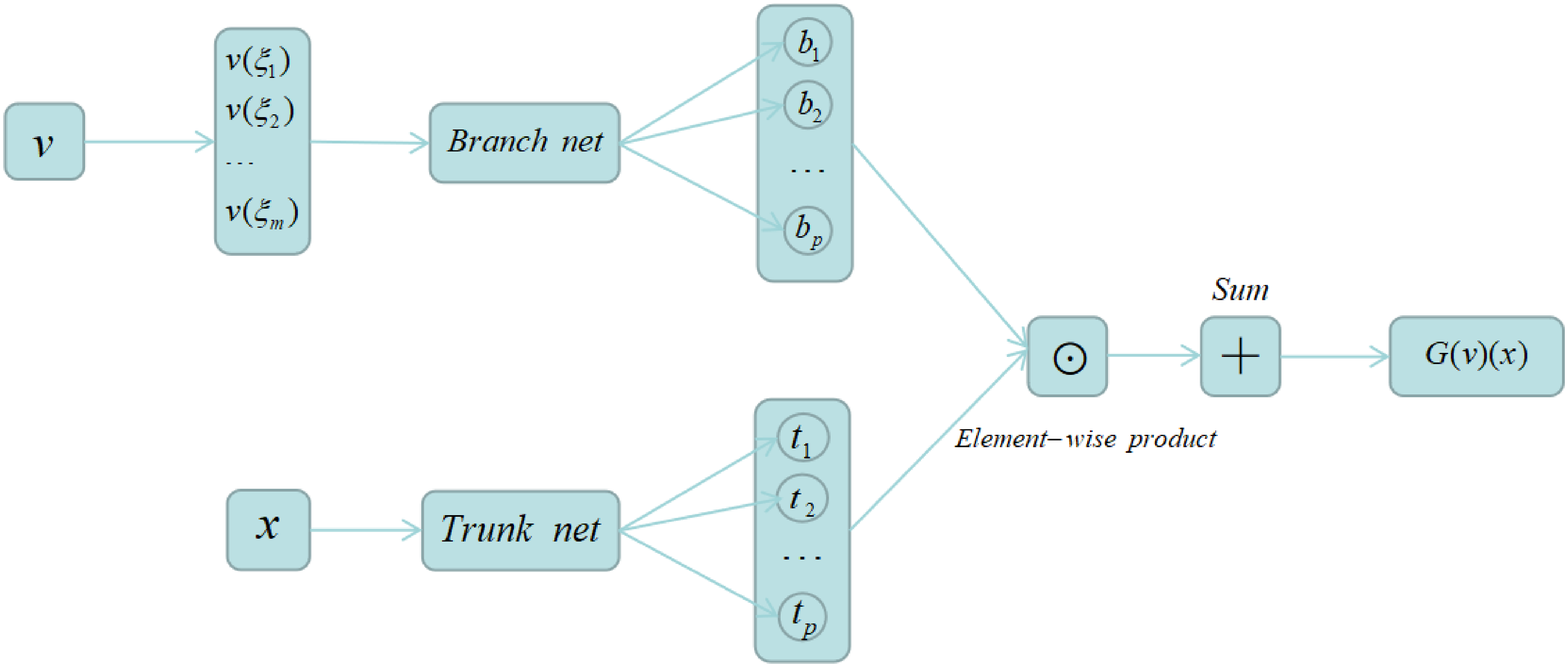}
		\end{minipage}
	}%
	\subfigure[MIONet]{
		\begin{minipage}[t]{0.5\linewidth}
			\centering
			\includegraphics[width=3in]{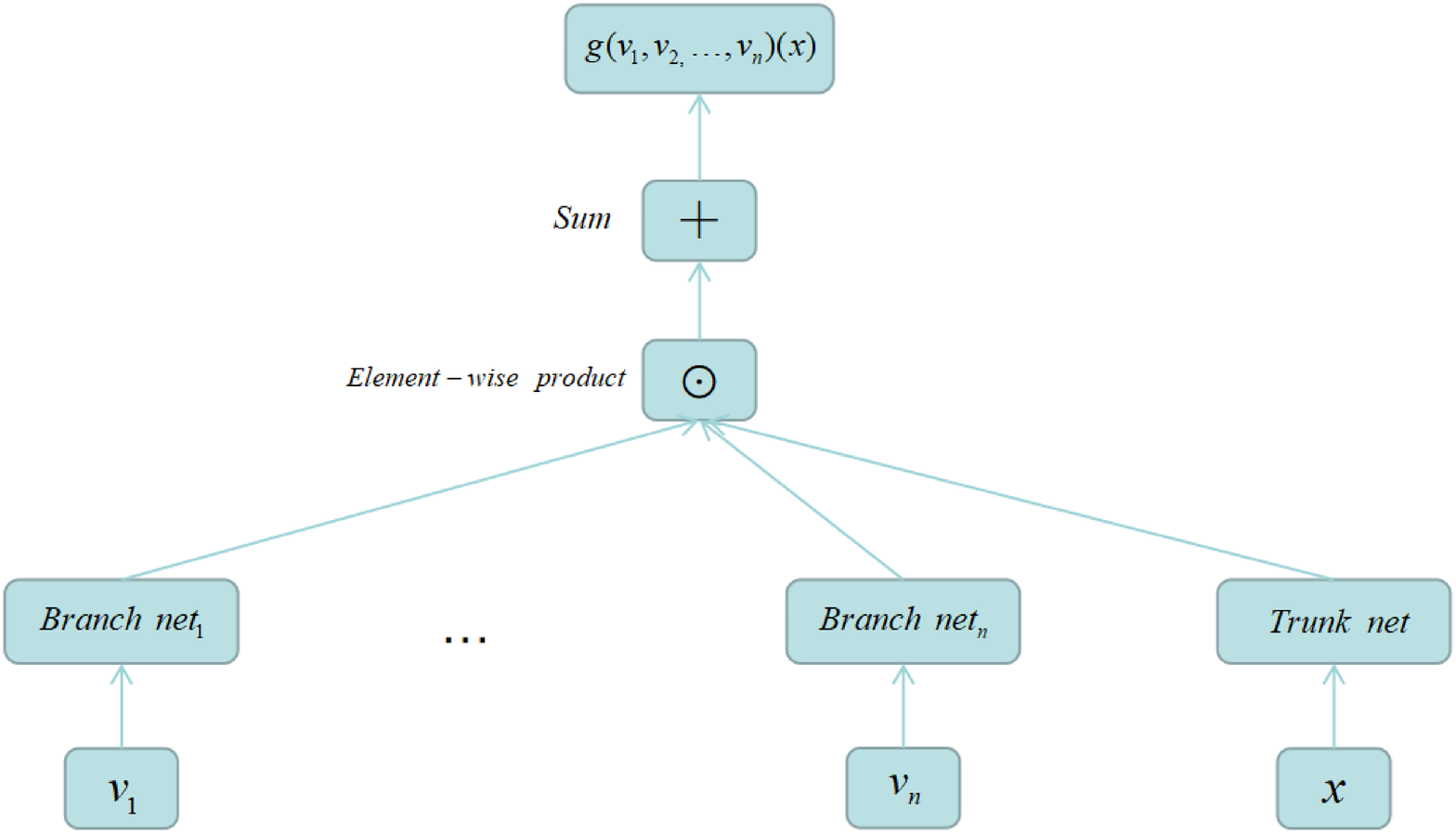}
		\end{minipage}
	}%
	\caption{The unstacked DeepONet and MIONet.}
	\label{fig_2.1}
\end{figure}

The deep operator networks is proposed by Lu et al. in \cite{Lu_Jin_Karniadakis_2019} and used to learn a  operator between infinite-dimensional function spaces. A specific neural network structure is shown in Figure \ref{fig_2.1} (a). The outputs of the network are aggregated by element-wise product and summation of outputs of two sub-networks, named branch network and trunk network, respectively. Further, in order to learn a multiple inputs map, Jin et al. in \cite{Jin_Meng_Lu_2022} proposes a multiple inputs operator network (shown in Figure \ref{fig_2.1} (b)), to learn the solution operator with multiple inputs. The operator network encodes the input functions by several branch nets and encodes the coordinate of the domain by trunk net, and all branch and trunk nets have the same dimensional outputs, and they are aggregated by element-wise product and summation.

In the next section, we will apply deep operator networks to learn the solution operators of the subdiffusion problem (\ref{1.1}). To sort out the experimental details clearly, we describe basic settings below.

We define a PDE solution operator $ G:v\in\mathcal{V}\to u\in\mathcal{U}$ as
\begin{equation*}
	G(v)=u,
\end{equation*}
where $\mathcal{V}$ is a parameter space and $\mathcal{U}$ is a solution space. Now, one can represent the solution map by deep operator network $G_{NN}(v;\theta)$, the $\theta$ denotes all trainable parameters. Define
\begin{equation*}
	l(v,\theta)=\frac{1}{P}\sum_{i=1}^{P}\big|G_{NN}(v;\theta)(x_i)-u(x_i)\big|^2,
\end{equation*}
where $u(x_i)$ denotes the PDE solution evaluated at location $x_i$. Then, the deep operator network can be trained
by minimizing the following loss function:
\begin{eqnarray*}
	\mathcal{L}(\theta)&=&\frac{1}{N}\sum_{j=1}^{N}l(v^{j},\theta)\\
	&=&\frac{1}{NP} \sum_{j=1}^{N}\sum_{i=1}^{P}\big| G_{NN}(v^{j};\theta)(x_{i})-u(v^{j})(x_i)\big|^2.
\end{eqnarray*}

We choose fully connected neural networks as the approximate network. Unless otherwise specified, we use four hidden layers with size $128$ for all branch nets and trunk nets. In order to measure the accuracy of the approximate solutions with respect to the reference solutions, we calculate the relative $l_2$ error between the approximate solutions $\hat{u}_{a}$ and the reference solutions $\hat{u}_r$ defined by 
\begin{equation*}
    {\rm relative}~l_2~{\rm error} := \sqrt{\frac{\sum_{i\ge 1}|\hat{u}_{a}(x_i)-\hat{u}_r(x_i)|^2}{\sum_{i\ge 1}|\hat{u}_r(x_i)|^2}}.
\end{equation*}
Based on the above introduction, we will study several kinds of operator learning problems in subdiffusion model (\ref{1.1}).

\subsection{Task 1. Learning solution operator mapping $(\alpha,a)$ to solution $u$.}\label{sub_sec3.1}
We begin with an example in learning the solution operator that maps fractional order $\alpha$ and diffusion coefficient $a(x,y)$ to the solution $u(x,y,t)$ in (\ref{3.1}), that is, learning the following map
\begin{equation}\label{3.1.1}
	G:(\alpha, a)\to u(x,y,t), ~(x,y,t)\in\bar{\Omega}\times[0,T],
\end{equation}
where we denote the approximate deep operator network as $G_{NN}^{(\alpha,a)}$. To obtain the training/testing dataset, we fix
$c(x,y)=-xy-4$, $u_{0}(x,y) = 6\sin(2\pi x)\sin(3\pi y)$ and source $f(x,y)=\sin(3\pi x)\sin(\pi y)+6\exp(x^2+y^2)$. For the diffusion coefficient $a(x, y)$, we randomly sample different functions from a mean-$a_{0}$ GRF:
\begin{equation*}
	a \sim \mathcal{G} (a_{0}, k_{l}(x^{[1]}, x^{[2]})),
\end{equation*}
where the covariance kernel $k_l(x^{[1]}, x^{[2]})=\exp(-\|x^{[1]} - x^{[2]} \|^2/{2l^2})$ is the radial-basis function (RBF) kernel with a length-scale parameter $l > 0$, and $x^{[1]}=(x_1,y_1)$, $x^{[2]}=(x_2, y_2)$.
The length-scale $l$ determines the smoothness of the sampled function,
and a larger $l$ leads to a smoother diffusion coefficient $a$. In this section, we take $l=0.3$, $a_{0}=5$. To get the training/testing dataset of $\alpha$, we first divide $[\epsilon,1-\epsilon]$ into 20 equal parts, and then randomly select a sample from these points, where we take the $\epsilon=0.001$.
For the corresponding numerical solution $u$ which obtained by a $L_{1}$-type finite difference scheme \cite{Sun_Wu_2006}, the grid resolutions of time and space variables in the finite difference method are 51 and $101\times 101$.

\begin{figure}[H]
	\centering
	\subfigure[Training loss]
	{\epsfxsize 0.49\hsize \epsfbox{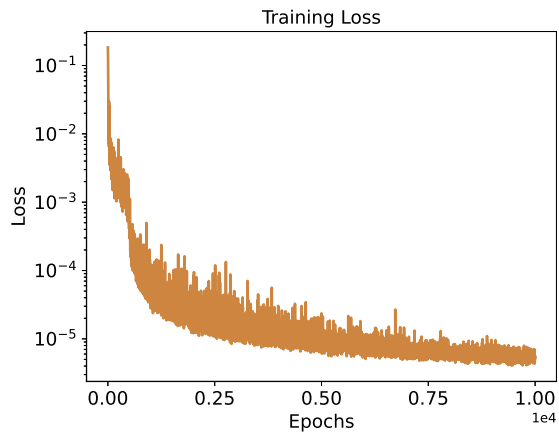}}
	\caption{Training loss of the deep operator network $G_{NN}^{(\alpha,a)}$ for 10000 epochs. The optimizer is Adam with learning rate $1\times 10^{-4}$.}\label{fig2.1_loss}
\end{figure}

\begin{figure}[H]
	\centering
	{\epsfxsize 0.3\hsize \epsfbox{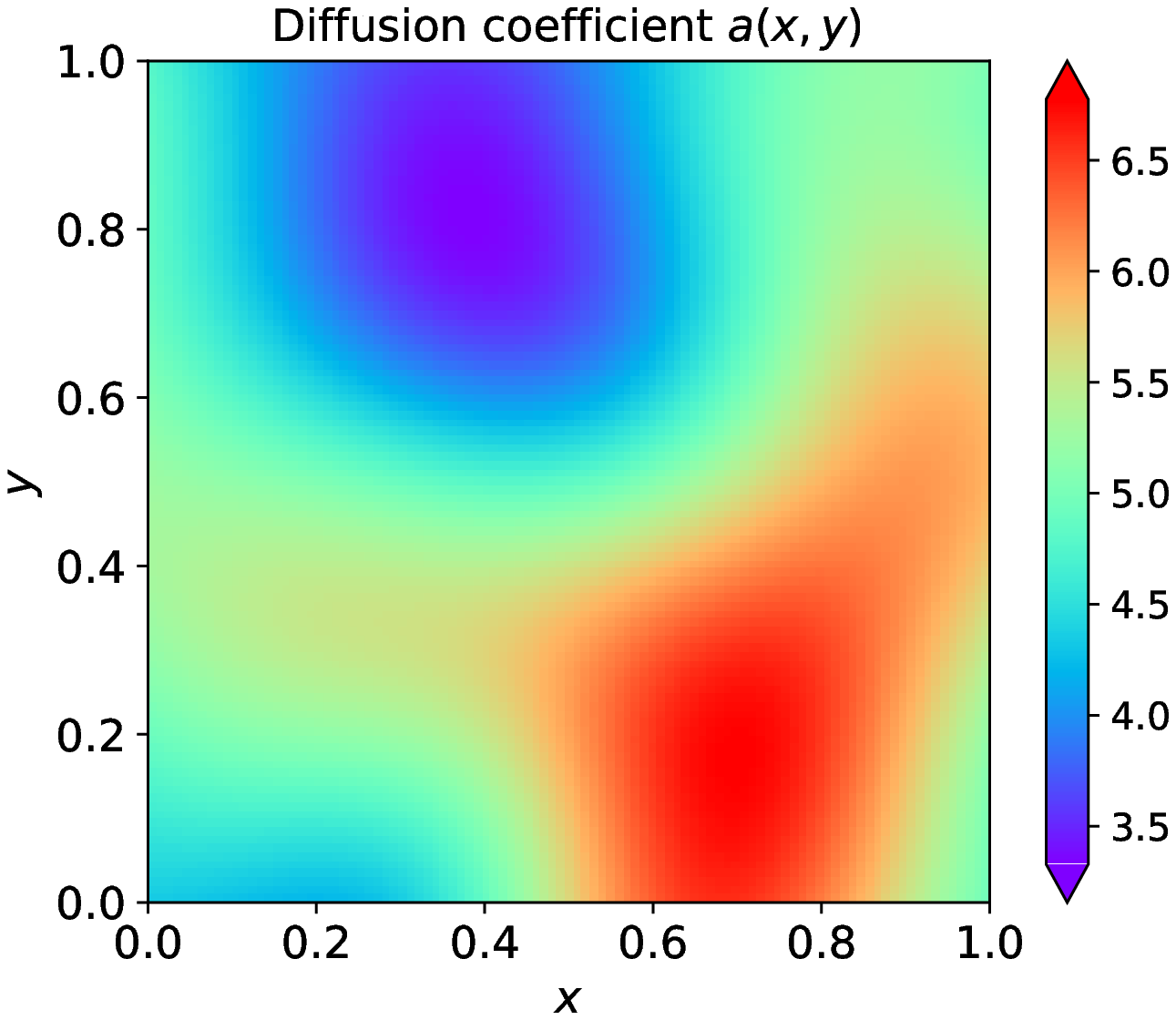}}
	\caption{Test sample 1 of the task 1. The input diffusion coefficient $a(x,y)$ of the operator network $G_{NN}^{(\alpha,a)}$.}\label{fig2.1_sample1_diffu}
\end{figure}

\begin{figure}[H]
	\centering
	\subfigure[FDM solution at $t=0$]{
		\begin{minipage}[t]{0.3\linewidth}
			\centering
			\includegraphics[width=2in]{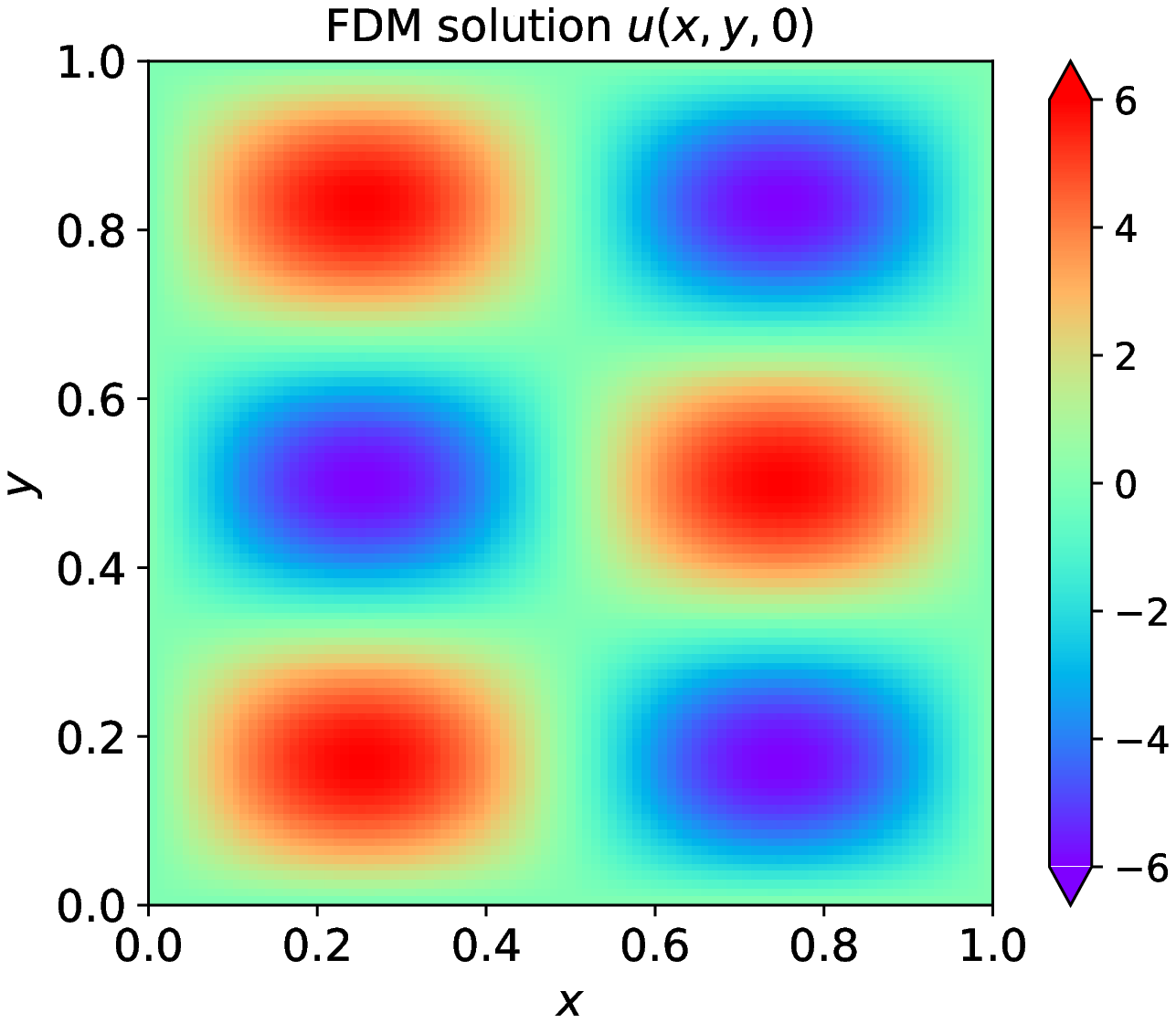}
		\end{minipage}
	}%
	\subfigure[$G_{NN}^{(\alpha,a)}$ solution at $t=0$]{
		\begin{minipage}[t]{0.3\linewidth}
			\centering
			\includegraphics[width=2in]{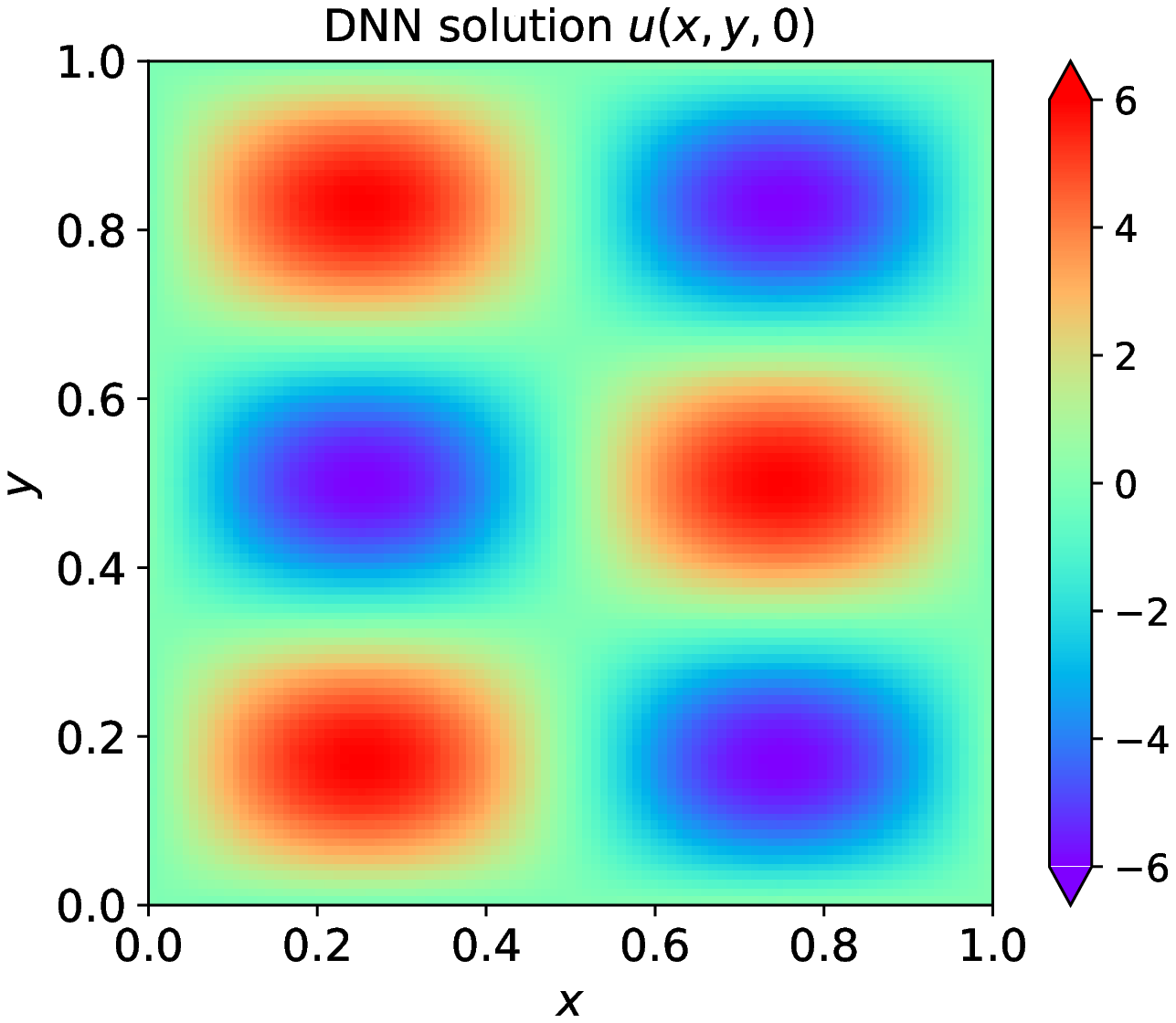}
		\end{minipage}
	}%
	\subfigure[Point-wise errors at $t=0$]{
		\begin{minipage}[t]{0.3\linewidth}
			\centering
			\includegraphics[width=2in]{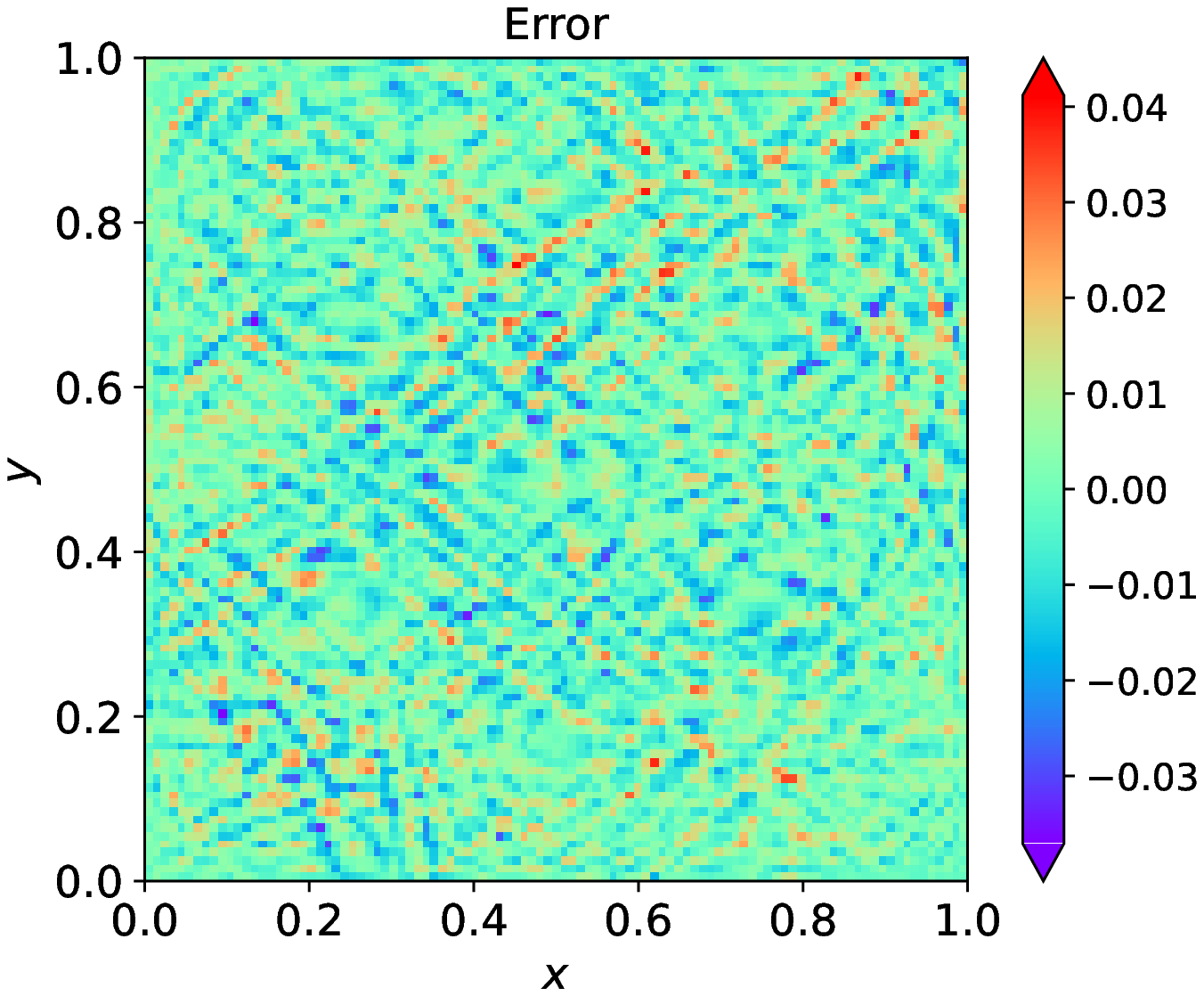}
		\end{minipage}
	}%

	\subfigure[FDM solution at $t=0.04$]{
		\begin{minipage}[t]{0.3\linewidth}
			\centering
			\includegraphics[width=2in]{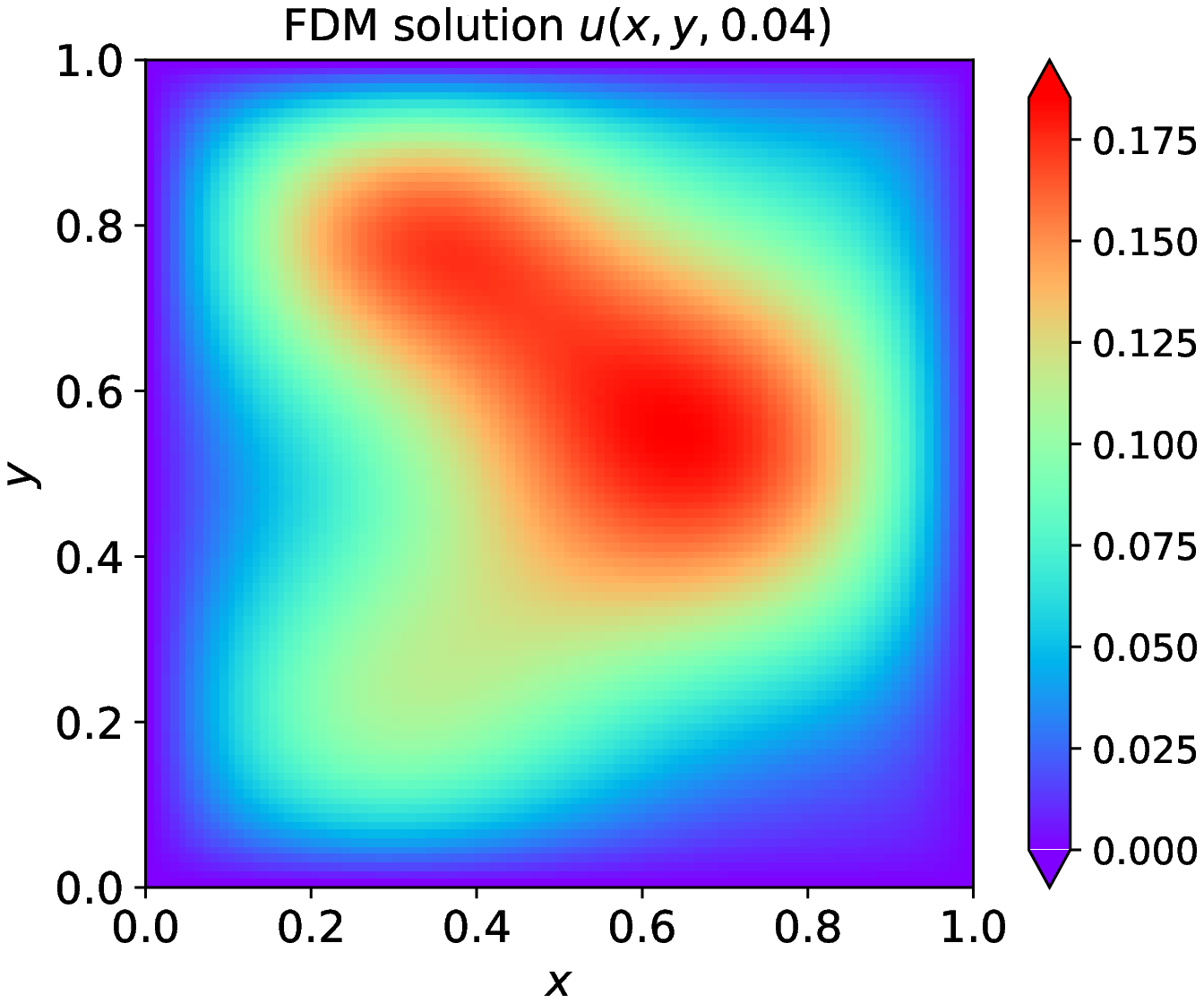}
		\end{minipage}
	}%
	\subfigure[$G_{NN}^{(\alpha,a)}$ solution at $t=0.04$]{
		\begin{minipage}[t]{0.3\linewidth}
			\centering
			\includegraphics[width=2in]{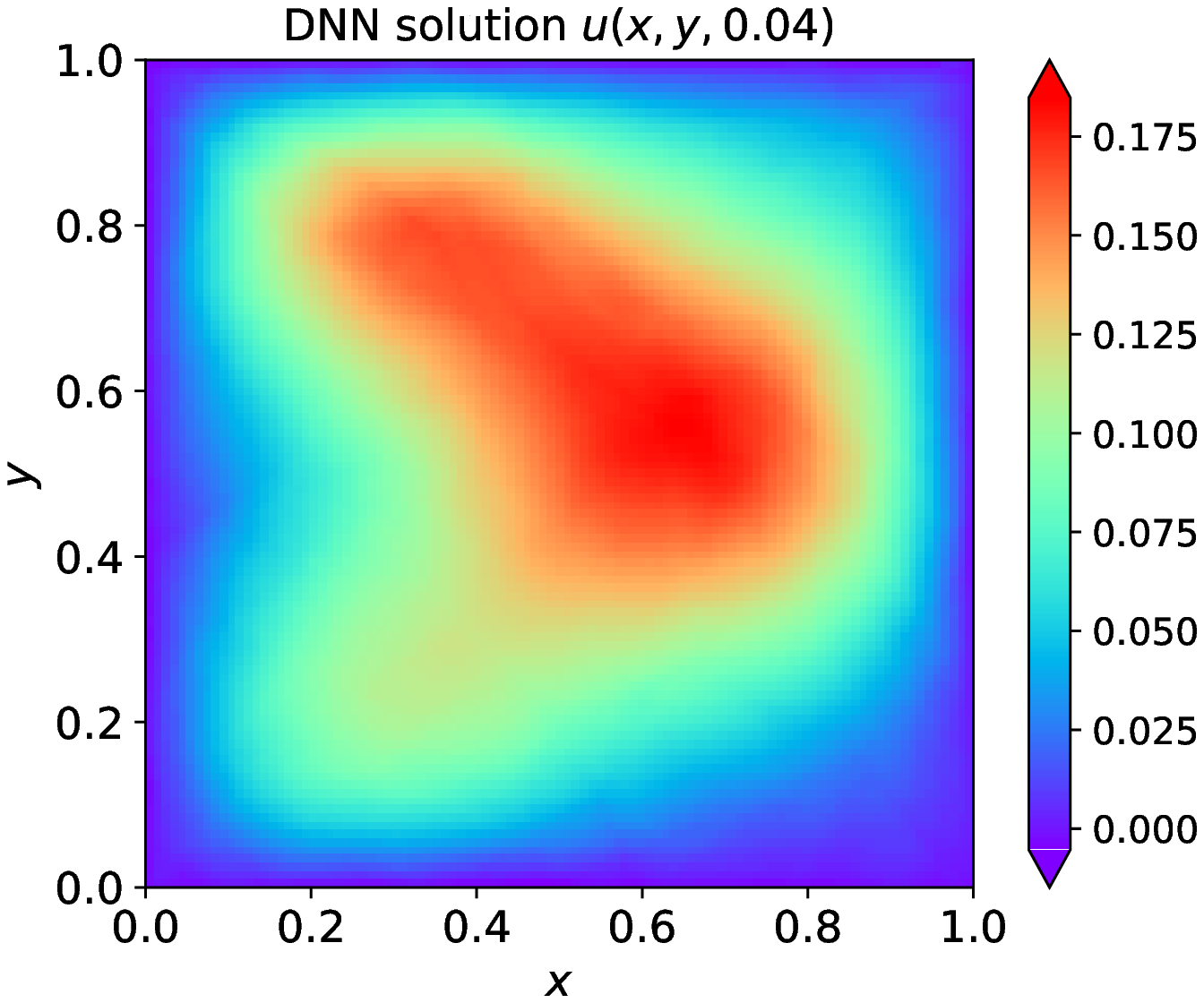}
		\end{minipage}
	}%
	\subfigure[Point-wise errors at $t=0.04$]{
		\begin{minipage}[t]{0.3\linewidth}
			\centering
			\includegraphics[width=2in]{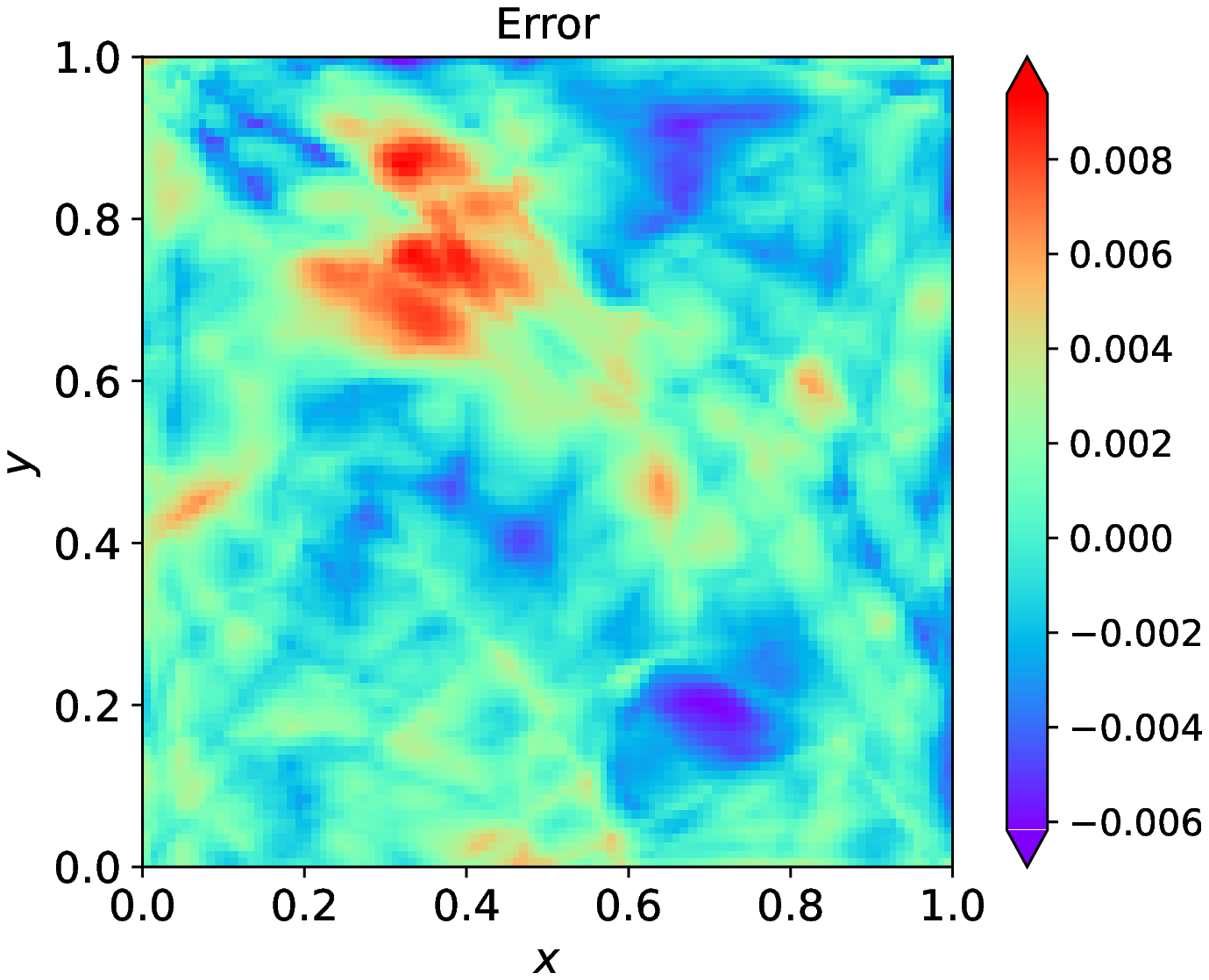}
		\end{minipage}
	}%

	\subfigure[FDM solution at $t=1$]{
		\begin{minipage}[t]{0.3\linewidth}
			\centering
			\includegraphics[width=2in]{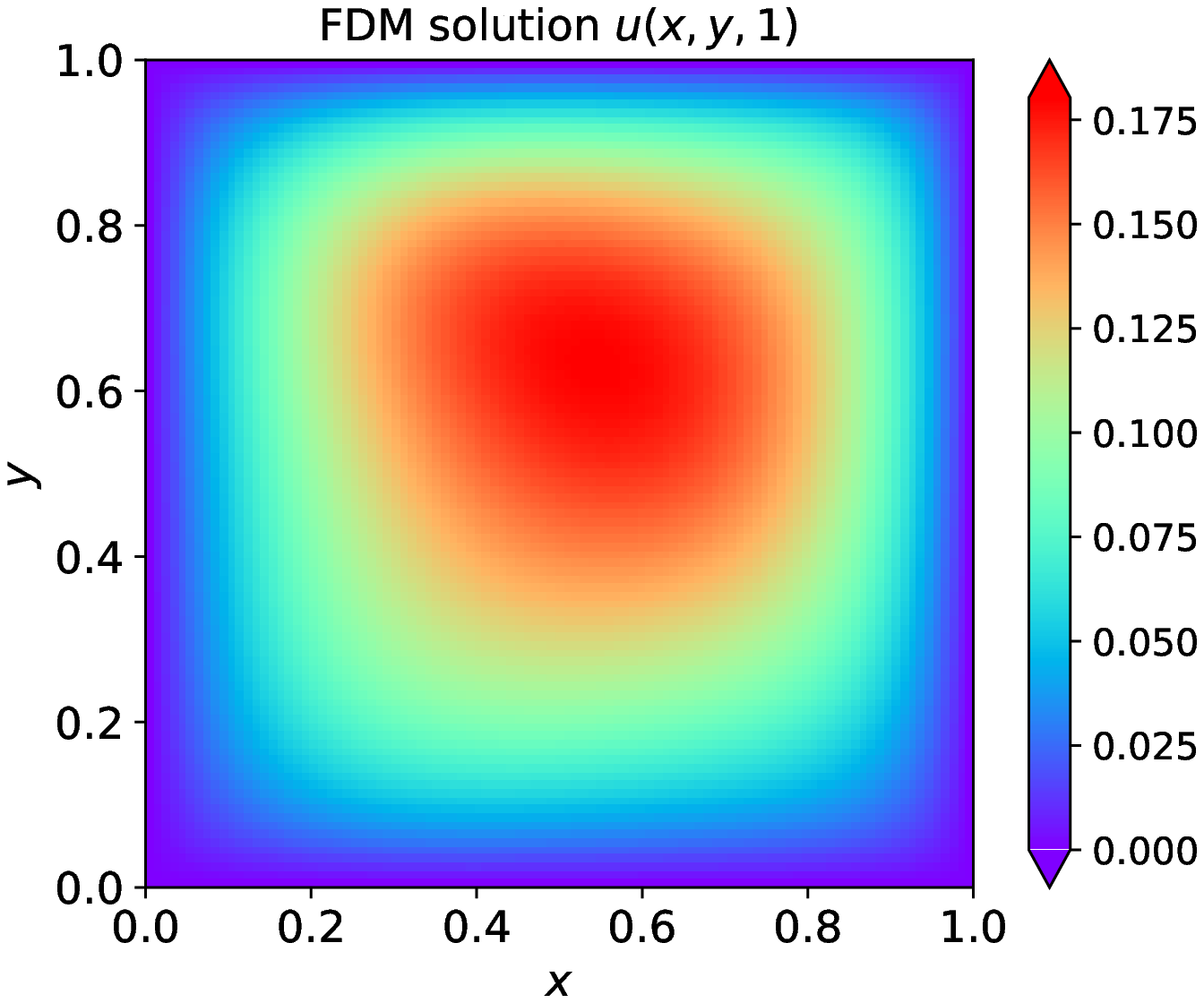}
		\end{minipage}
	}%
	\subfigure[$G_{NN}^{(\alpha,a)}$ solution at $t=1$]{
		\begin{minipage}[t]{0.3\linewidth}
			\centering
			\includegraphics[width=2in]{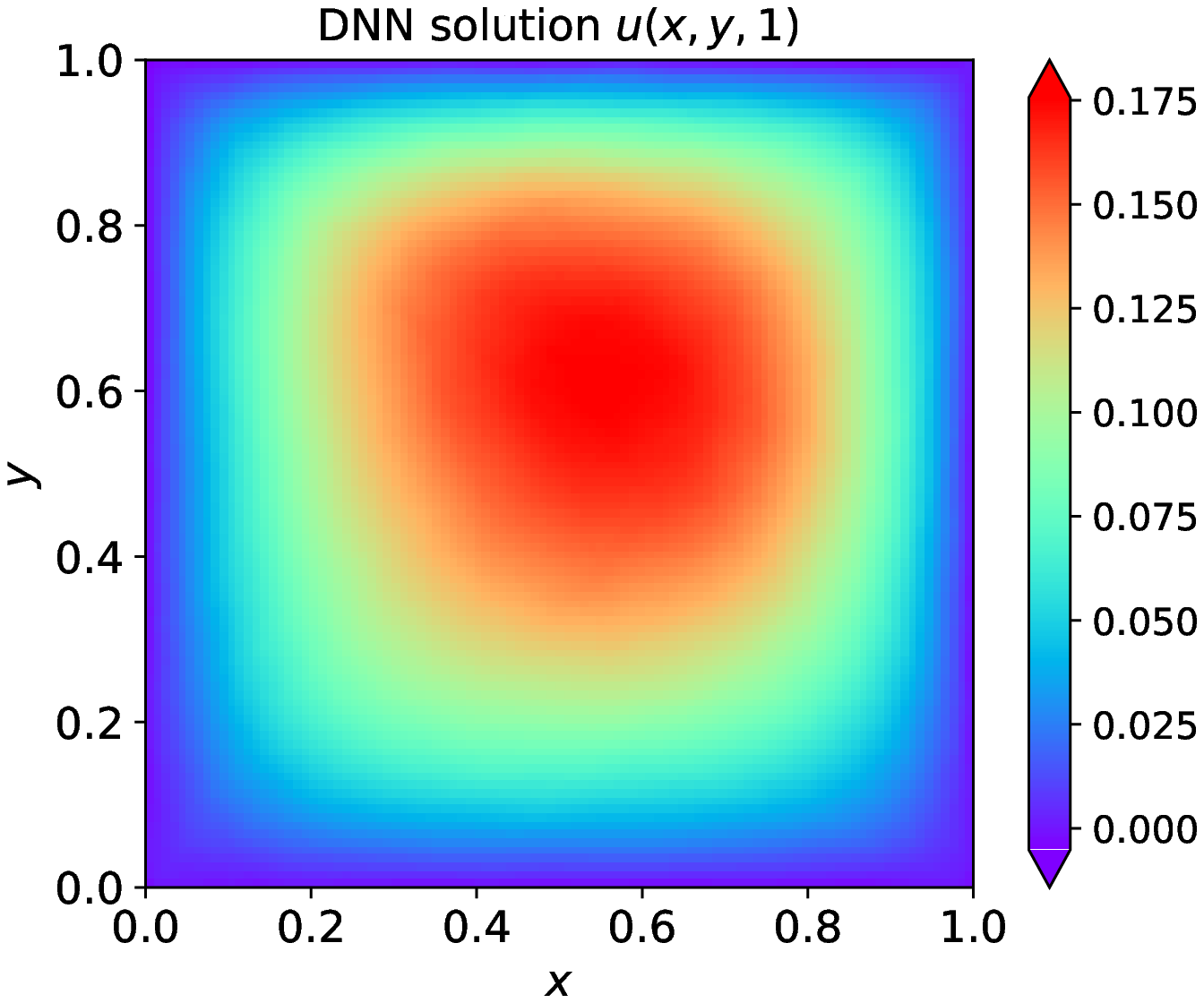}
		\end{minipage}
	}%
	\subfigure[Point-wise errors at $t=1$]{
		\begin{minipage}[t]{0.3\linewidth}
			\centering
			\includegraphics[width=2in]{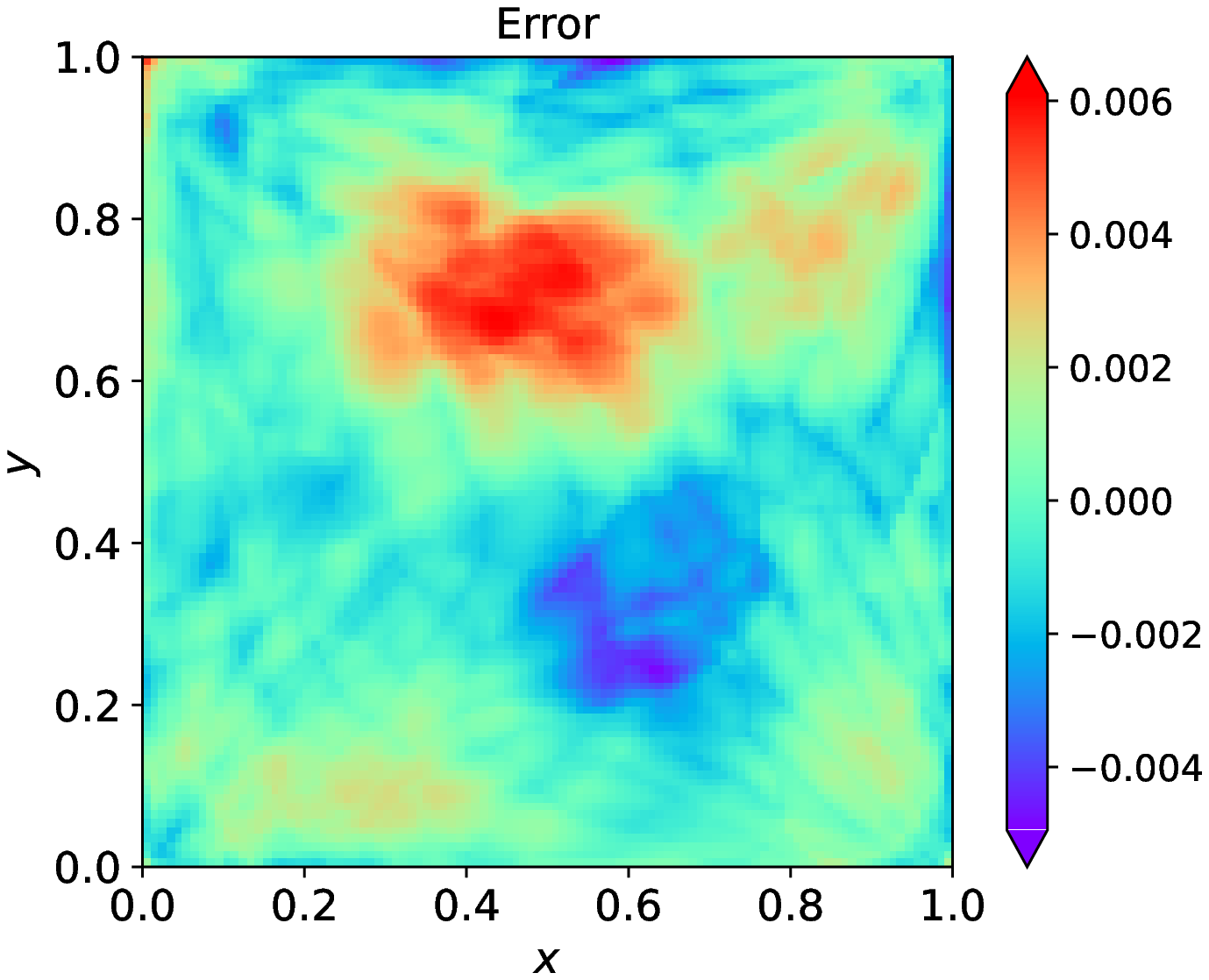}
		\end{minipage}
	}%
	\centering
	\caption{Test sample 1 of the task 1. Comparison between the finite difference solution and the outputs of the deep operator network $G_{NN}^{(\alpha,a)}$ on $101\times 101$ grid points corresponding to the fixed reaction coefficient $c(x,y)=-xy-4$, initial value $u_0(x,y)=6\sin(2\pi x)\sin(3\pi y)$, source $f(x,y)=\sin(3\pi x)\sin(\pi y)+6\exp(x^2+y^2)$, and the fractional order $\alpha=0.6$.}
	\label{fig2.1_sample1}
\end{figure}

\begin{figure}[H]
	\centering
	{\epsfxsize 0.3\hsize \epsfbox{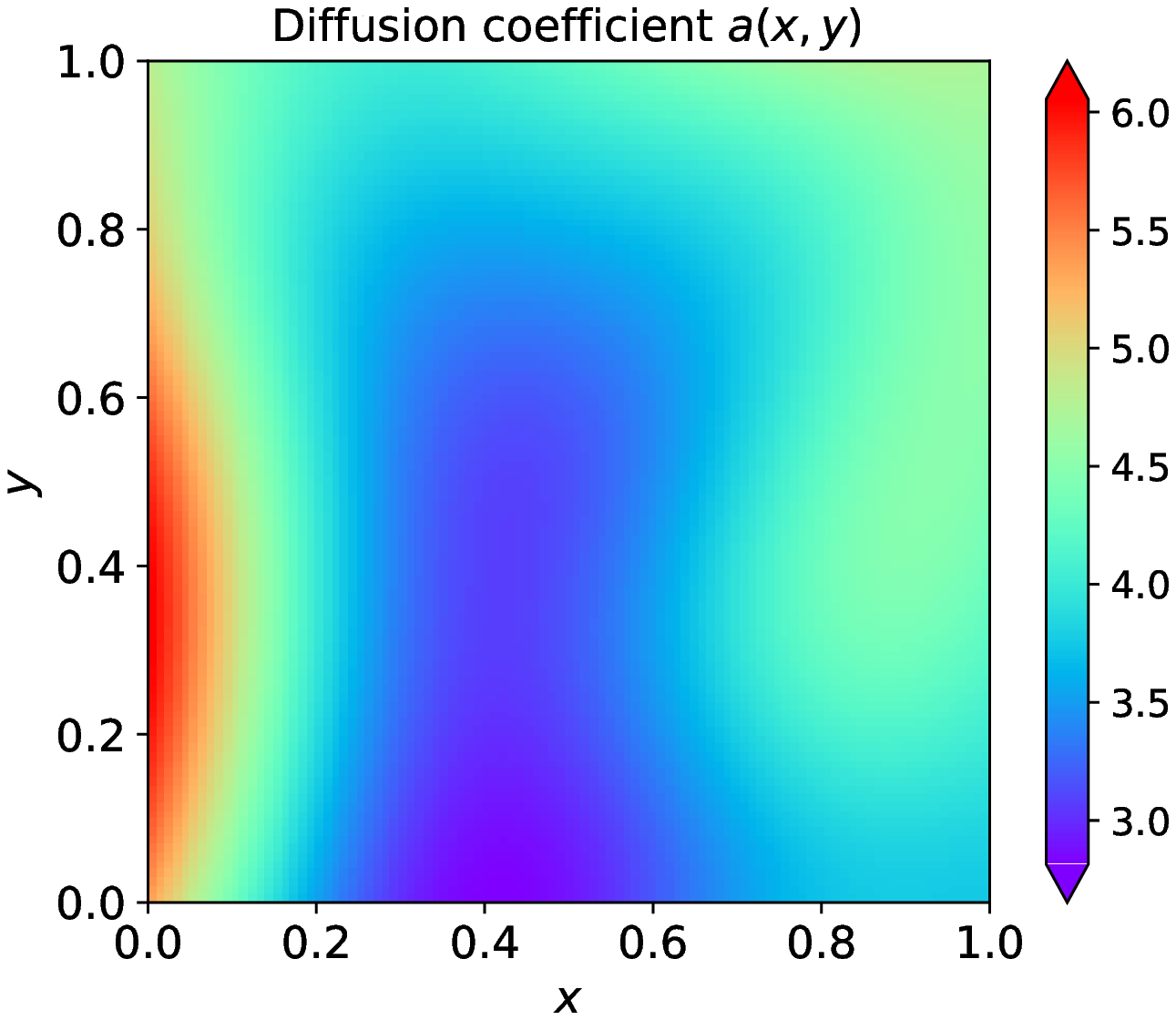}}
	\caption{Test sample 2 of the task 1. The input diffusion coefficient $a(x,y)$ of the operator network $G_{NN}^{(\alpha,a)}$.}\label{fig2.1_sample2_diffu}
\end{figure}

\begin{figure}[H]
	\centering
	\subfigure[FDM solution at $t=0$]{
		\begin{minipage}[t]{0.3\linewidth}
			\centering
			\includegraphics[width=2in]{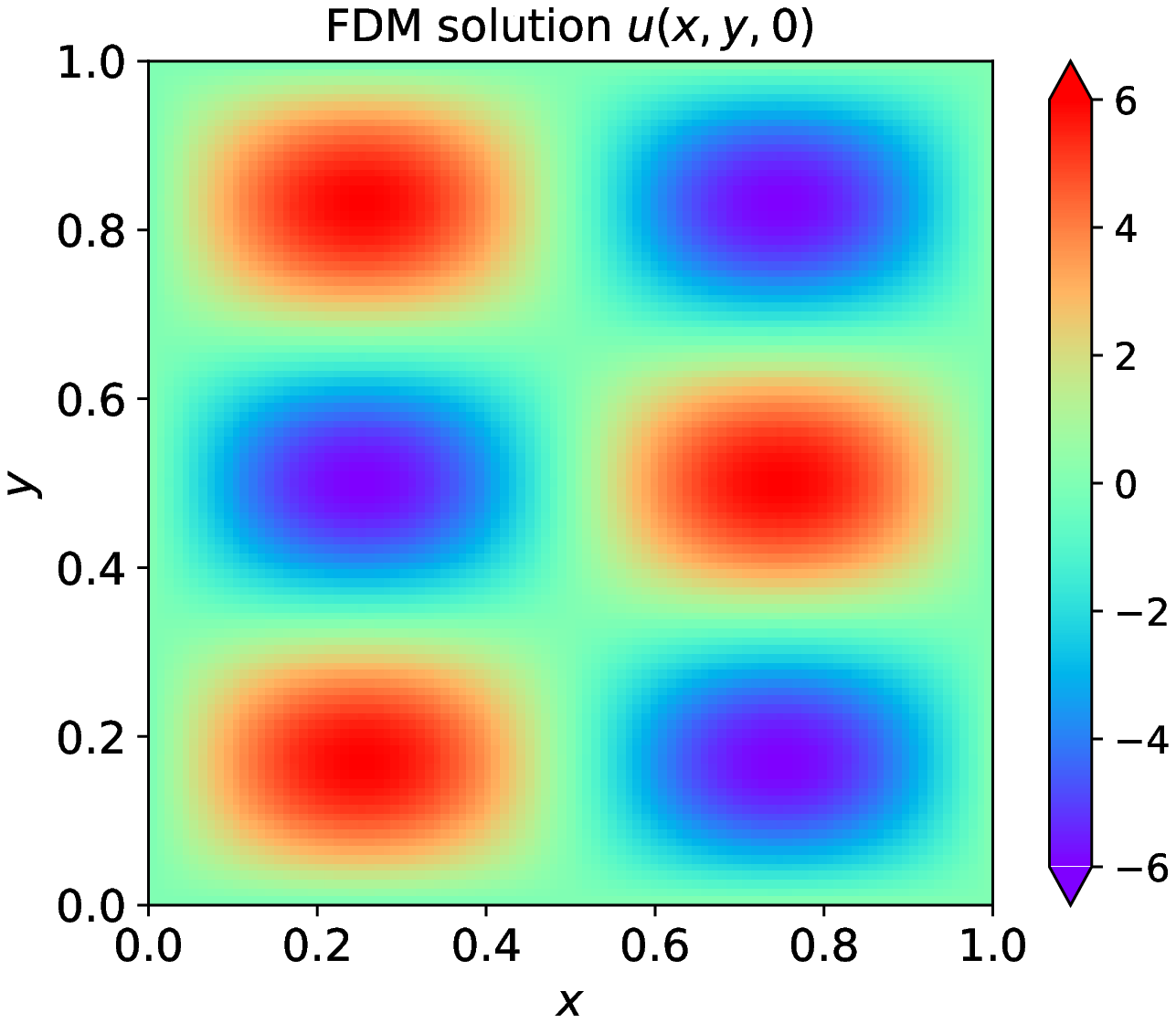}
		\end{minipage}
	}%
	\subfigure[$G_{NN}^{(\alpha,a)}$ solution at $t=0$]{
		\begin{minipage}[t]{0.3\linewidth}
			\centering
			\includegraphics[width=2in]{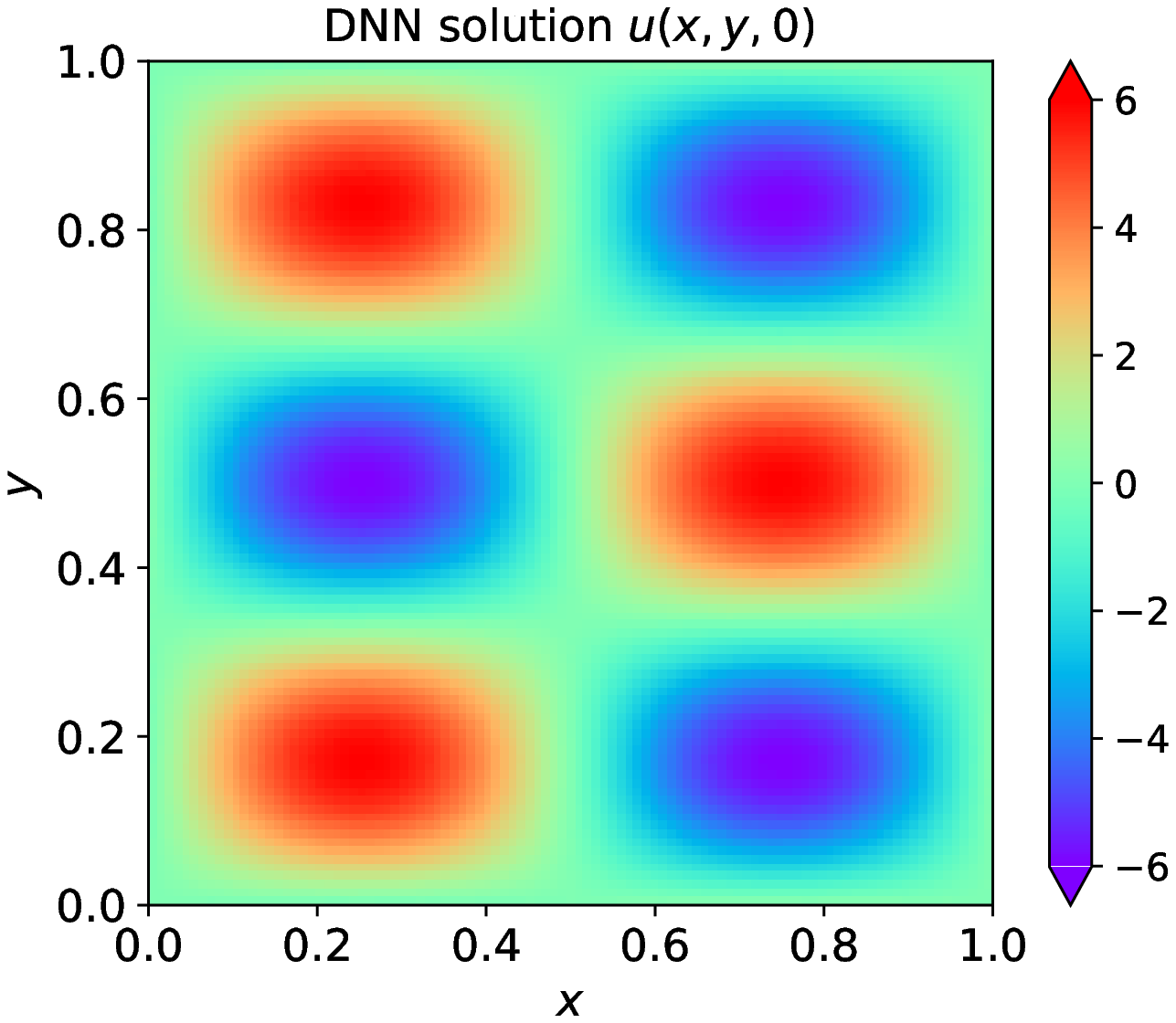}
		\end{minipage}
	}%
	\subfigure[Point-wise errors at $t=0$]{
		\begin{minipage}[t]{0.3\linewidth}
			\centering
			\includegraphics[width=2in]{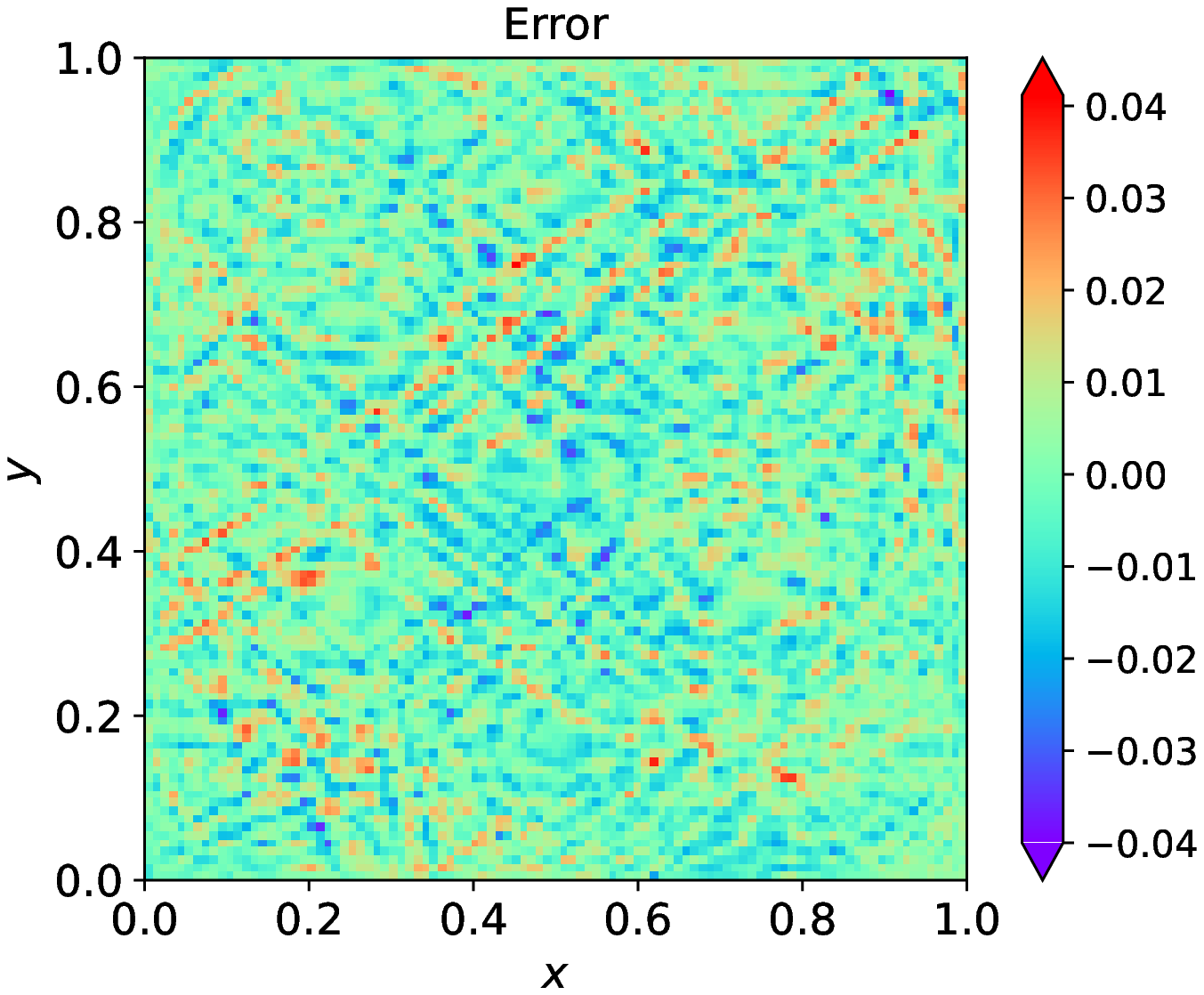}
		\end{minipage}
	}%

	\subfigure[FDM solution at $t=0.04$]{
		\begin{minipage}[t]{0.3\linewidth}
			\centering
			\includegraphics[width=2in]{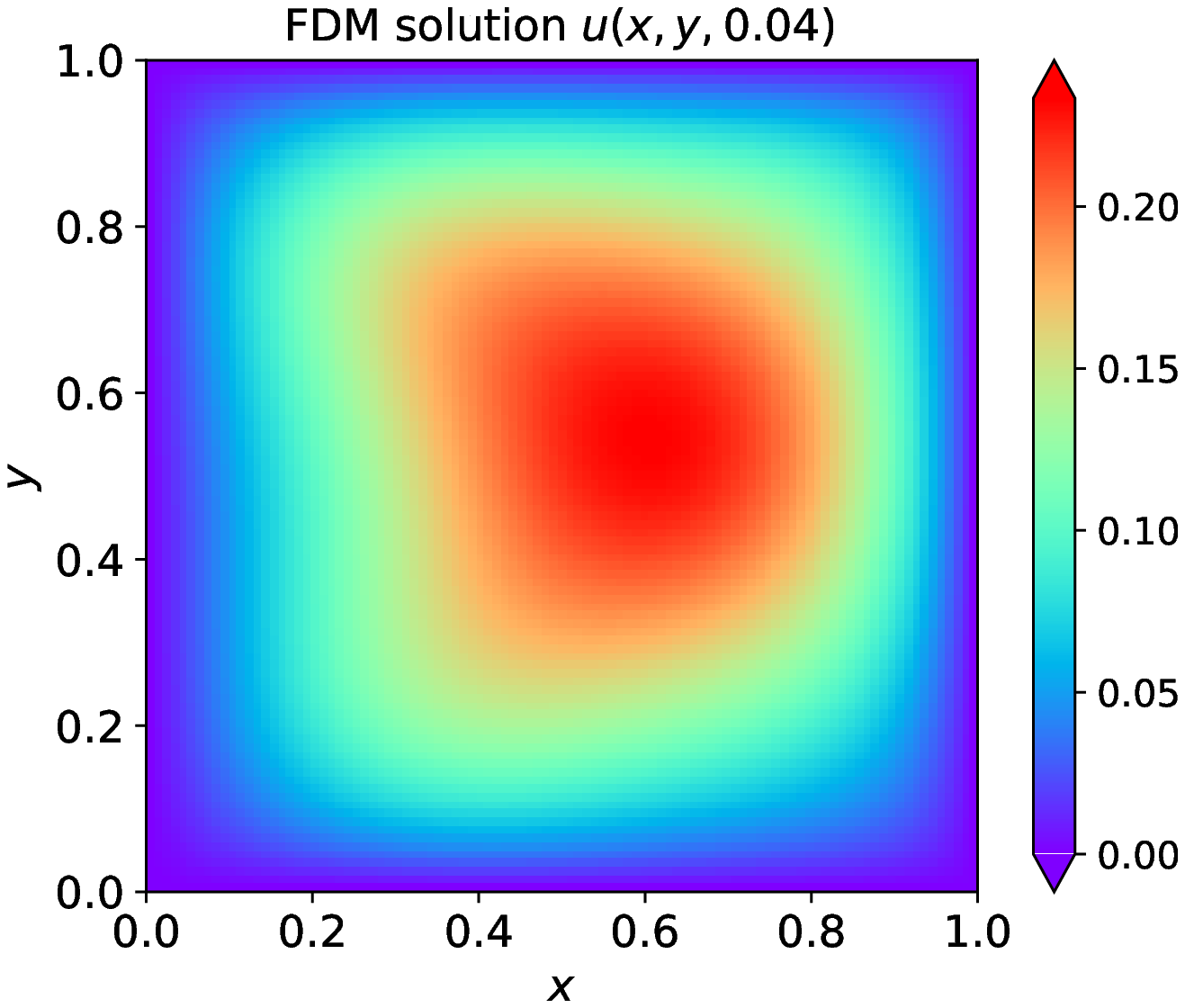}
		\end{minipage}
	}%
	\subfigure[$G_{NN}^{(\alpha,a)}$ solution at $t=0.04$]{
		\begin{minipage}[t]{0.3\linewidth}
			\centering
			\includegraphics[width=2in]{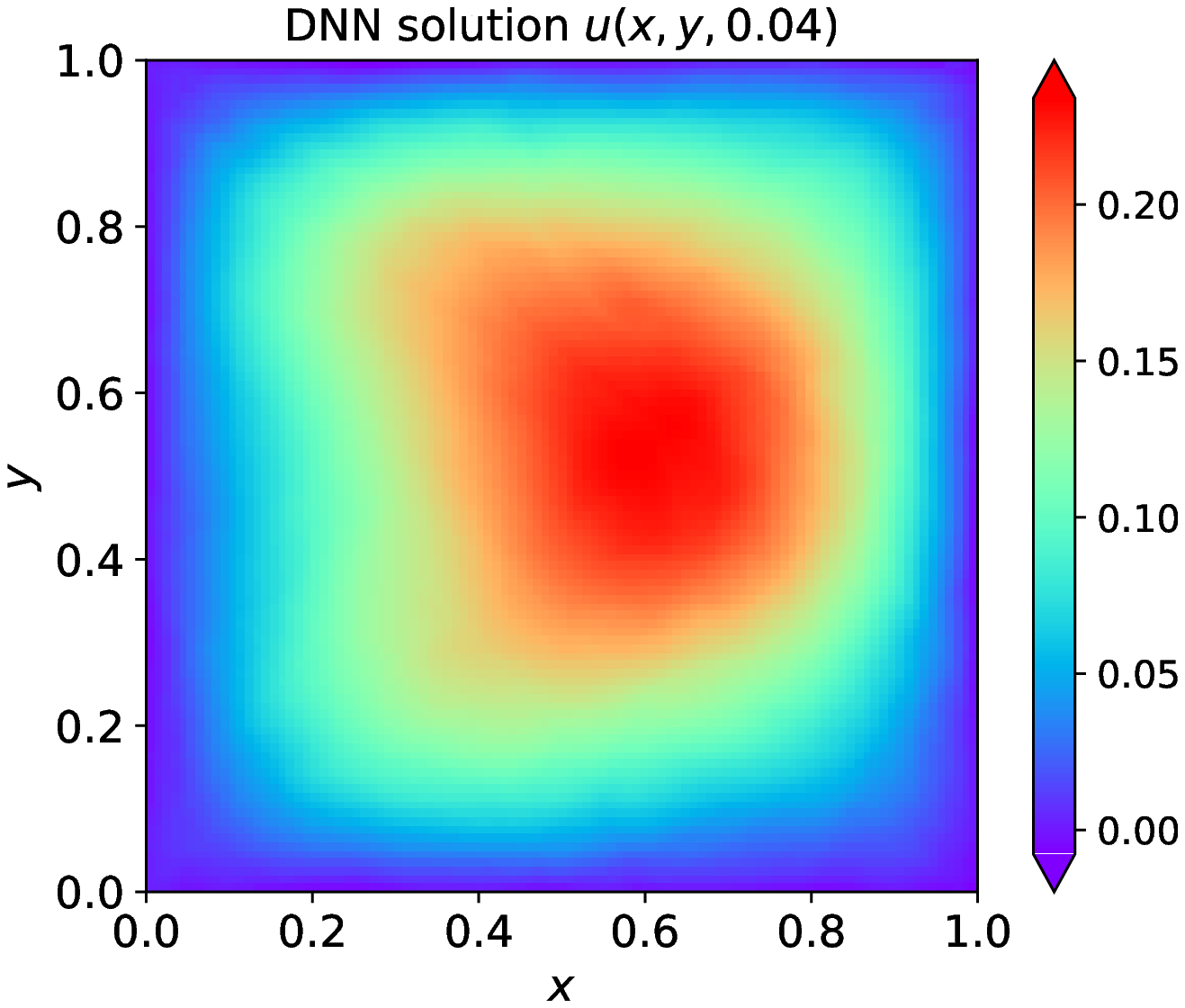}
		\end{minipage}
	}%
	\subfigure[Point-wise errors at $t=0.04$]{
		\begin{minipage}[t]{0.3\linewidth}
			\centering
			\includegraphics[width=2in]{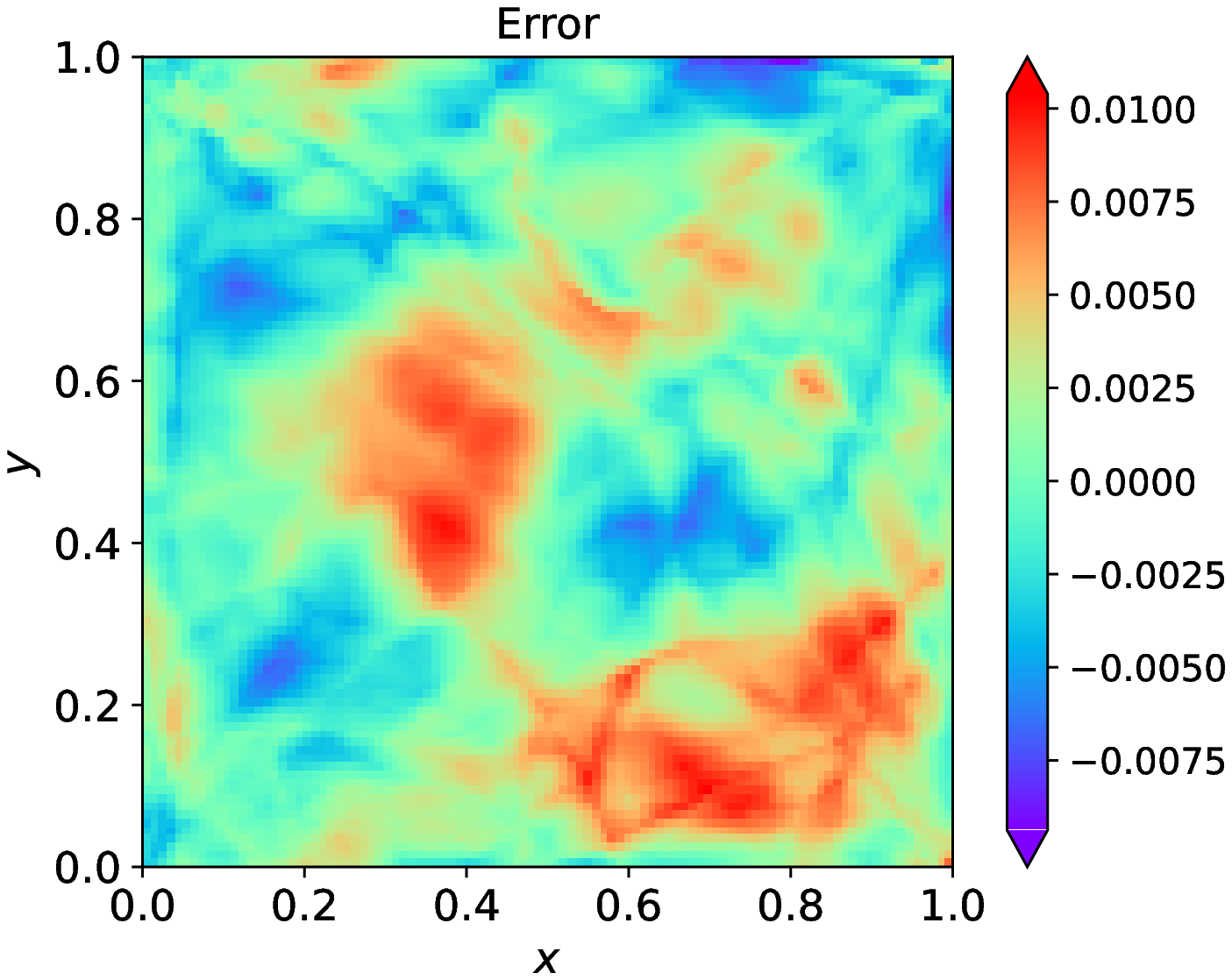}
		\end{minipage}
	}%

	\subfigure[FDM solution at $t=1$]{
		\begin{minipage}[t]{0.3\linewidth}
			\centering
			\includegraphics[width=2in]{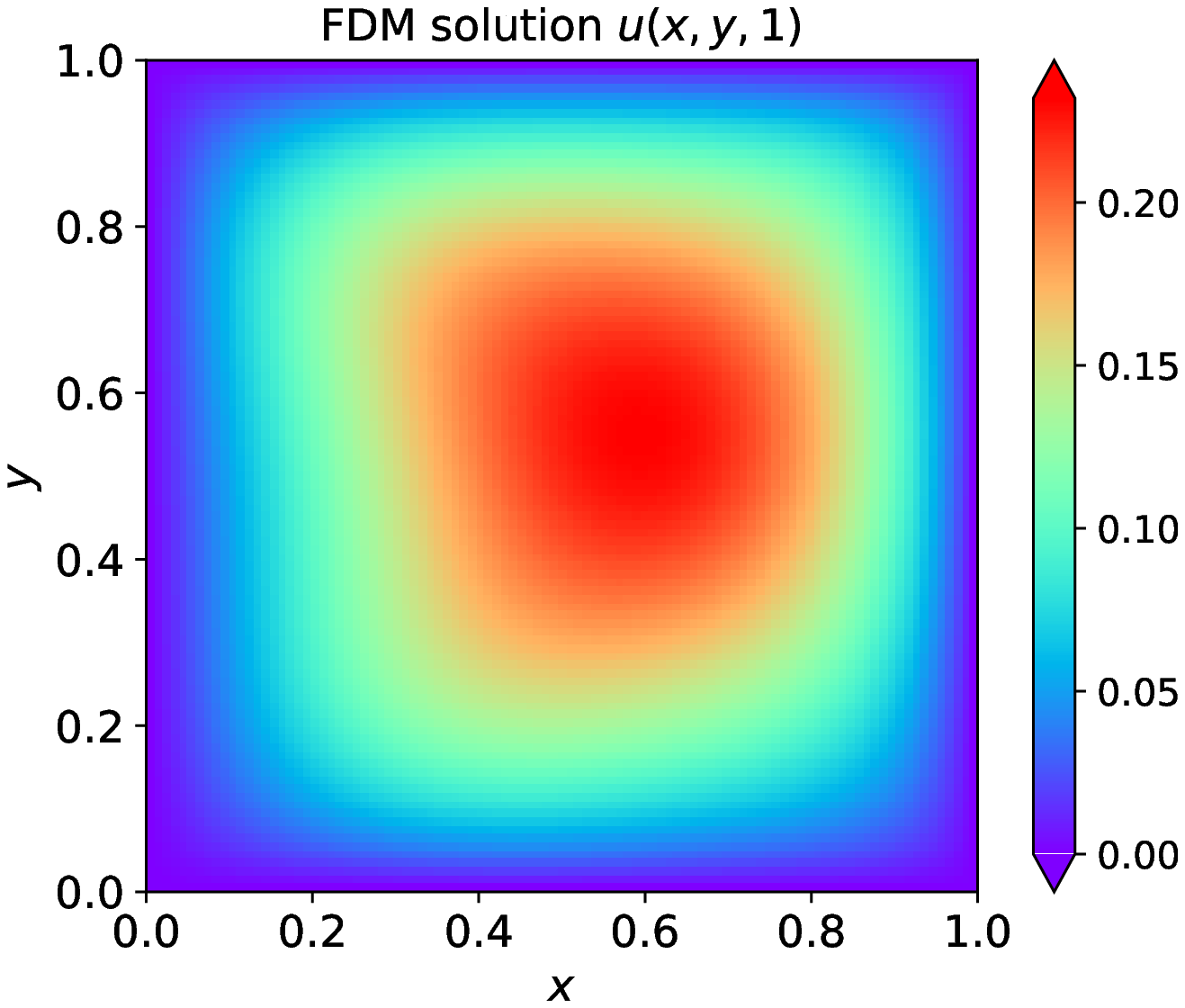}
		\end{minipage}
	}%
	\subfigure[$G_{NN}^{(\alpha,a)}$ solution at $t=1$]{
		\begin{minipage}[t]{0.3\linewidth}
			\centering
			\includegraphics[width=2in]{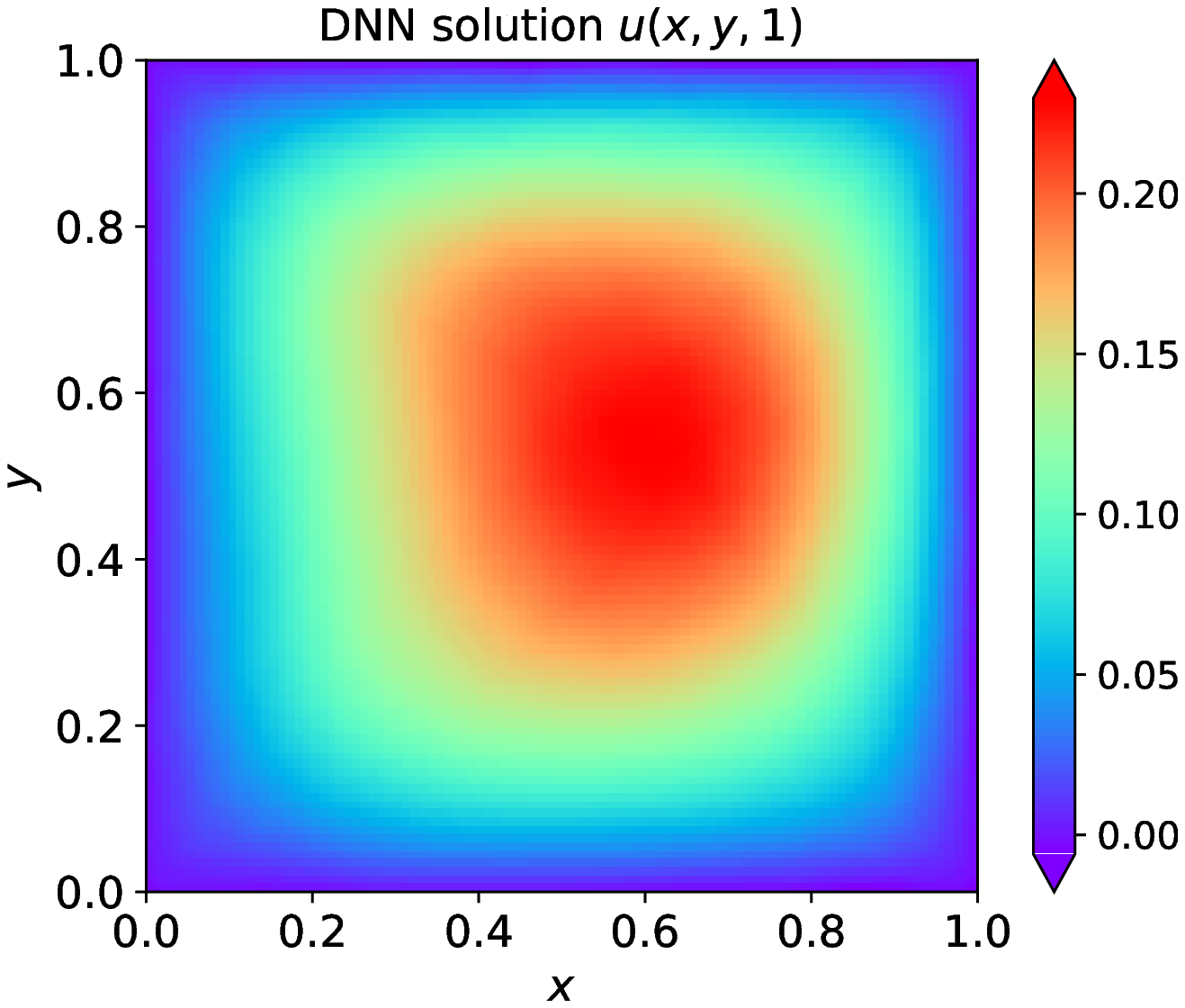}
		\end{minipage}
	}%
	\subfigure[Point-wise errors at $t=1$]{
		\begin{minipage}[t]{0.3\linewidth}
			\centering
			\includegraphics[width=2in]{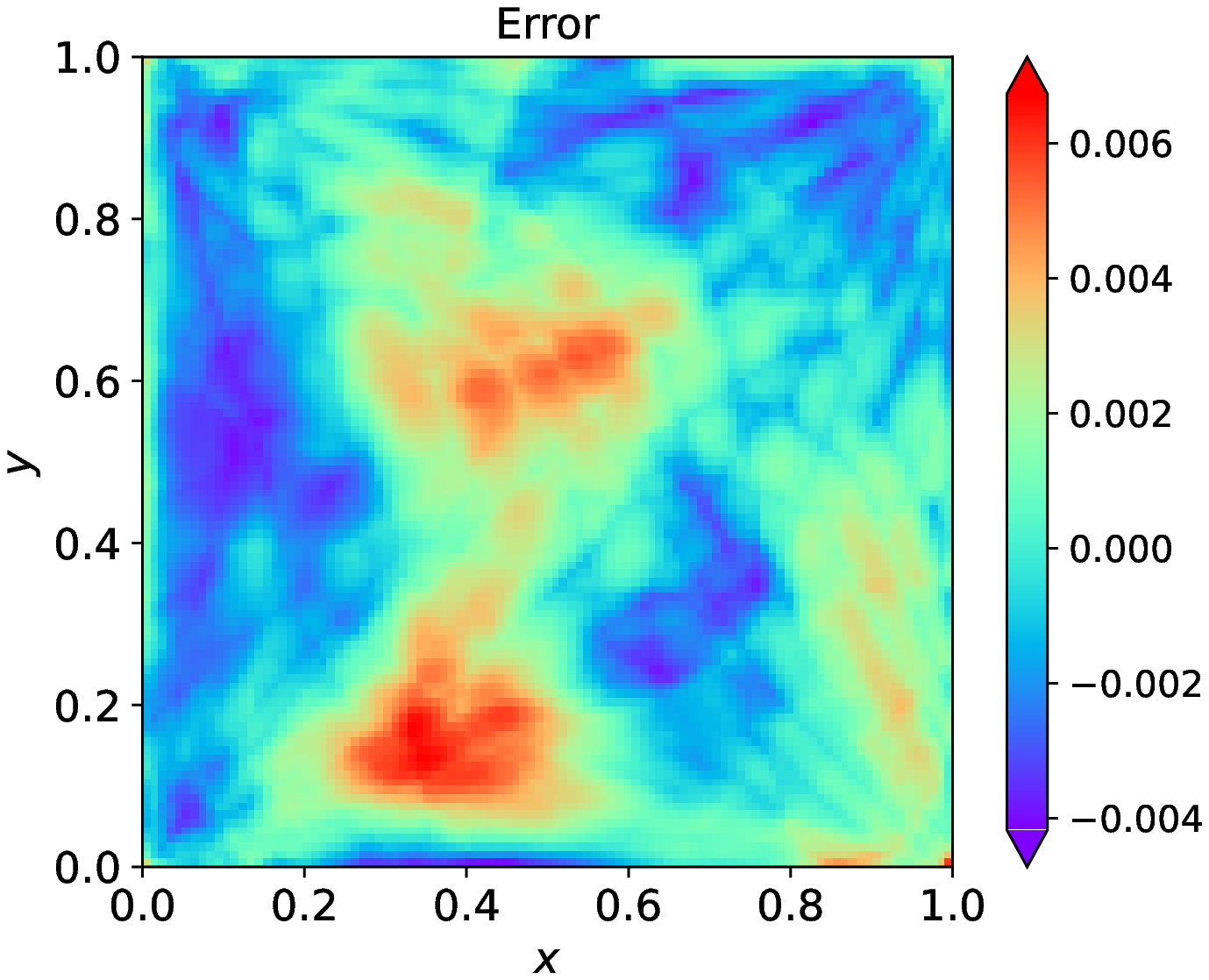}
		\end{minipage}
	}%
	\centering
	\caption{Test sample 2 of the task 1. Comparison between the finite difference solution and the outputs of the deep operator network $G_{NN}^{(\alpha,a)}$ on $101\times 101$ grid points corresponding to the fixed reaction coefficient $c(x,y)=-xy-4$, initial value $u_0(x,y)=6\sin(2\pi x)\sin(3\pi y)$, source $f(x,y)=\sin(3\pi x)\sin(\pi y)+6\exp(x^2+y^2)$, and the fractional order $\alpha=0.1$.}
	\label{fig2.1_sample2}
\end{figure}

We select 1000 training samples and 500 test samples  to train  and test the neural network. The training loss reaches a small value after training as shown in Figure \ref{fig2.1_loss}. We also calculate the relative $l_2$ error between the outputs of the deep operator network $G_{NN}^{(\alpha,a)}$ and the exact solutions on the test dataset, and its average is 0.005912, which means that the trained deep operator network $G_{NN}^{(\alpha,a)}$ can accurately match the numerical solution corresponding to the diffusion coefficient on the test dataset. To visually present the accuracy of outputs of $G_{NN}^{(\alpha,a)}$, two different input samples 
$a$ are shown in Figure \ref{fig2.1_sample1_diffu} and Figure \ref{fig2.1_sample2_diffu}, and its corresponding outputs of the deep operator network as well as the exact solution (finite difference solution) on the test dataset are presented in Figure \ref{fig2.1_sample1} and Figure \ref{fig2.1_sample2}, respectively. It again shows that the prediction of the deep operator network $G_{NN}^{(\alpha,a)}$ is well agreement with the ground truth.

\subsection{Task 2. Learning solution operator mapping $(a,f)$ to solution $u(x,y,t)$.}
In this section, we will learn another solution operator that maps any diffusion coefficients $a(x,y)$ and sources $f(x,y)$ to the solution (\ref{3.1}),
i.e., 
\begin{equation}\label{a_f_t1}
	G: (a, f) \to u(x,y,t),~(x,y,t)\in \bar{\Omega}\times [0,T],
\end{equation}
where we denote the approximate deep operator network as $G_{NN}^{(a,f)}$. 
Similarly, in this part, to get training and testing dataset, we fixed
$\alpha=0.5$, $c(x,y)=0$, $u_{0}(x,y)=\sin(\pi x)\sin(\pi y)$. The source $f$ is generated according to $f\sim \mu$ where $\mu=\mathcal{N}(0,(-\Delta)^{-s})$ with the homogeneous Dirichlet boundary conditions on the operator $-\Delta$, $\mathcal{N}(0,(-\Delta)^{-s})$ denotes the Gaussian measure with mean function 0 and covariance operator $(-\Delta)^{-s}$. In practice, the random function $f$ given by the
Karhunen-Lo${\rm \grave{e}}$ve expansion 
\begin{equation*}
    f=\sum_{k=1}^{\infty}\sqrt{\gamma_k}\zeta_k\phi_k,
\end{equation*}
where the $\{\gamma_k,\phi_k\}_{k=1}^{\infty}$ denotes an orthonormal set of eigenvalues/eigenvectors for the operator $(-\Delta)^{-s}$ ordered so that $\gamma_1\ge\gamma_2\ge\cdots$. The $\{\zeta_k\}_{k=1}^{\infty}$ to be an i.i.d. sequence with $\zeta_1\sim N(0,1)$.
The diffusion coefficient $a$, and the numerical solution $u$ 
of the training/testing dataset are obtained with the same settings as the dataset in the previous section. 

\begin{figure} 
	\centering
	\subfigure[Training loss]
	{\epsfxsize 0.49\hsize \epsfbox{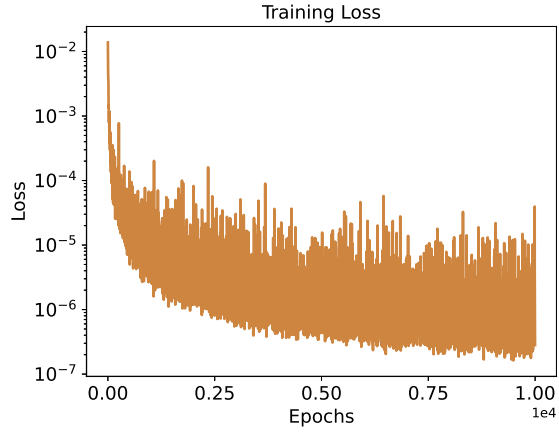}}
	\caption{Training loss of the deep operator network $G_{NN}^{(a,f)}$ for 10000 epochs. The optimizer is Adam with the learning rate is $1\times 10^{-4}$.}\label{fig3.3_loss}
\end{figure}

\begin{figure} 
	\centering
	\subfigure[Diffusion coefficient]{
		\begin{minipage}[t]{0.3\linewidth}
			\centering
			\includegraphics[width=2in]{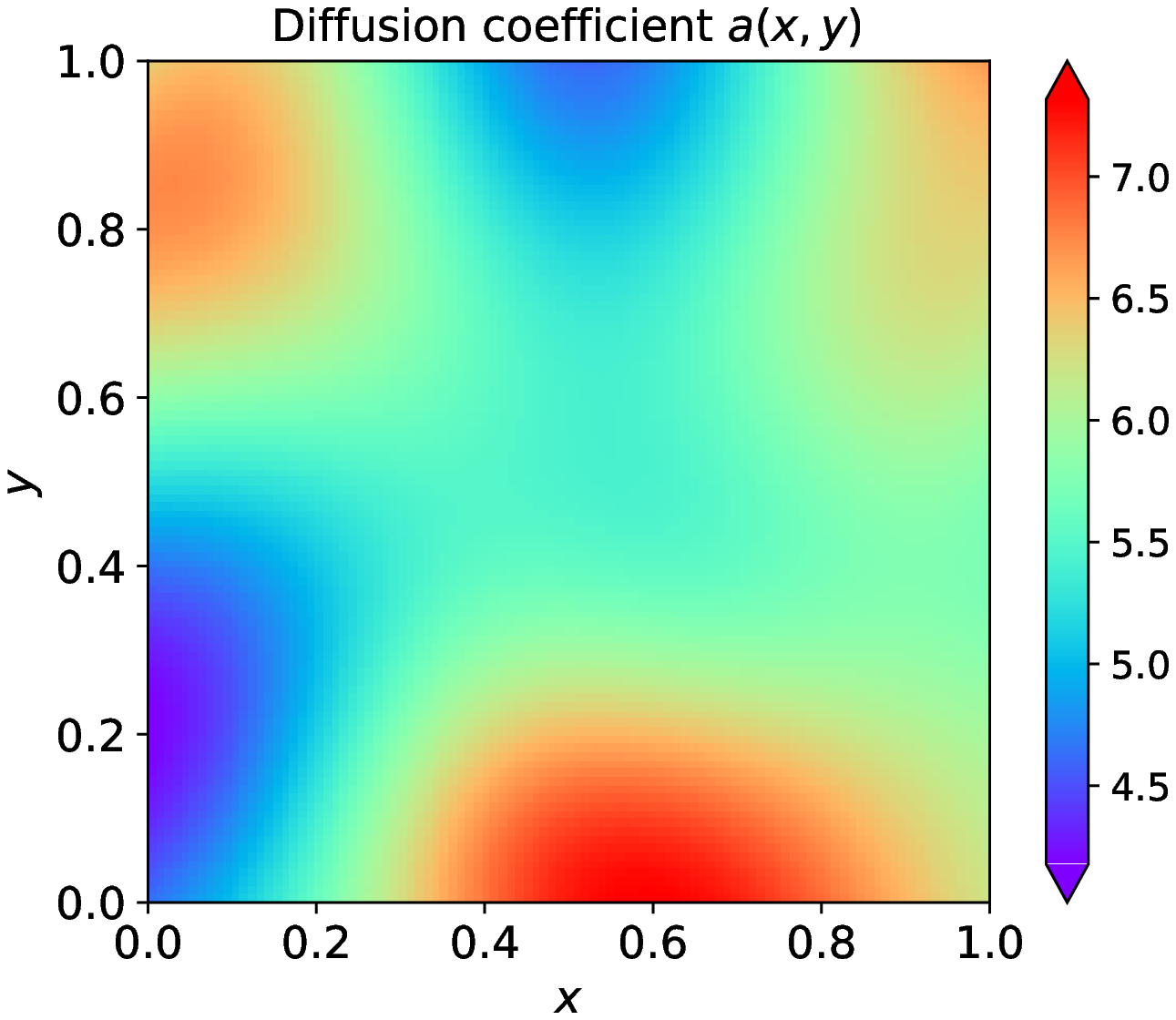}
		\end{minipage}
	}%
	\subfigure[Source]{
		\begin{minipage}[t]{0.3\linewidth}
			\centering
			\includegraphics[width=2in]{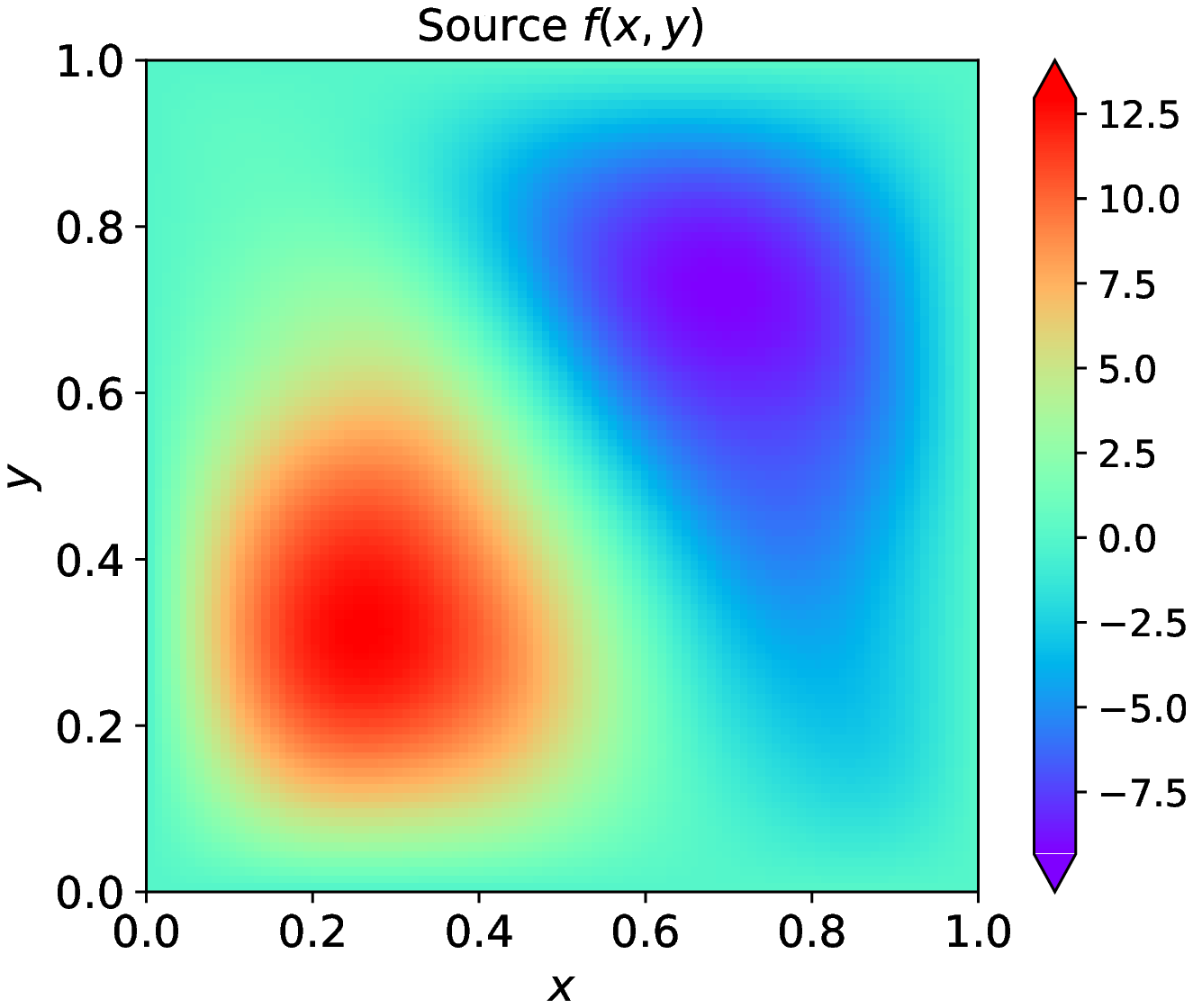}
		\end{minipage}
	}%
		\centering
	\caption{Test sample 1 of the task 2. The input diffusion coefficient $a(x,y)$ and source $f(x,y)$ of the operator network $G_{NN}^{(a,f)}$.}
	\label{fig2.3_sample1_a_f}
\end{figure}

\begin{figure} 
	\centering

	\subfigure[FDM solution at $t=0$]{
		\begin{minipage}[t]{0.3\linewidth}
			\centering
			\includegraphics[width=2in]{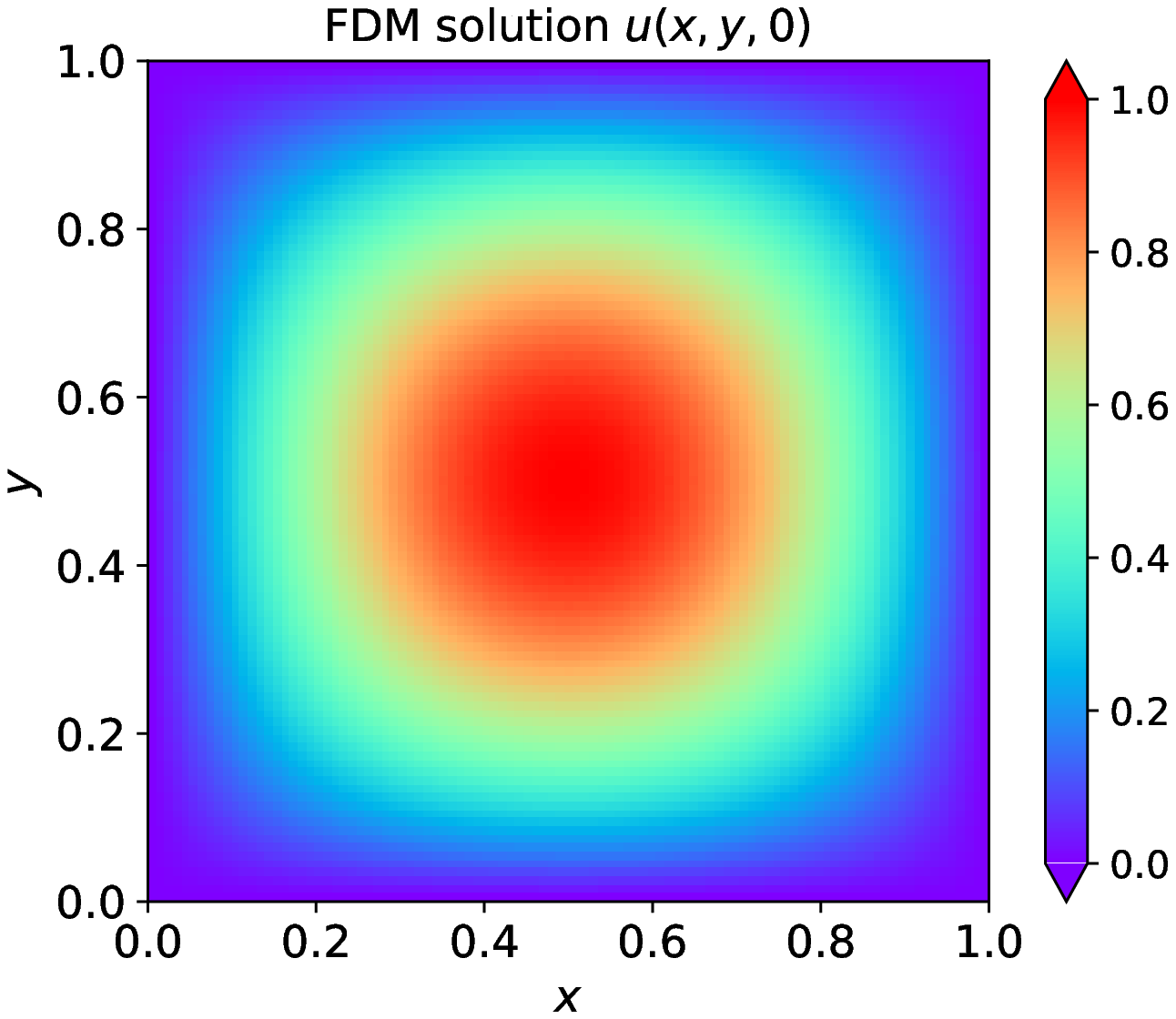}
		\end{minipage}
	}%
	\subfigure[$G_{NN}^{(a,f)}$ solution at $t=0$]{
		\begin{minipage}[t]{0.3\linewidth}
			\centering
			\includegraphics[width=2in]{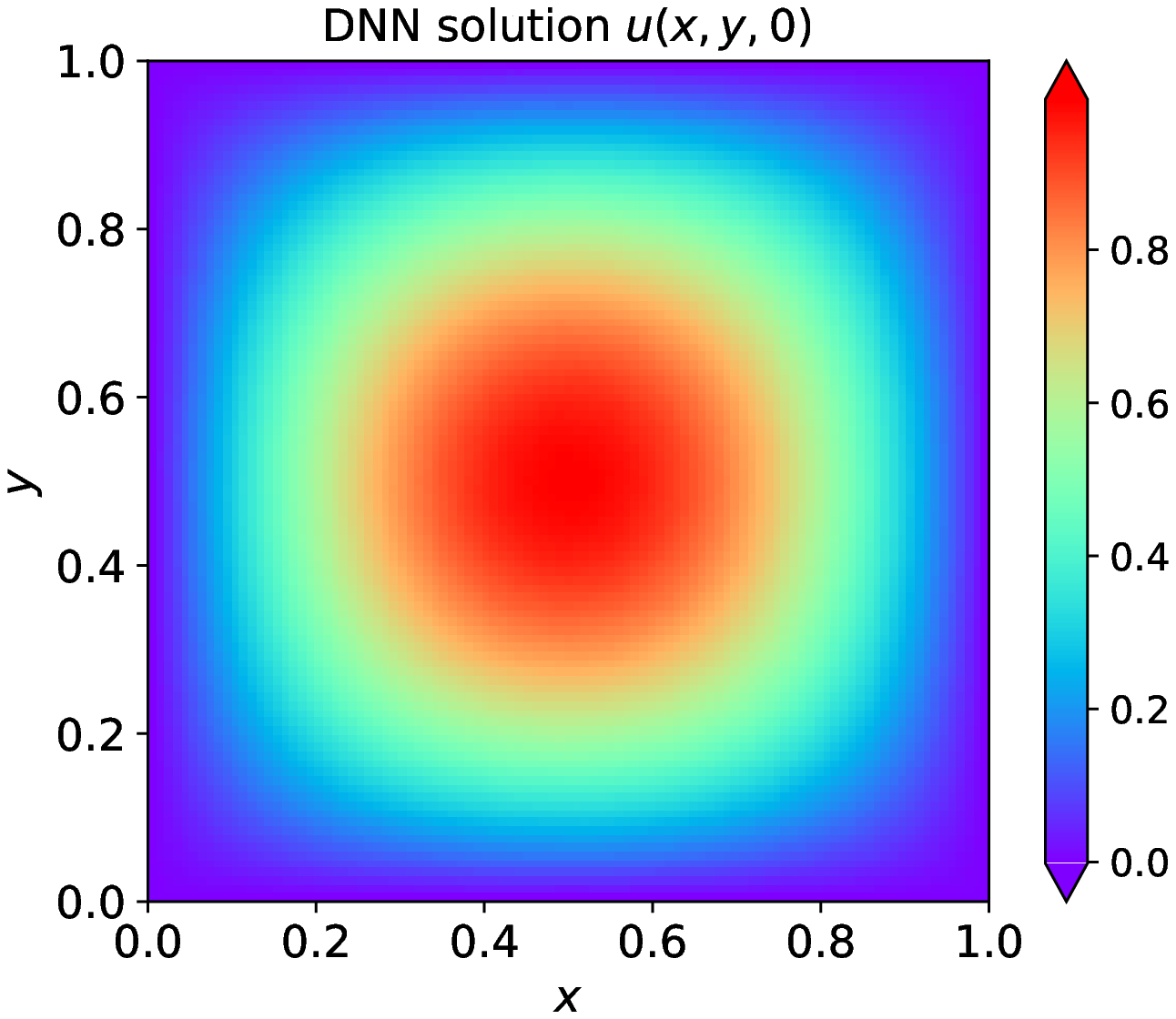}
		\end{minipage}
	}%
	\subfigure[Point-wise errors at $t=0$]{
		\begin{minipage}[t]{0.3\linewidth}
			\centering
			\includegraphics[width=2in]{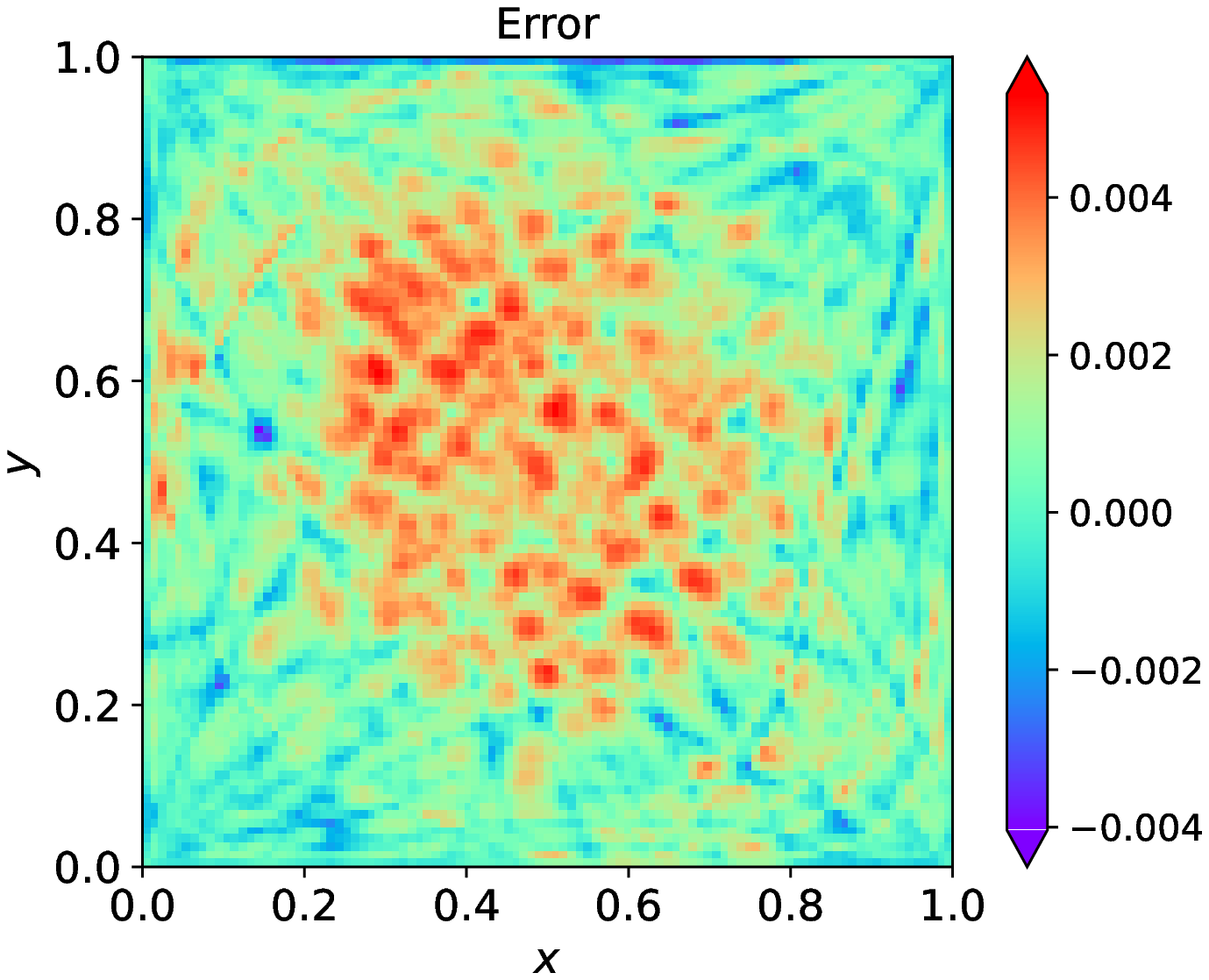}
		\end{minipage}
	}%

	\subfigure[FDM solution at $t=0.04$]{
		\begin{minipage}[t]{0.3\linewidth}
			\centering
			\includegraphics[width=2in]{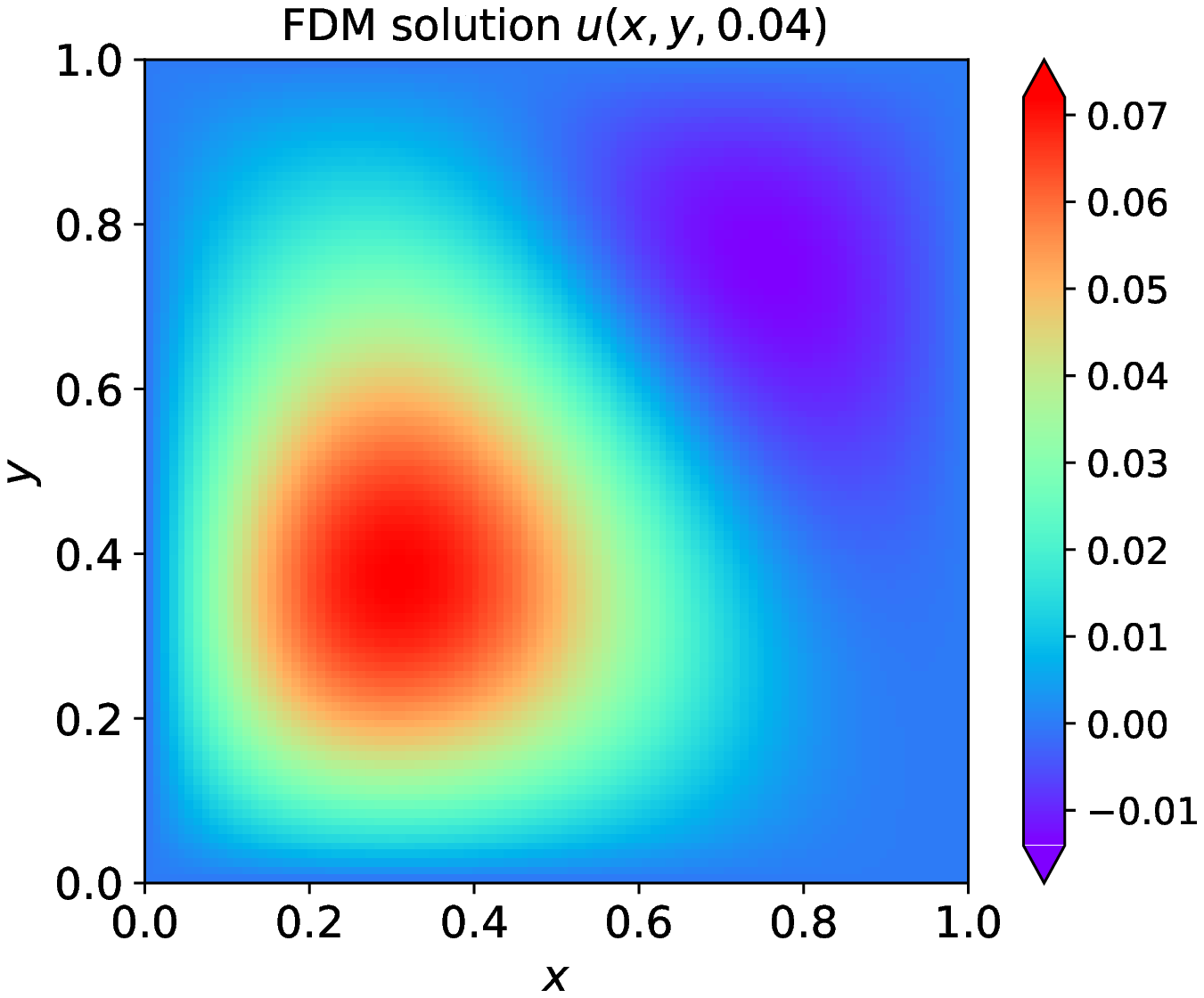}
		\end{minipage}
	}%
	\subfigure[$G_{NN}^{(a,f)}$ solution at $t=0.04$]{
		\begin{minipage}[t]{0.3\linewidth}
			\centering
			\includegraphics[width=2in]{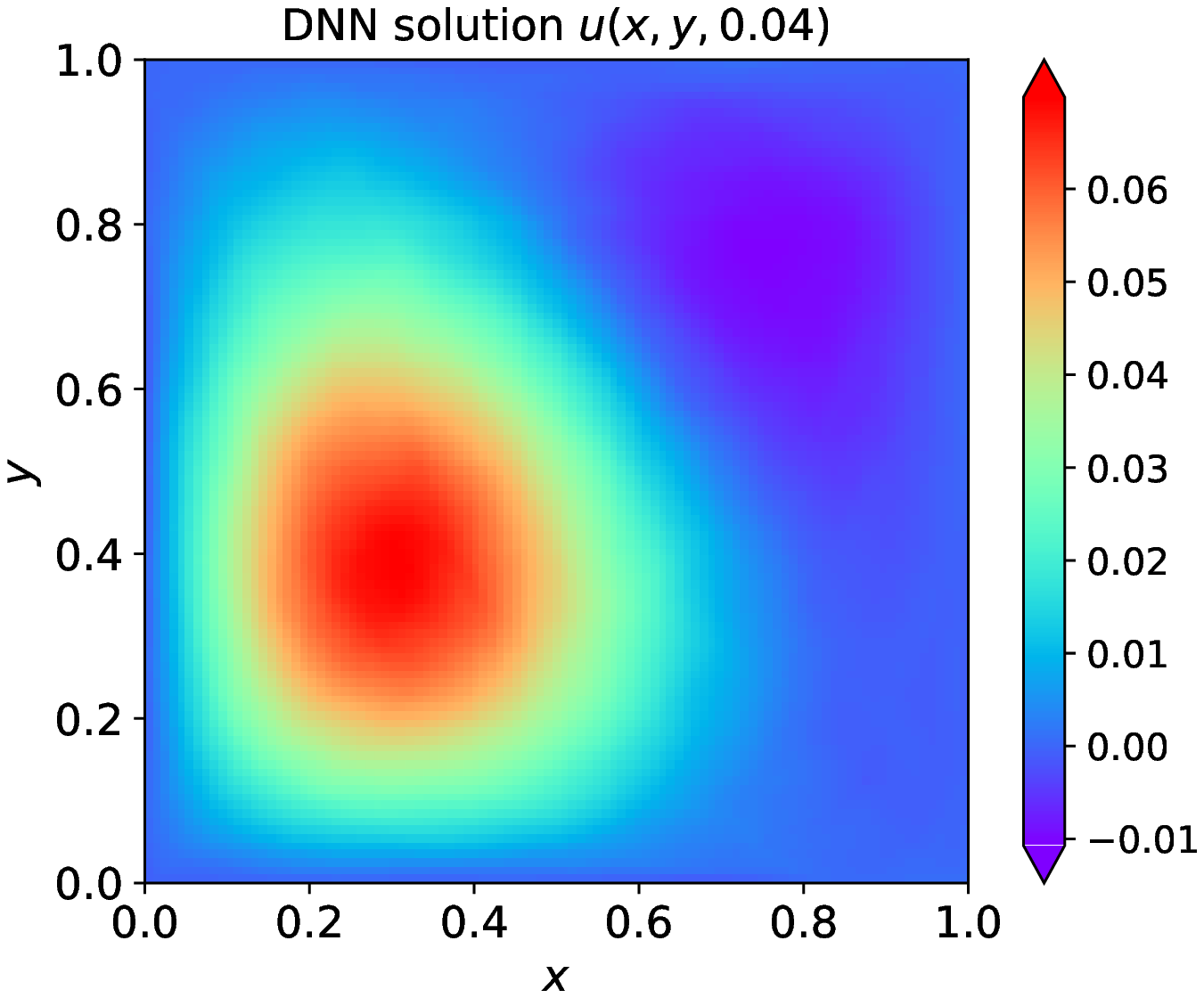}
		\end{minipage}
	}%
	\subfigure[Point-wise errors at $t=0.04$]{
		\begin{minipage}[t]{0.3\linewidth}
			\centering
			\includegraphics[width=2in]{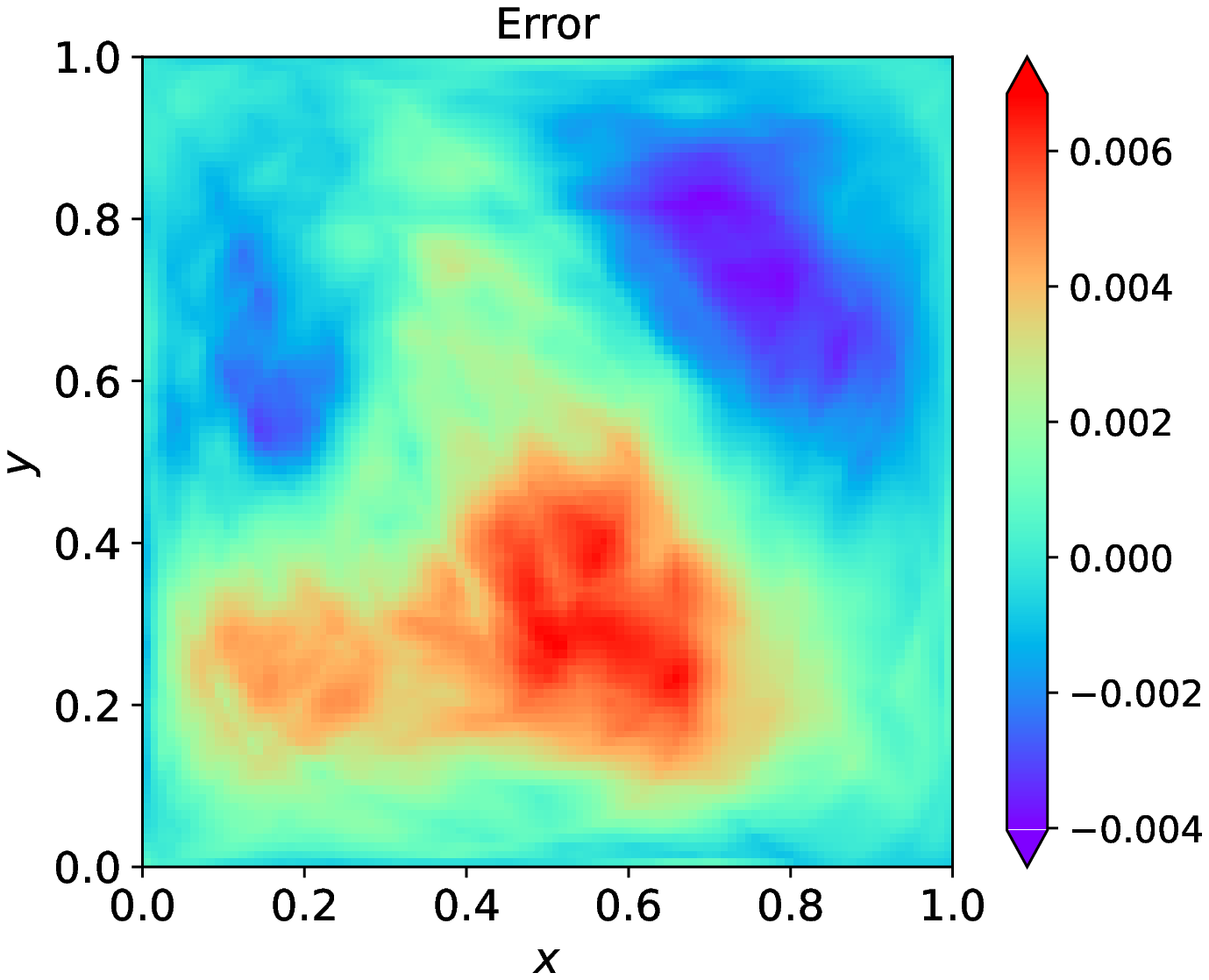}
		\end{minipage}
	}%

	\subfigure[FDM solution at $t=1$]{
		\begin{minipage}[t]{0.3\linewidth}
			\centering
			\includegraphics[width=2in]{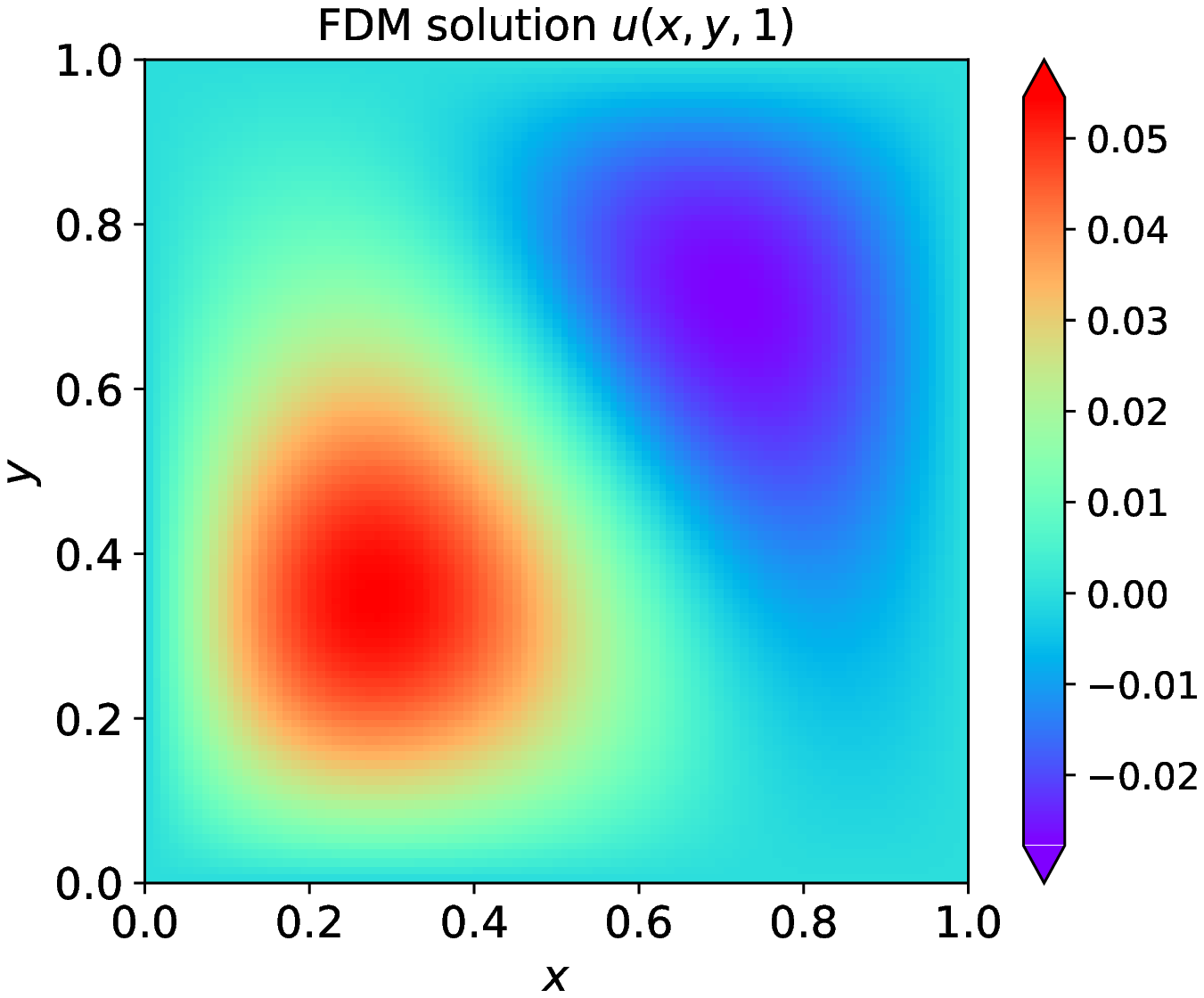}
		\end{minipage}
	}%
	\subfigure[$G_{NN}^{(a,f)}$ solution at $t=1$]{
		\begin{minipage}[t]{0.3\linewidth}
			\centering
			\includegraphics[width=2in]{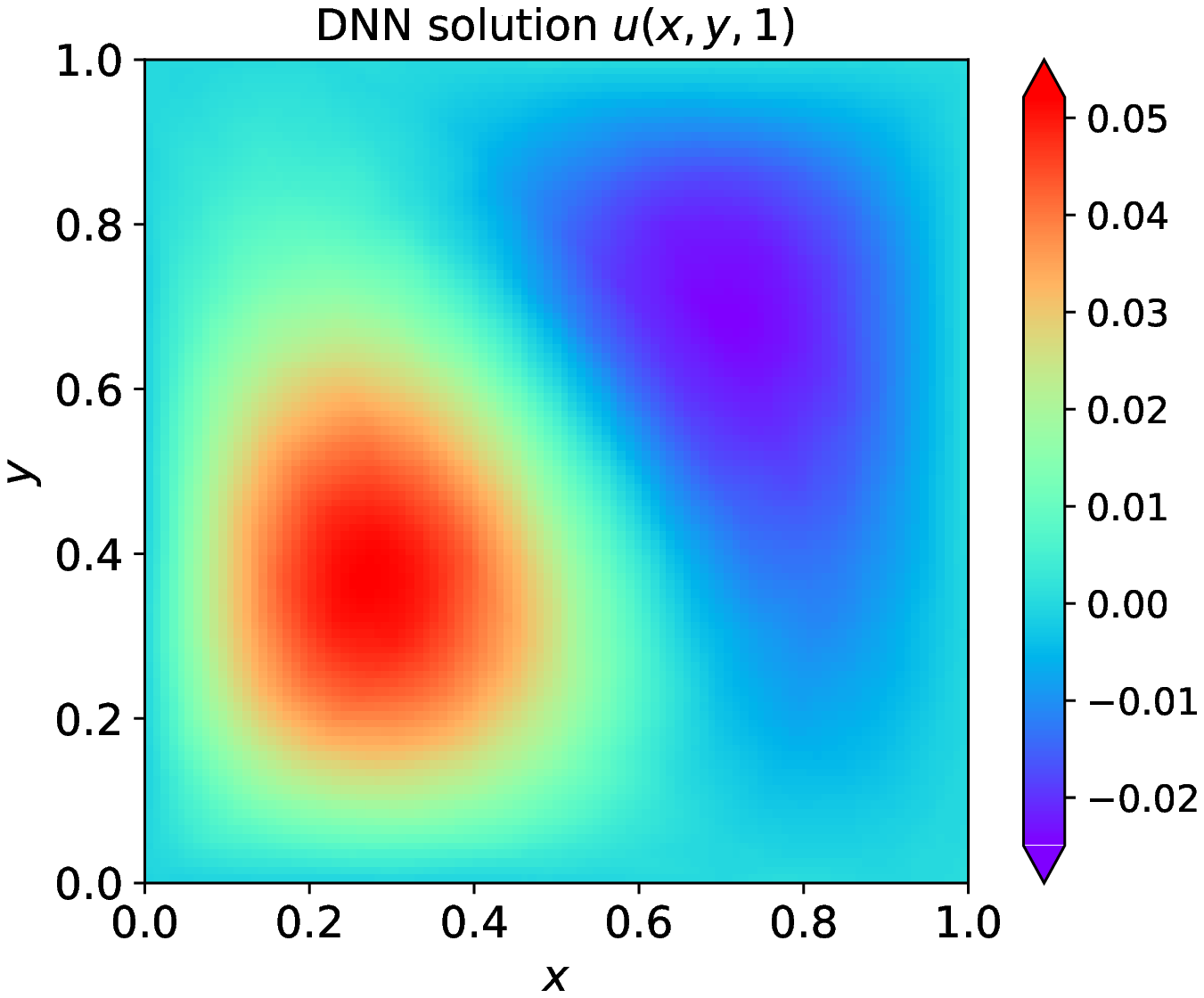}
		\end{minipage}
	}%
	\subfigure[Point-wise errors at $t=1$]{
		\begin{minipage}[t]{0.3\linewidth}
			\centering
			\includegraphics[width=2in]{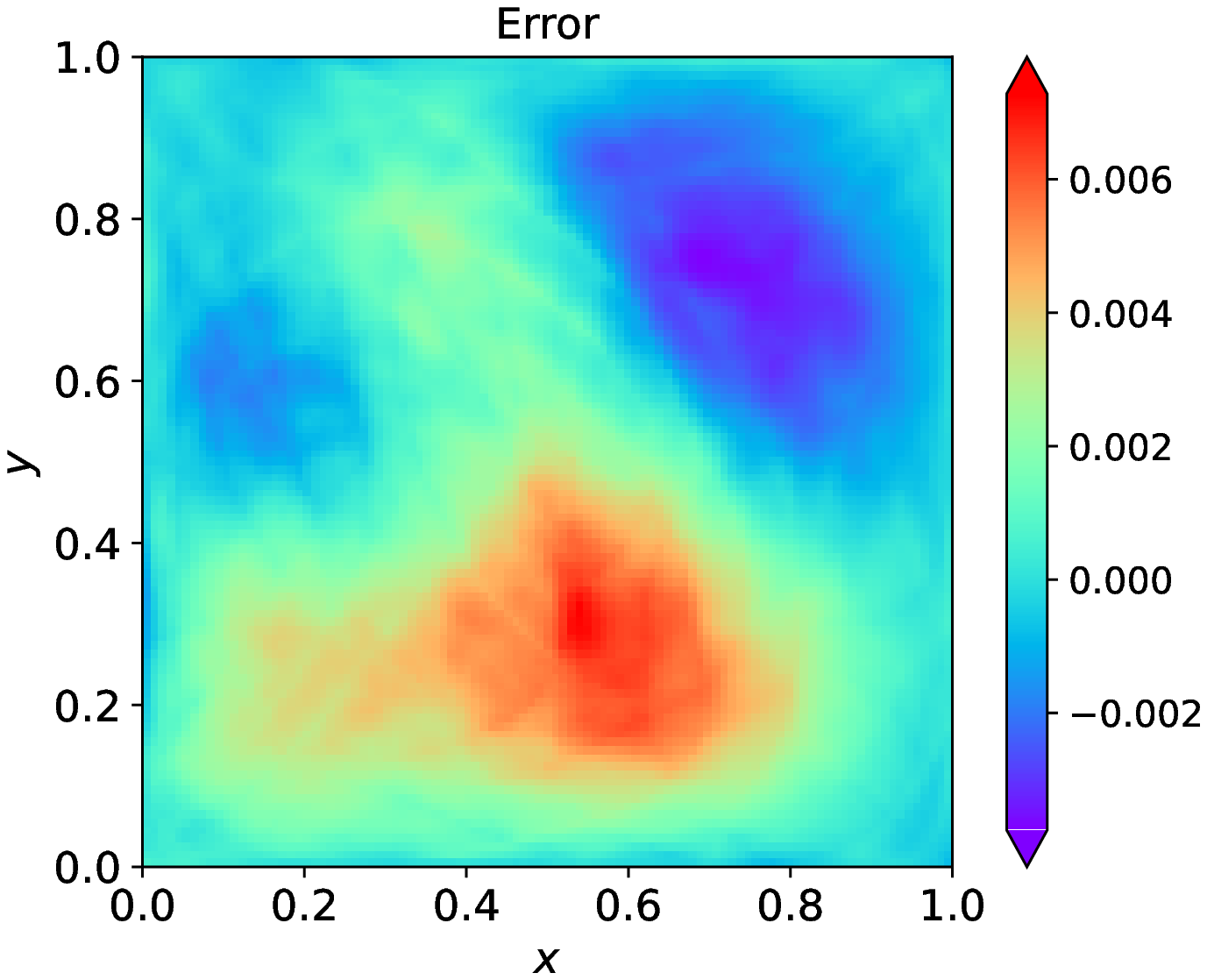}
		\end{minipage}
	}%
	\centering
	\caption{Test sample 1 of the task 2. Comparison between the finite difference solution and the outputs of deep operator network $G_{NN}^{(a,f)}$ on $101\times 101$ grids points corresponding to the fixed fractional order $\alpha=0.5$, reaction coefficient $c(x,y)=0$, and initial value $u_0(x,y)=\sin(\pi x)\sin(\pi y)$.}
	\label{fig2.3_sample1}
\end{figure}

\begin{figure} 
	\centering
	\subfigure[Diffusion coefficient]{
		\begin{minipage}[t]{0.3\linewidth}
			\centering
			\includegraphics[width=2in]{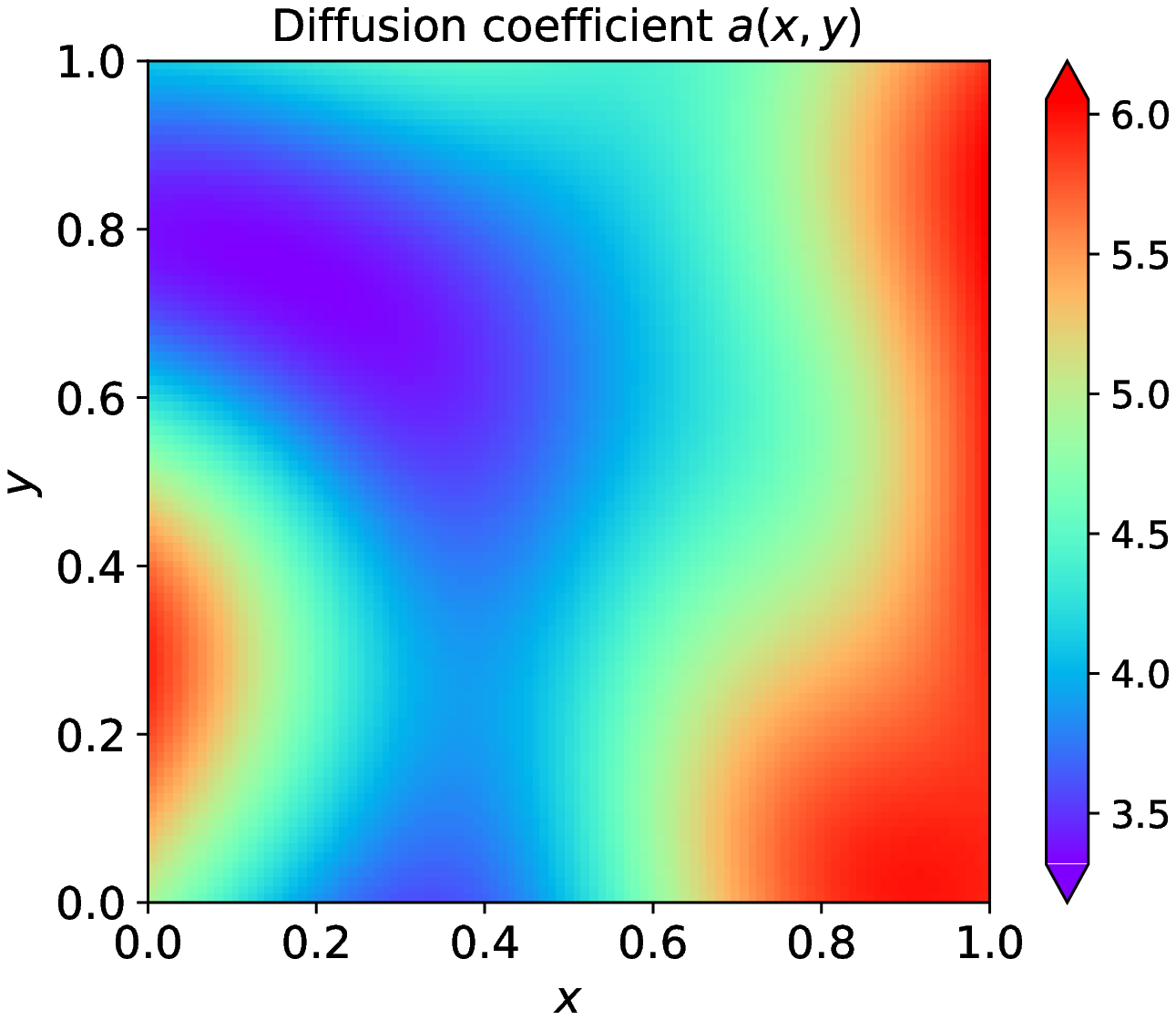}
		\end{minipage}
	}%
	\subfigure[Source]{
		\begin{minipage}[t]{0.3\linewidth}
			\centering
			\includegraphics[width=2in]{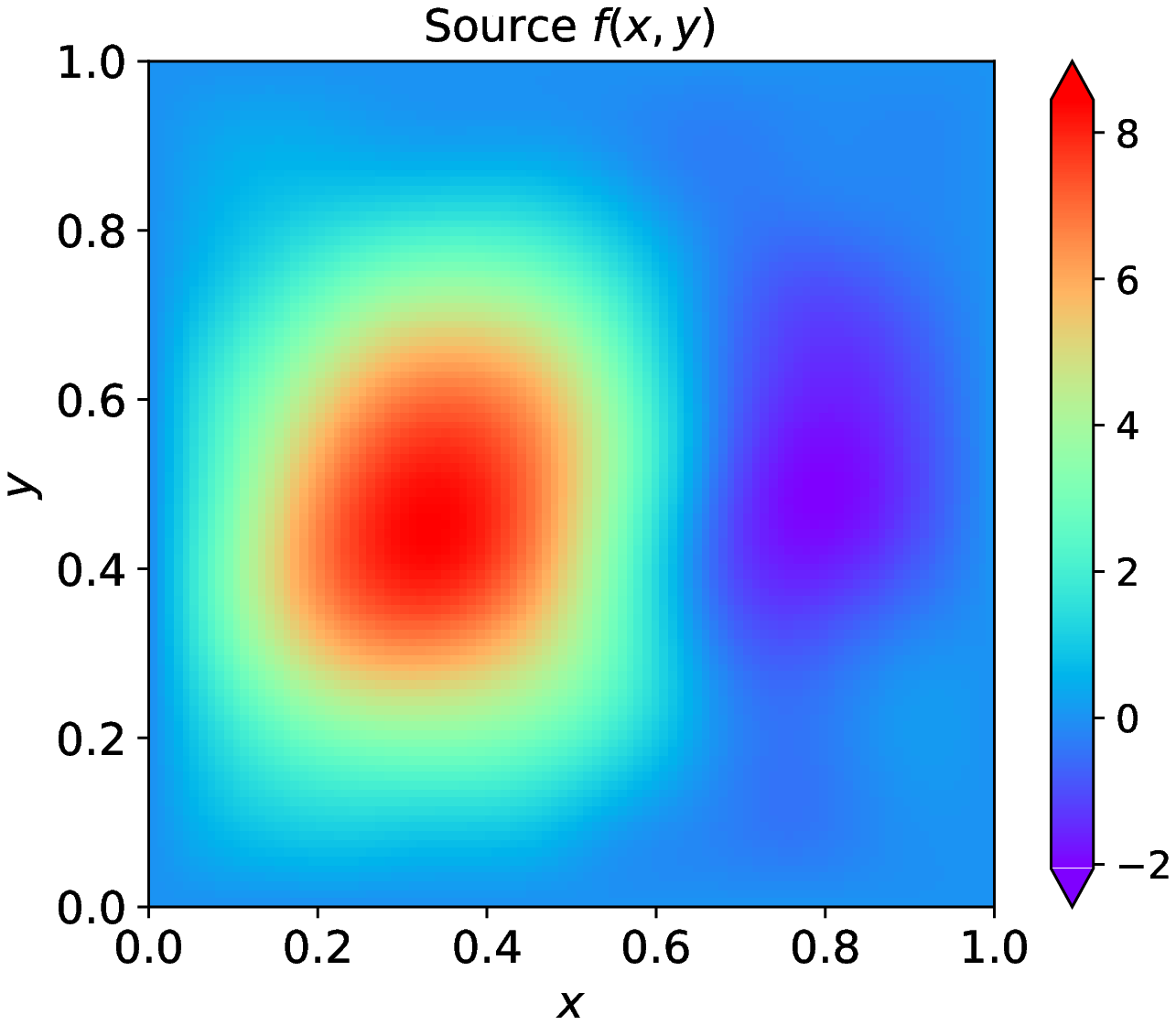}
		\end{minipage}
	}%
		\centering
	\caption{Test sample 2 of the task 2. The input diffusion coefficient $a(x,y)$ and source $f(x,y)$ of the operator network $G_{NN}^{(a,f)}$.}
	\label{fig2.3_sample2_a_f}
\end{figure}

\begin{figure} 
	\centering
	\subfigure[FDM solution at $t=0$]{
		\begin{minipage}[t]{0.3\linewidth}
			\centering
			\includegraphics[width=2in]{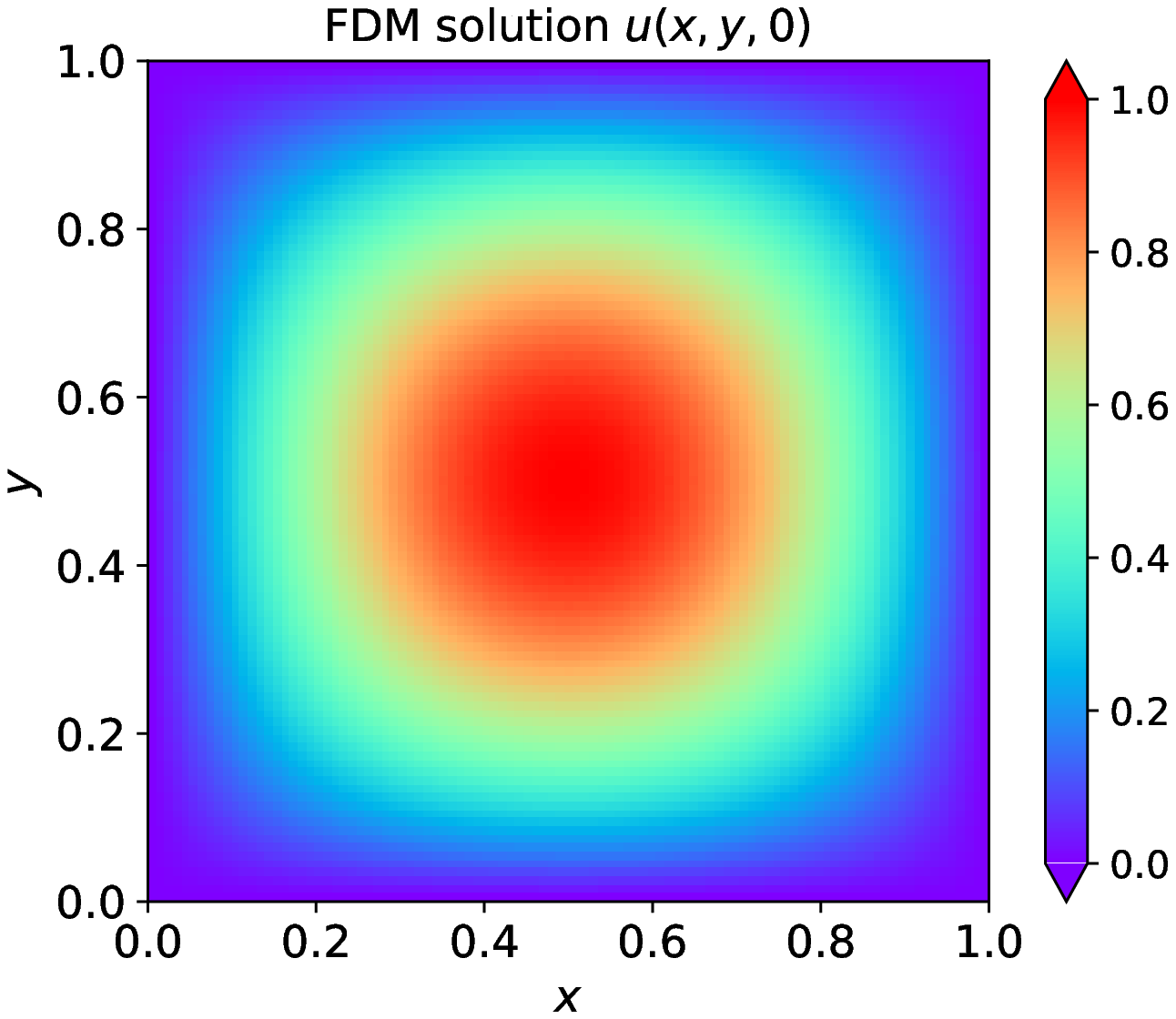}
		\end{minipage}
	}%
	\subfigure[$G_{NN}^{(a,f)}$ solution at $t=0$]{
		\begin{minipage}[t]{0.3\linewidth}
			\centering
			\includegraphics[width=2in]{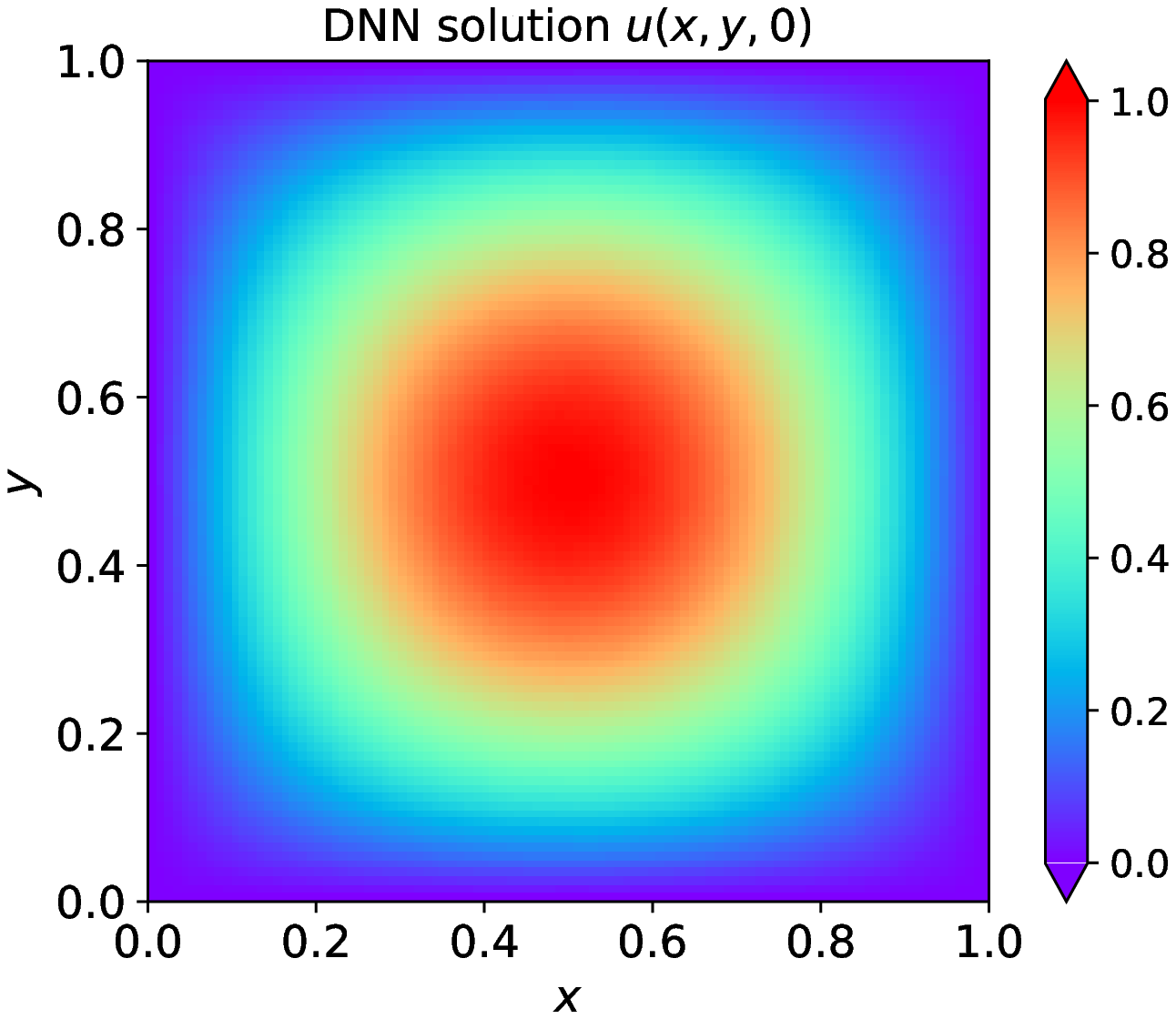}
		\end{minipage}
	}%
	\subfigure[Point-wise errors at $t=0$]{
		\begin{minipage}[t]{0.3\linewidth}
			\centering
			\includegraphics[width=2in]{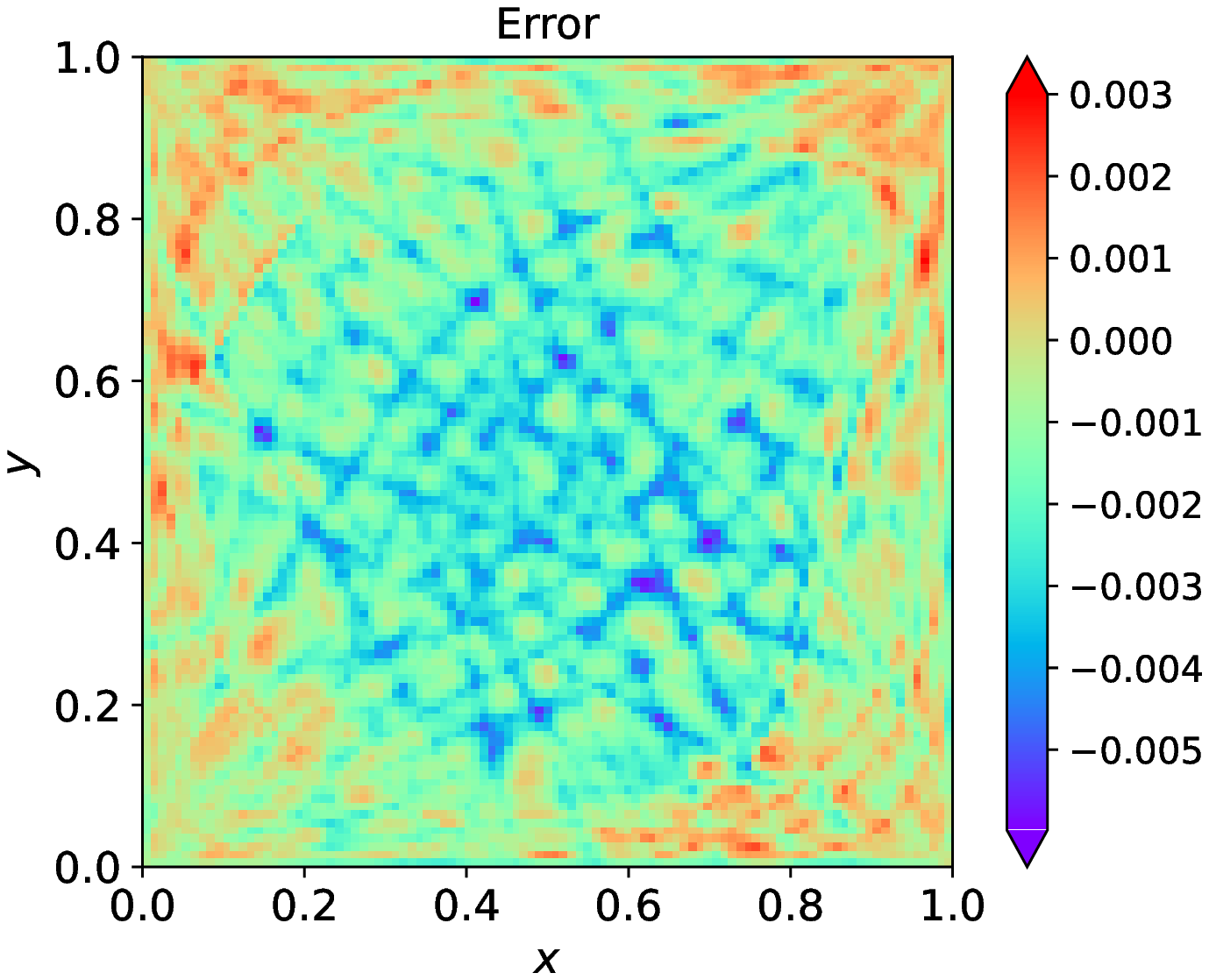}
		\end{minipage}
	}%

	\subfigure[FDM solution at $t=0.04$]{
		\begin{minipage}[t]{0.3\linewidth}
			\centering
			\includegraphics[width=2in]{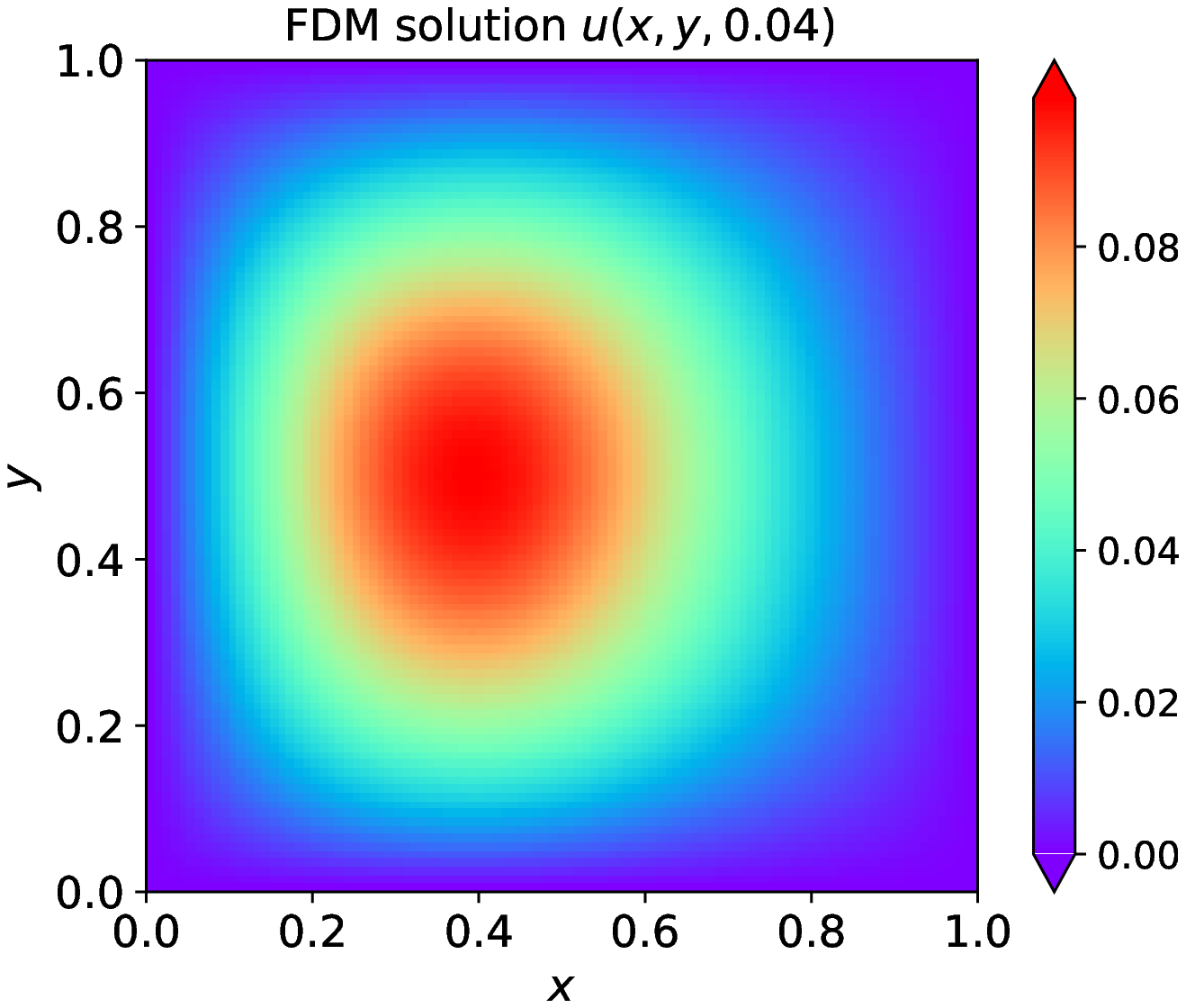}
		\end{minipage}
	}%
	\subfigure[$G_{NN}^{(a,f)}$ solution at $t=0.04$]{
		\begin{minipage}[t]{0.3\linewidth}
			\centering
			\includegraphics[width=2in]{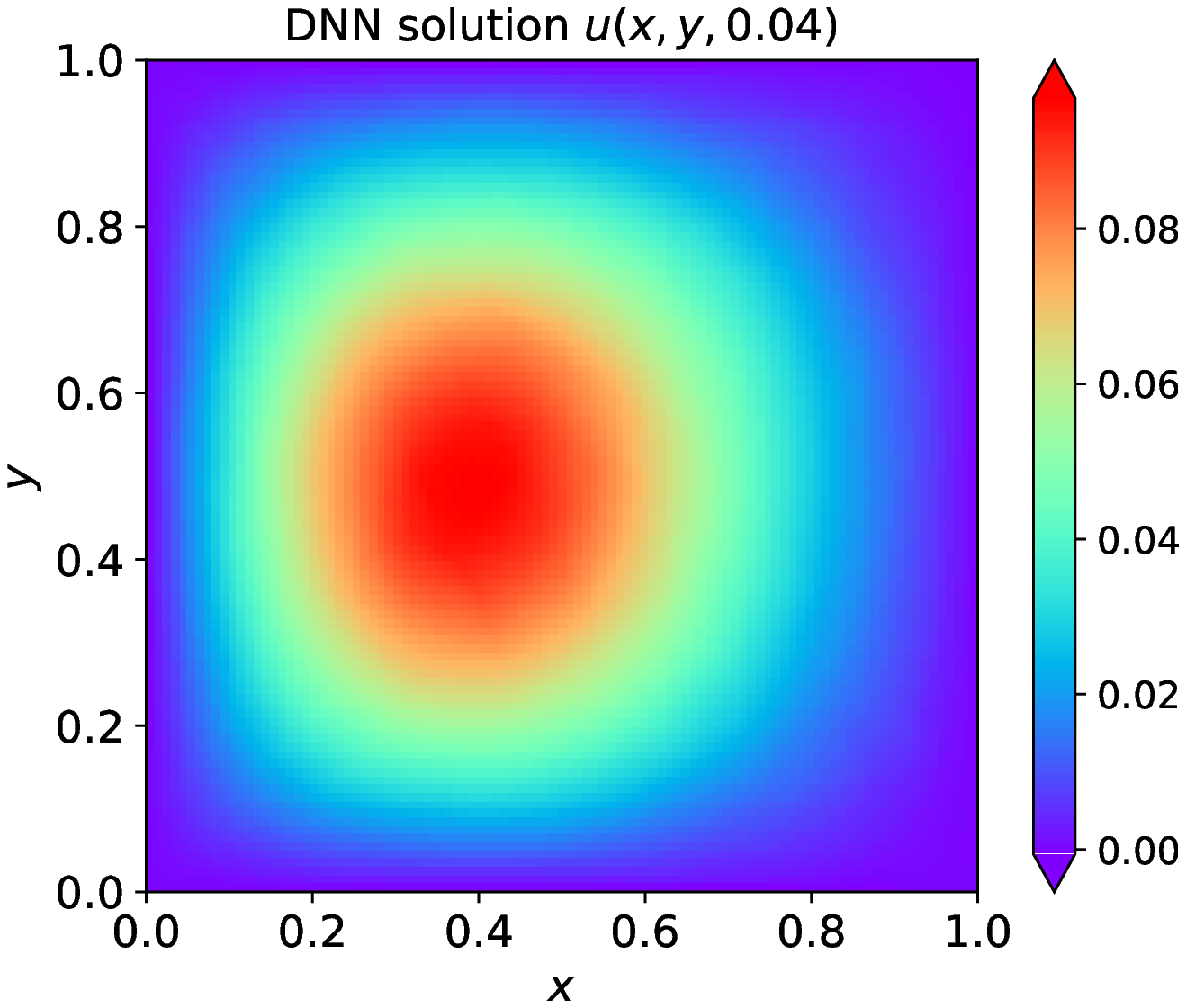}
		\end{minipage}
	}%
	\subfigure[Point-wise errors at $t=0.04$]{
		\begin{minipage}[t]{0.3\linewidth}
			\centering
			\includegraphics[width=2in]{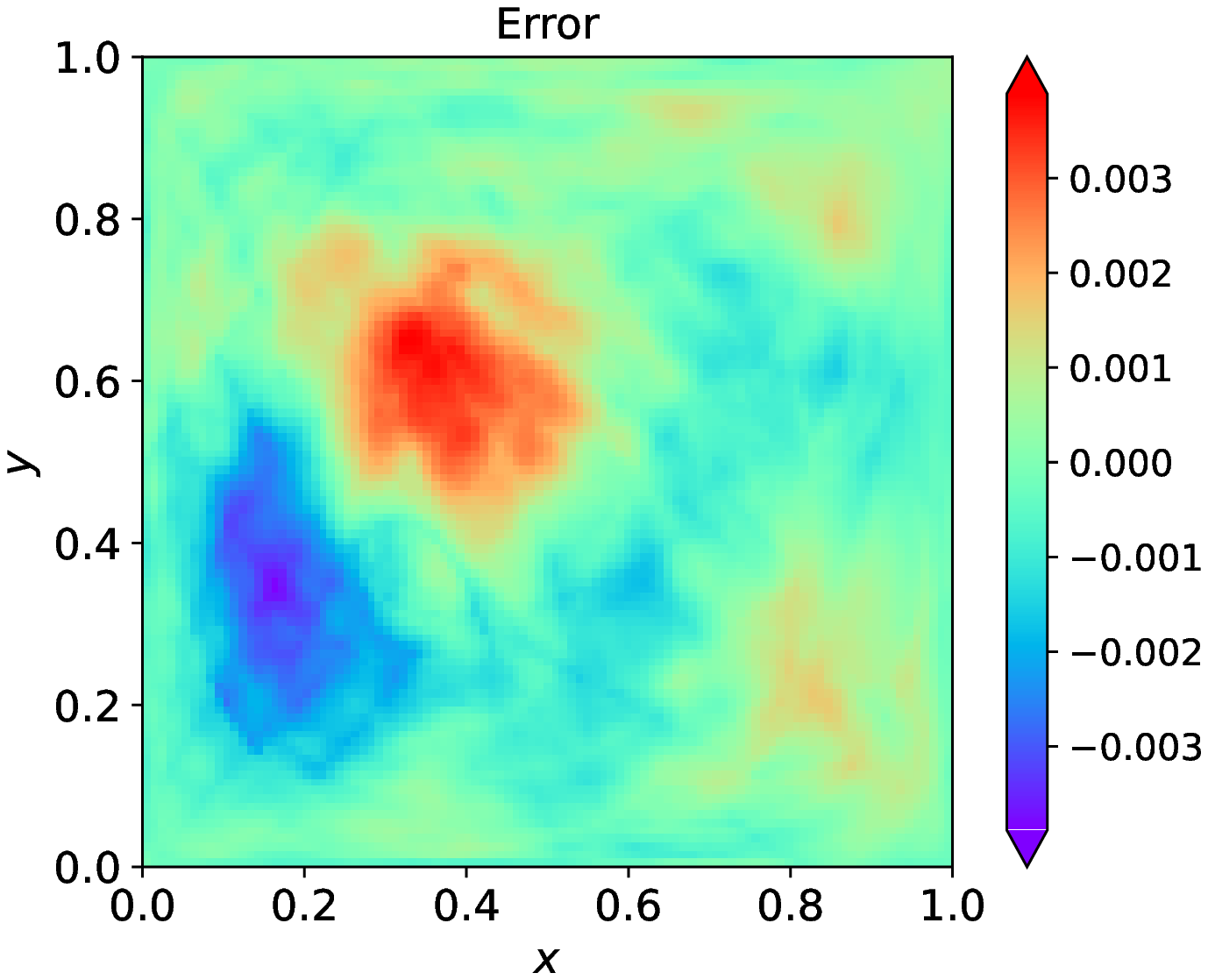}
		\end{minipage}
	}%

	\subfigure[FDM solution at $t=1$]{
		\begin{minipage}[t]{0.3\linewidth}
			\centering
			\includegraphics[width=2in]{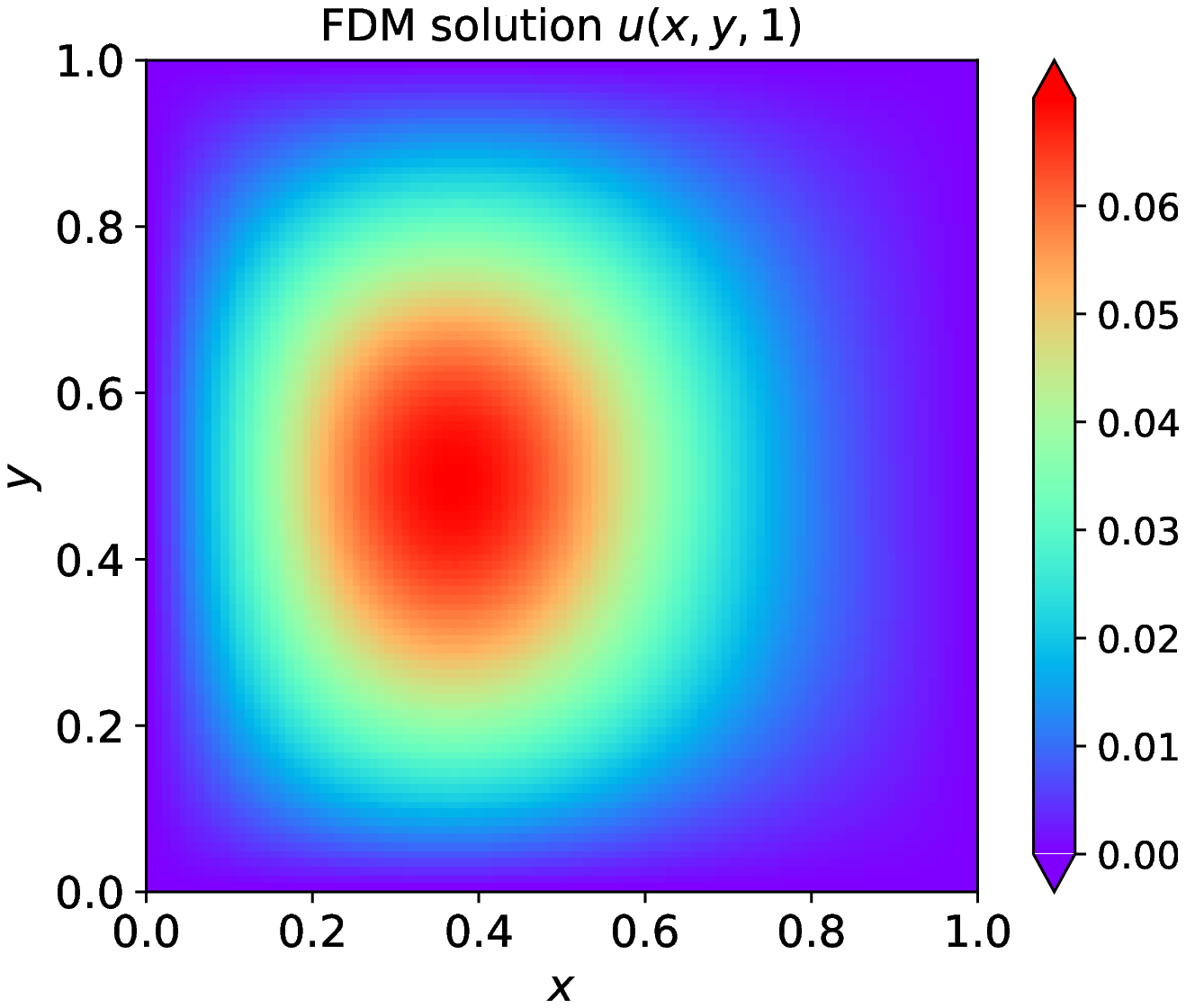}
		\end{minipage}
	}%
	\subfigure[$G_{NN}^{(a,f)}$ solution at $t=1$]{
		\begin{minipage}[t]{0.3\linewidth}
			\centering
			\includegraphics[width=2in]{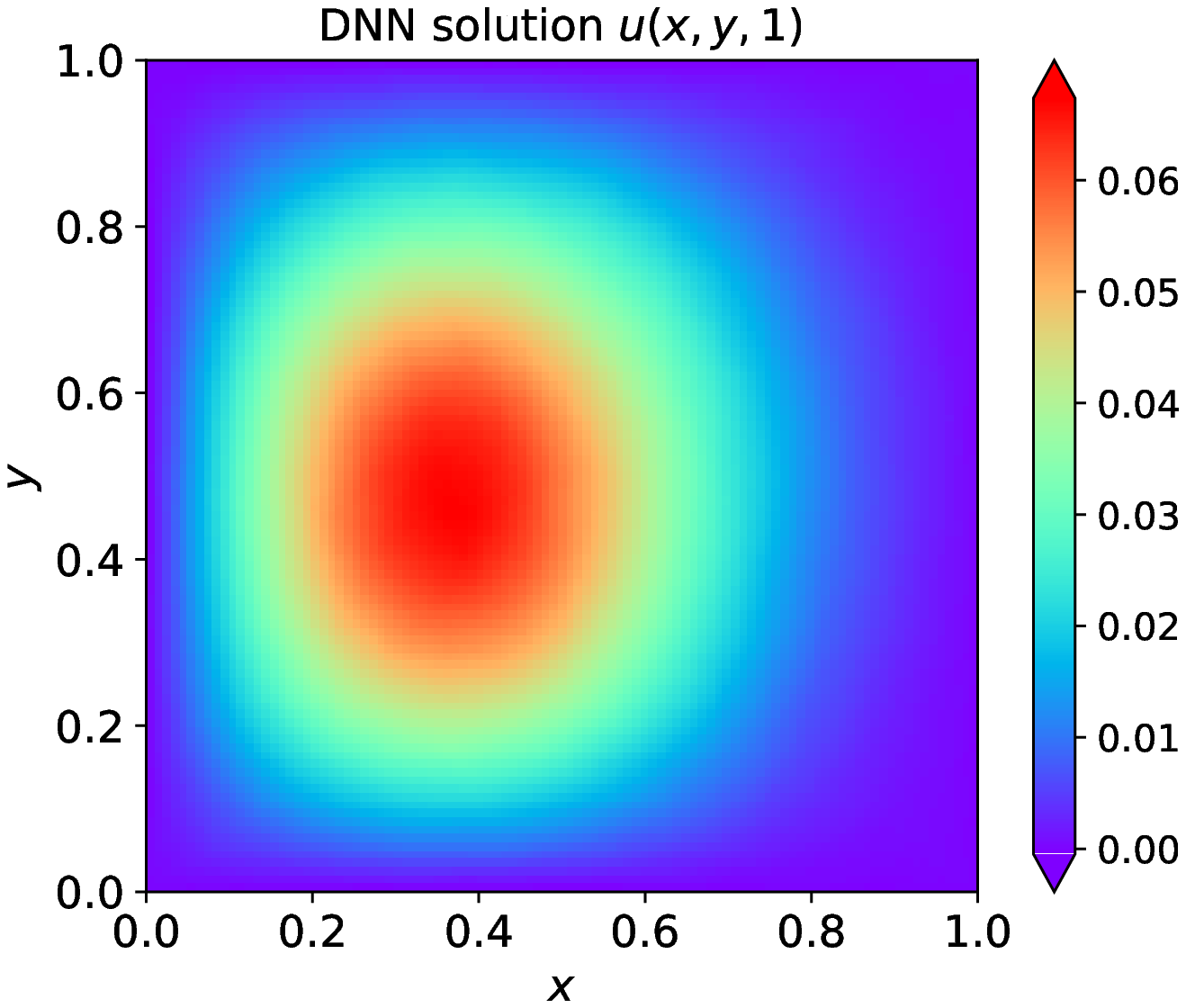}
		\end{minipage}
	}%
	\subfigure[Point-wise errors at $t=1$]{
		\begin{minipage}[t]{0.3\linewidth}
			\centering
			\includegraphics[width=2in]{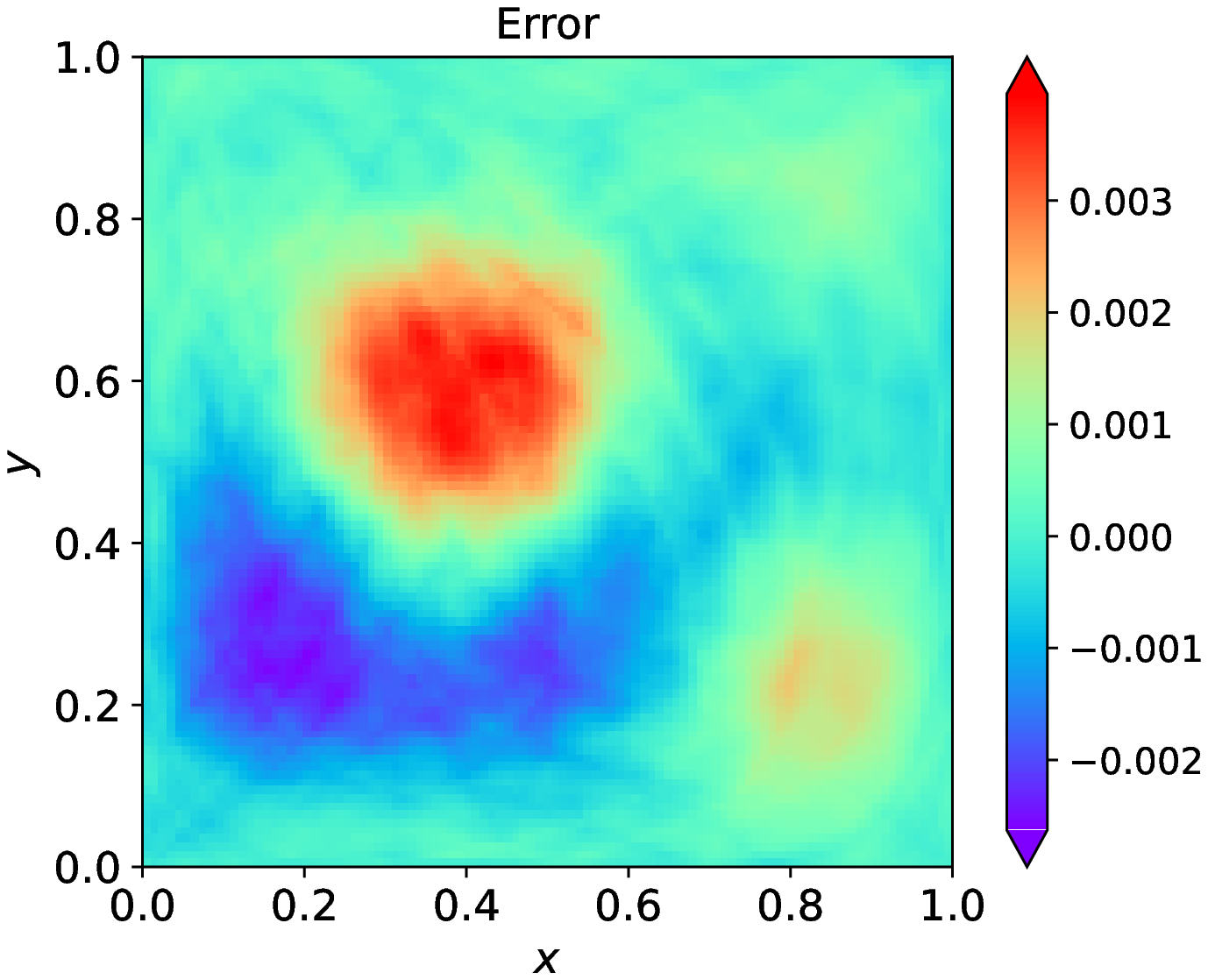}
		\end{minipage}
	}%
	\centering
	\caption{Test sample 2 of the task 2. Comparison between the finite difference solution and the outputs of deep operator network $G_{NN}^{(a,f)}$ on $101\times 101$ grids points corresponding to the fixed fractional order $\alpha=0.5$, reaction coefficient $c(x,y)=0$, and initial value $u_0(x,y)=\sin(\pi x)\sin(\pi y)$.}
	\label{fig2.3_sample2}
\end{figure}

We take 1000 and 500 samples for aims of training and testing, respectively. Figure \ref{fig3.3_loss} shows that the training loss reaches a small value after training. To assess the accuracy of the approximate operator network $G_{NN}^{(a,f)}$ in the subdiffusion problem (\ref{1.1}), two different input samples $(a,f)$ of the test dataset are shown in Figure \ref{fig2.3_sample1_a_f} and Figure \ref{fig2.3_sample2_a_f} and the predicted outputs of the operator network $G_{NN}^{(a,f)}$, the exact numerical solutions, and the errors at different time are shown in Figure \ref{fig2.3_sample1} and Figure \ref{fig2.3_sample2}.
Moreover, we calculate the relative $l_2$ error between the outputs of the deep operator network $G_{NN}^{(a,f)}$ and the exact solution of the test dataset, and its average is 0.022617. These numerical experiments illustrate that the trained deep operator network can be used as an accurate surrogate of solution operator of the subdiffusion problem (\ref{1.1}).

\section{Deep Bayesian inversion in subdiffusion}\label{sec4}
In this section, we will study a few inverse problems that are ill-posed and difficult to solve in the subdiffusion problem (\ref{1.1}). These difficulties stem from two sources. One is that inverse problems are ill-posed. The essence of inverse problems is to solve the inverse of the forward map. Because the inverse operator is unbounded, the noise in the measurement data is amplified when solving an inverse problem, making it difficult to obtain a stable solution. The complexity of the physical model poses the second challenge. To solve the inverse problems, the forward or adjoint problems of the physical model have to be solved repeatedly. At the same time, for the subdiffusion (\ref{1.1}), the fractional derivative includes a convolution, which requires much time in discretizing the integral when we simulating a forward problem. To summarize, in order to solve the inverse problems of the subdiffusion (\ref{1.1}), it is necessary to combine effective regularization methods with fast numerical solution methods.

In this paper, we apply a Bayesian inversion method to solve inverse problems of the subdiffusion problem (\ref{1.1}). 
The following provides a brief introduction to the Bayesian inversion.

First, we need to define a forward map:
\begin{equation}\label{IP_4.0}
	F(m) = u(x,y,t), \quad (x,y,t)\in \Omega\times[0,T],
\end{equation}
where $m$ denotes the unknown parameter, $u(x,y,t)$ represents the solution of the problem (\ref{1.1}). Then the noisy observed data generated by
\begin{equation}
	d = g (m) + \eta,
\end{equation}
here the $g=\mathcal{O}\circ F$, $\mathcal{O}$ is an observation operator, and $\eta$ be the observed noise which is statistically independent with the input $m$
and $\eta\sim \mathcal{N}(0,\Sigma)$, where $\Sigma=\sigma^2 I$ and $\sigma$ denotes the standard deviation of the noise. In the Bayesian framework,  the prior of unknown parameter $m$ is described in terms of a probability density function $p_0(m)$ and the goal is to find
the posterior probability measure $p_{d}(m)$. The probability density function of $d$ given $m$ can be denoted by
\begin{equation*}
	p(d|m)=p(d-g(m)),
\end{equation*}
where    $p$ represents the probability density function of the noise $\eta$.
Combining the Bayes formula, then the posterior probability density function of $p_{d}(m)$ is given via
\begin{equation}\label{IP_4.1}
	p_{d}(m)\propto p(d-g(m))p_{0}(m).
\end{equation}

Obtaining the posterior probability distribution is the first step in solving Bayesian inverse problems, and how to solve it becomes a complicated and crucial issue. There are many methods available today to solve the Bayesian posterior distribution, including the Kalman filter, Markov chain Monte Carlo, and others. However, these methods need to solve a large number of forward problems repeatedly. At the same time, traditional numerical methods such as the finite element method and the finite difference method are time-consuming, which makes solving the Bayesian inverse problems extremely difficult. Therefore, it is necessary to find a fast and effective numerical method to accelerate the Bayesian inversion. In this paper, we use the trained deep operator network as an approximation operator of the forward operator $F$ to accelerate the Bayesian inversion. This allows us to solve the inverse problems in the subdiffusion problem (\ref{1.1}) more efficiently.

\subsection{Inverse the order of fractional derivative}
The rate of particle diffusion is influenced by the fractional order $\alpha$ in the time-fractional subdiffusion model (\ref{1.1}). In many practical applications, however, we cannot know the exact value of the fractional order $\alpha$ and need to infer it from the observed data. In the following, we will use the Bayesian inversion method to identify the fractional order $\alpha$ based on the given measurement data $u(x,y, T)$. In Bayesian settings, we take $m=\alpha$, $d=(u(x_1,y_1,T),u(x_2,y_2,T),\cdots, u(x_M,y_M,T))+\eta$, where $M$ is the total number of sensors. And we solve the Bayesian posterior distribution (\ref{IP_4.1}) by an iterative regularization ensemble Kalman method (IREKM), which has been widely used in the identification of fractional order in anomalous diffusion problems, such as \cite{Zhang_Jia_Yan_2018,Yan_Zhang_Wei_2021}.  The idea of the IREKM is to estimate the Bayesian posterior distribution by generating an ensemble of the unknown, where each ensemble member is updated by iterative formula.
Algorithm 1 summarizes the IREKM method's pseudocodes.

\begin{algorithm}[htbp]
\caption{IREKM}
\label{alg2}
\begin{algorithmic}
\STATE \textbf{Input:}\\
~~~~~~~~Initial ensemble of inputs $\{\alpha_{0}^{j}\}_{j=1}^{J}$ from the prior $p_0(\alpha)$. Let $\nu\in(0,1)$ and $\tau>\frac{1}{\nu}$.\\
~~~~~~~~Measurements $d$ and covariance of measurement errors $\Sigma$.\\
~~~~~~~~For $n=0,1,\cdots$
\STATE \textbf{Prediction:}\\
~~~~~~~~Compute $G_n^{j}=G(\alpha_n^{j})$ for $j=\{1,\cdots, J\}$.\\
~~~~~~~~Compute $\bar{G}_n=\frac{1}{J}\sum_{j=1}^{J}G_n^{j}$.\\
\STATE \textbf{Discrepancy principle:}\\
~~~~~~~~Let $\delta=|d-G(\alpha^{\dagger})|_{\Sigma}$, where $\alpha^{\dagger}$ denotes the exact fractional order and $|\cdot|_{\Sigma}=|\Sigma^{-\frac{1}{2}}\cdot|$.\\
~~~~~~~~If $|d-\bar{G}_n|_{\Sigma}\le \tau\delta$, stop and output:\\
\begin{equation*}
\bar{\alpha}_n = \frac{1}{J}\sum_{j=1}^{J}\alpha_n^{j}.
\end{equation*}
\STATE \textbf{Analysis:}\\
~~~~~~~~Update each ensemble member using the following formula
\begin{eqnarray*}
\tilde{\alpha}_{n+1}^{j}&=&\alpha_{n}^{j}+C_{n}^{\alpha G}(C_{n}^{GG}+\mu_n\Sigma)^{-1}(d-G_{n}^{j}), ~j=\{1,2,\cdots, J\}\\
\alpha_{n+1}^{j} &=&\min\{0.999,\max(0.001,\tilde{\alpha}_{n+1}^{j})\}, 
\end{eqnarray*}
~~~~~~~~where
\begin{eqnarray*}
C_{n}^{GG}&=&\frac{1}{J-1}\sum_{n=1}^{J}(G_{n}^{j}-\bar{G}_{n})(G_{n}^{j}-\bar{G}_n)^{T}\\
C_{n}^{\alpha G}&=&\frac{1}{J-1}\sum_{n=1}^{J}(\alpha_{n}^{j}-\bar{\alpha}_{n})(G_{n}^{j}-\bar{G}_n)^{T}.
\end{eqnarray*}
~~~~~~~~~Here $\mu_n$ is chosen as follows:\\
~~~~~~~~~Let $\mu_0$ be an initial guess, and $\mu_n^{i+1}=2^i\mu_0$. Choose $\mu_n=\mu_{n}^{N}$, where $N$ is the first integer s.t.
\begin{equation*}
\mu_n^{N}|\Sigma^{\frac{1}{2}}(C_{n}^{GG}+\mu_{n}^{N}\Sigma)^{-1}(d-\bar{G}_{n})|\ge \nu |\Sigma^{-\frac{1}{2}}(d-\bar{G}_{n})|.
\end{equation*}
\end{algorithmic}
\end{algorithm}

\begin{figure}[H]
	\centering
	\subfigure[]{
		\begin{minipage}[t]{0.3\linewidth}
			\centering
			\includegraphics[width=2in]{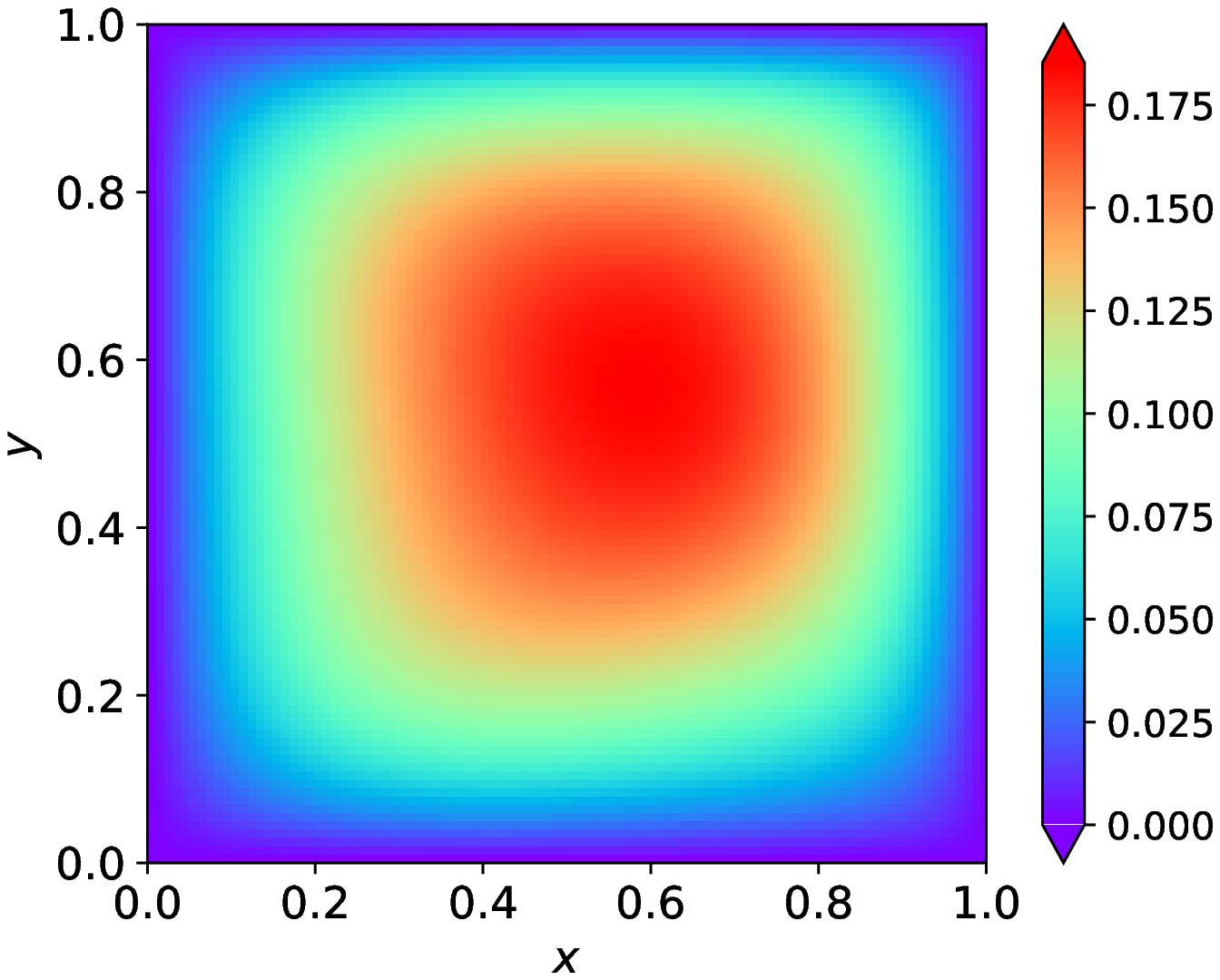}
		\end{minipage}
	}%
	\subfigure[]{
		\begin{minipage}[t]{0.3\linewidth}
			\centering
			\includegraphics[width=2in]{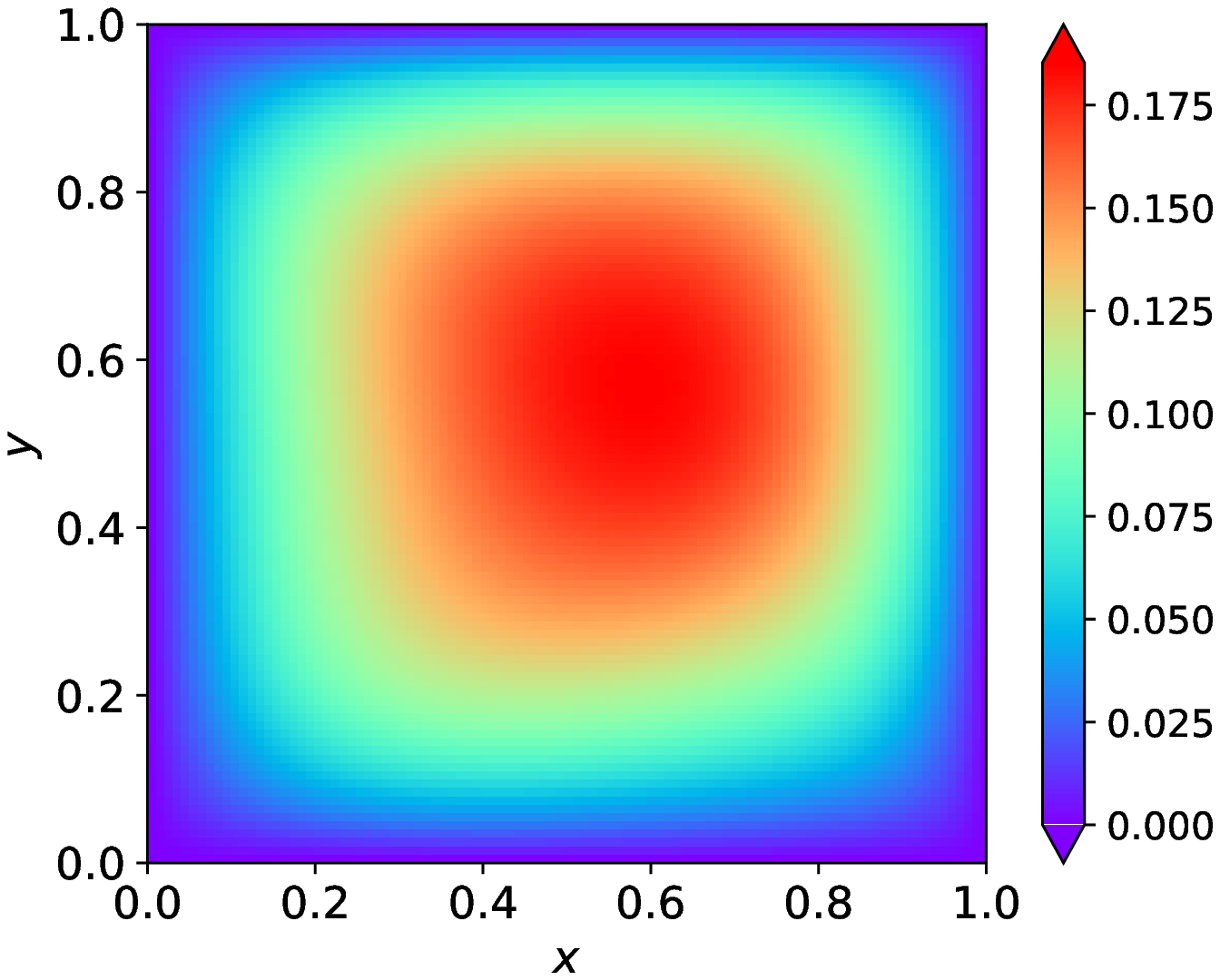}
		\end{minipage}
	}%
	\centering
	\caption{(a) The finite difference solution $u(x,y,t=1)$ when the exact fractional order $\alpha=0.5$ as an input. (b) The finite difference solution $u(x,y,t=1)$ when the approximate fractional order $\alpha=0.4869$ in case of $\sigma=0.003$ as an input.}
	\label{fig_ip_data_compare}
\end{figure}

We fix $c(x,y)=-(xy+4)$, $u_{0}(x,y)=6\sin(2\pi x)\sin(3\pi y)$, and source $f(x,y)=\sin(3\pi x)\sin(\pi y)+6\exp(x^2+y^2)$ in the experiment. In Algorithm 1, we take the prior as $p_0(\alpha)=U(\underline{\alpha},\overline{\alpha})$, where $\underline{\alpha}=0.001$, $\overline{\alpha}=0.999$. In this study, we use the trained deep operator network $G_{NN}^{(\alpha,a)}$ which is used to approximate the solution operator (\ref{3.1.1})
to simulate forward problems that take a long time to solve using traditional numerical methods, where the diffusion coefficient is fixed. Although we know all the information of region $\bar{\Omega}\times [0,T]$ for the approximate solution $u_{NN}(x,y,t;\alpha,a)$ (the outputs of the operator network $G_{NN}^{(\alpha,a)}$), in this inverse problem, we only need the terminal outputs of the approximate solution $u_{NN}(x,y,t;\alpha,a)$, i.e., $u_{NN}(x,y,T;\alpha,a)$. Table \ref{tab0} shows the inversion results of the fractional order $\alpha$ when the standard deviation of the noise is 0.001 and 0.003, respectively.
Furthermore, in Figure \ref{fig_ip_data_compare}, we compare the finite difference
solution $u(x,y,t)$ evaluated at $t=1$ with the exact $\alpha$ and the approximate $\alpha$ as an input, respectively. According to the numerical results, the deep Bayesian inversion method can identify the fractional order $\alpha$ well, and the numerical solution $u(x,y,1)$ generated by the approximate $\alpha$  is nearly identical to the numerical solution $u(x,y,1)$ generated by the exact $\alpha$.

\begin{table}[H]
	\centering
	\setlength{\abovecaptionskip}{0.15cm}
	\setlength{\belowcaptionskip}{0.15cm}
	\caption{Identification of fractional order $\alpha$ in the  subdiffusion problem (\ref{1.1}).}
	\begin{tabular}{l|ccccccccc}
		\hline
		\diagbox{Noise}{$\alpha$} & 0.1    & 0.2    & 0.3    & 0.4    & 0.5    & 0.6    & 0.7    & 0.8    & 0.9    \\
		\hline
		$\sigma=0.001$            & 0.1047 & 0.2056 & 0.3399 & 0.4251 & 0.5076 & 0.6072 & 0.7057 & 0.7814 & 0.8303 \\
		$\sigma=0.003$            & 0.1071 & 0.1949 & 0.3252 & 0.4156 & 0.4869 & 0.5842 & 0.7088 & 0.7892 & 0.8628 \\
		\hline
	\end{tabular}

	\label{tab0}
\end{table}

\subsection{Inverse the diffusion coefficient from the terminal data}
Next, we will consider the inverse diffusion coefficient problem using the known terminal data. Let $m=a(x)$, $d=(u(x_1,y_1,T),u(x_2,y_2,T),\cdots,u(x_M,y_M,T))+\eta$, $M$ represents the number of measurement points. To solve the posterior probability distribution, we apply a function space Markov chain Monte Carlo (MCMC) approach in \cite{Cotter_Roberts_Stuart_White_2013}. 
Algorithm 2 summarizes the pseudocodes of the MCMC method. 

\begin{algorithm}[H]
\caption{MCMC-pCN}
\label{alg2}
\begin{algorithmic}
\STATE Let $X$ is a Hilbert space, $N(0,C)$ is a Gaussian measure. \\
\STATE Define $\rho(m, \tilde{m})=\min\{1,\exp(\Phi(m)-\Phi(\tilde{m}))\}$ where $\Phi(m)=\frac{1}{2}|d-g(m)|_{\Sigma}^2$. The sequences $\{m^{(0)}\}_{k\ge 0}$ are generated as follows:\\
\STATE 1. Set $k=0$ and pick $m^{0}\in X$;\\
\STATE 2. Propose $\tilde{m}^{k}=\sqrt{(1-\beta^2)}m^{k}+\beta \upsilon^{k}$, $\upsilon^{k}\sim N(0,C)$;\\
\STATE 3. Set $m^{k+1}=\tilde{m}^{k}$ with the probability $\rho(m^{k},\tilde{m}^{k})$;\\
\STATE 4. Set $m^{k+1}=m^{k}$ otherwise.
\end{algorithmic}
\end{algorithm}

It is well known that implementing the MCMC method necessitates repeatedly solving many forward problems. Furthermore, due to the influence of the time-fractional derivative in the subdiffusion problem (\ref{1.1}), traditional numerical methods require a long time to solve the forward problems, which poses significant challenges to the implementation of the MCMC method. Here, we seek to accelerate the implementation of the MCMC approach by employing the trained deep operator network $G_{NN}^{(\alpha,a)}$ as an approximate operator to the forward operator $F$.

\begin{figure} 
	\centering
	\subfigure[Diffusion coefficient]{
		\begin{minipage}[t]{0.3\linewidth}
			\centering
			\includegraphics[width=2in]{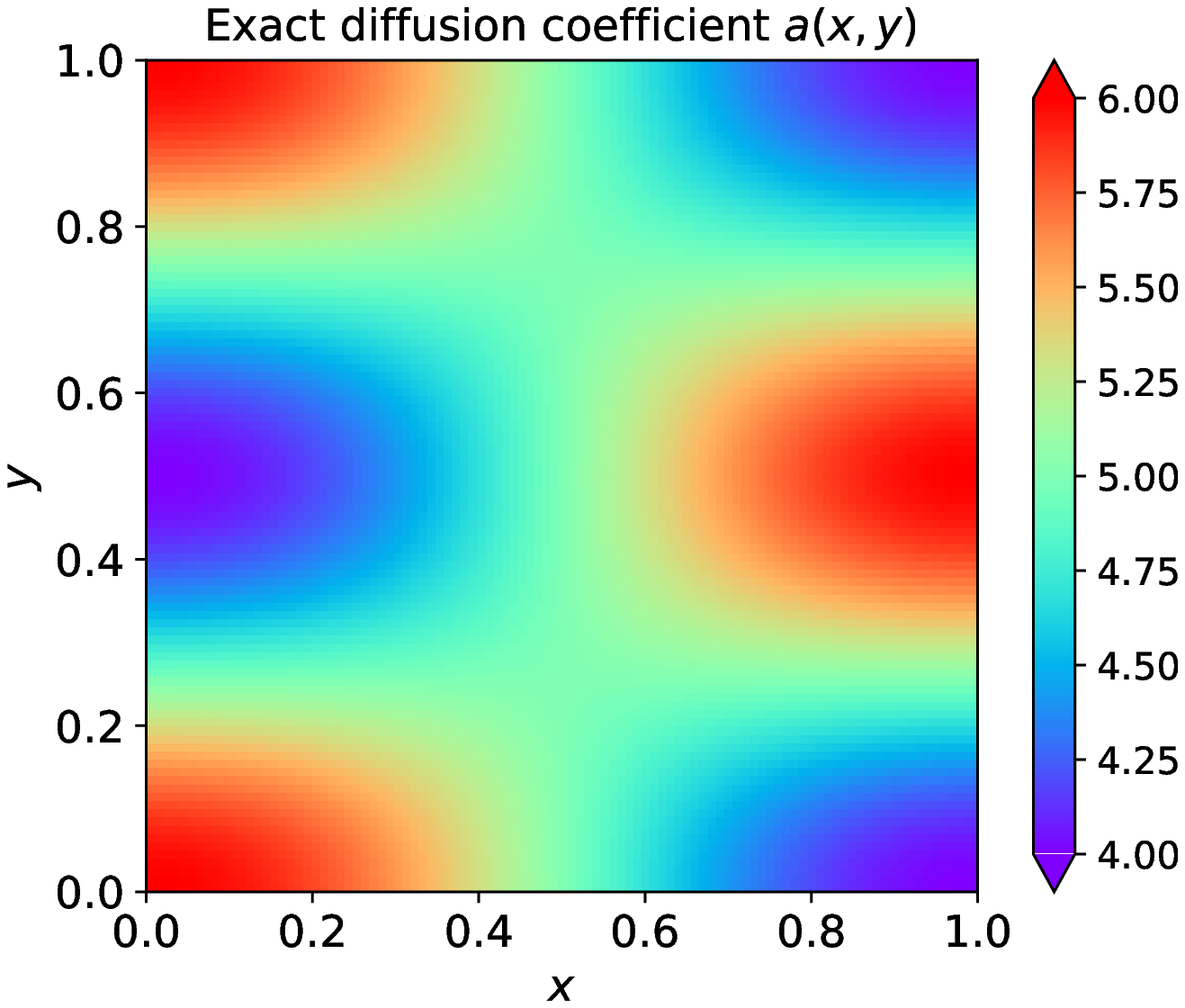}
		\end{minipage}
	}%
	\subfigure[Exact $u(x,y;T)$]{
		\begin{minipage}[t]{0.3\linewidth}
			\centering
			\includegraphics[width=2in]{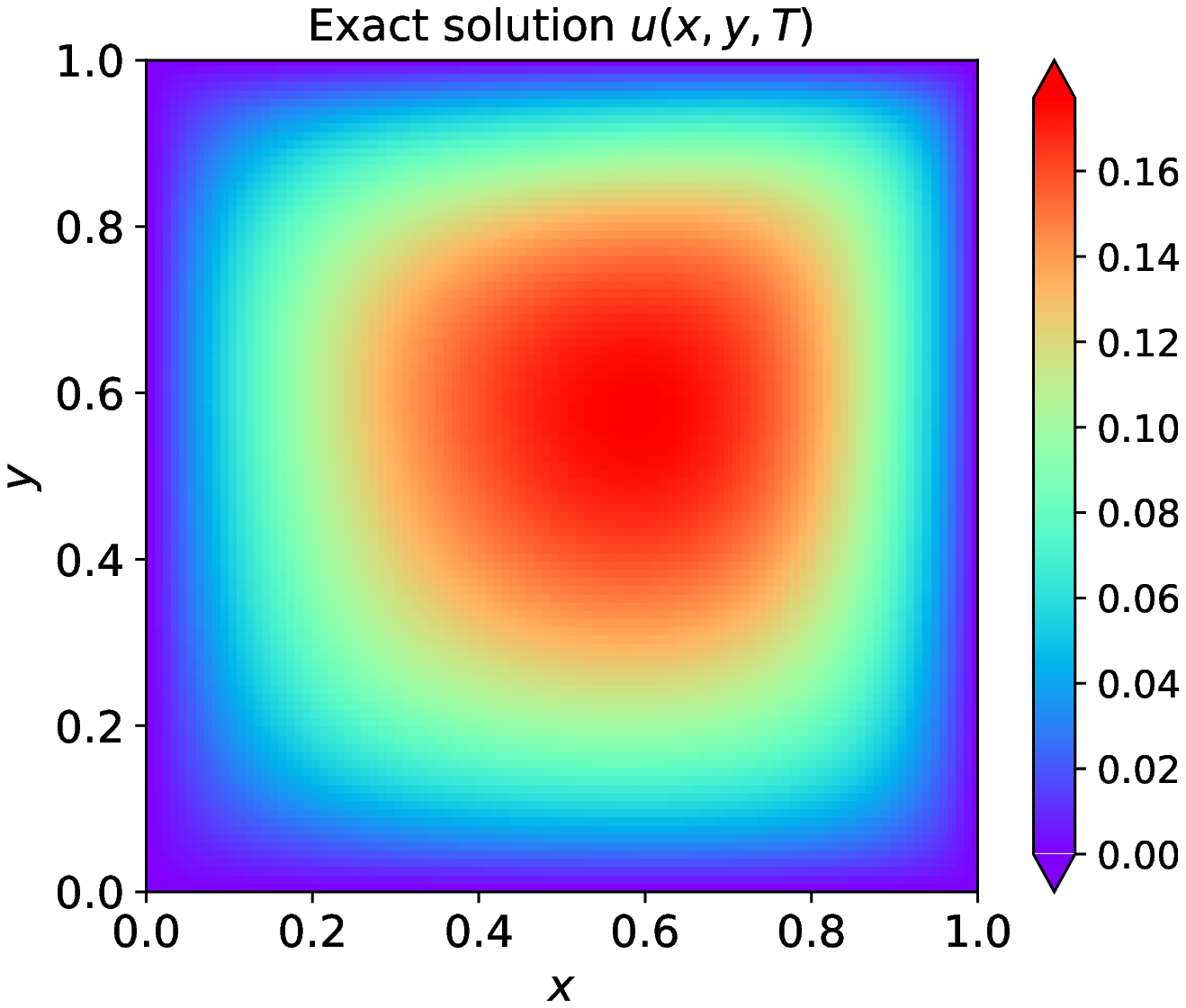}
		\end{minipage}
	}%

	\subfigure[Noise data with $\delta=0.001$]{
		\begin{minipage}[t]{0.3\linewidth}
			\centering
			\includegraphics[width=2in]{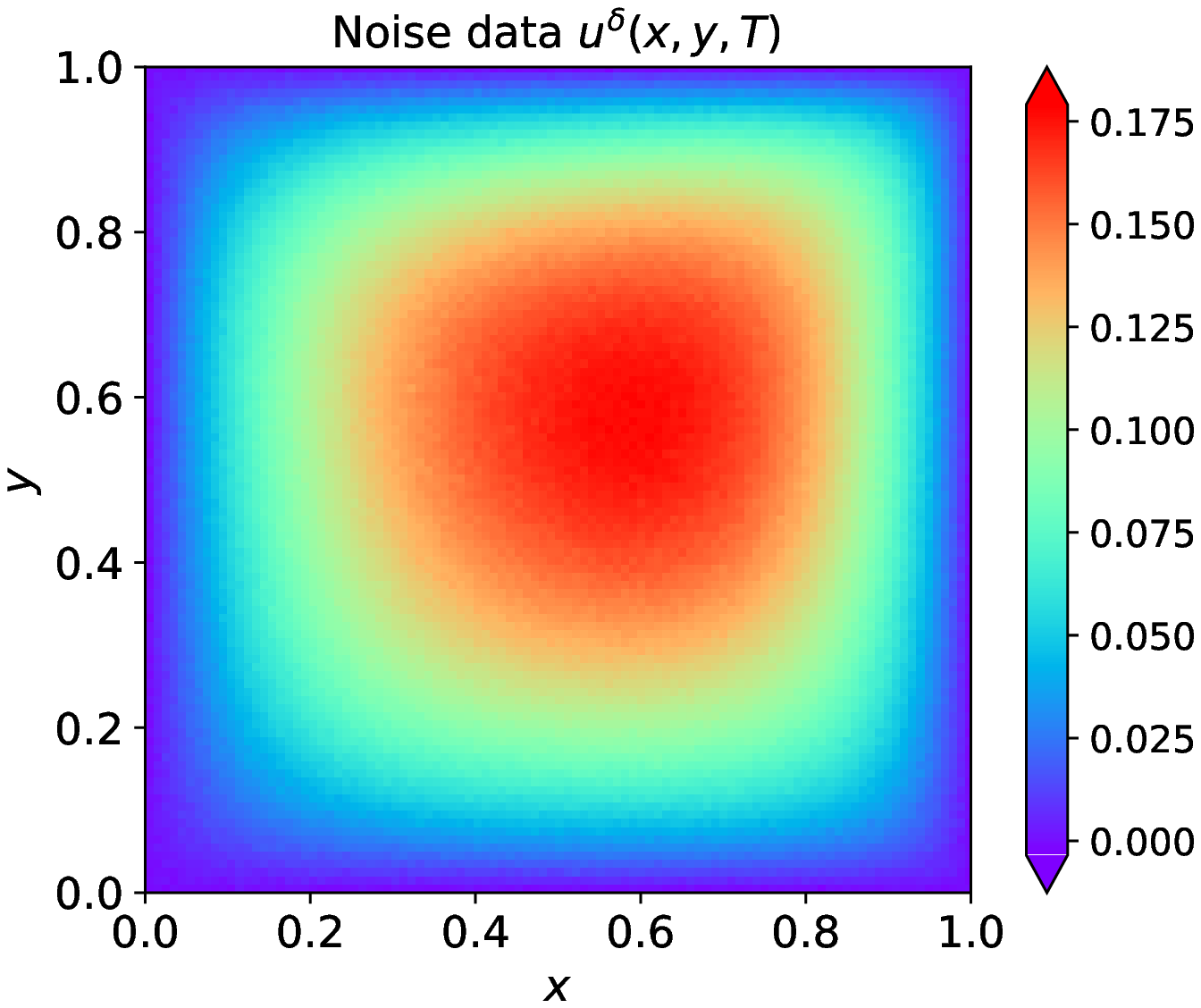}
		\end{minipage}
	}%
	\subfigure[Noise data with $\delta=0.005$]{
		\begin{minipage}[t]{0.3\linewidth}
			\centering
			\includegraphics[width=2in]{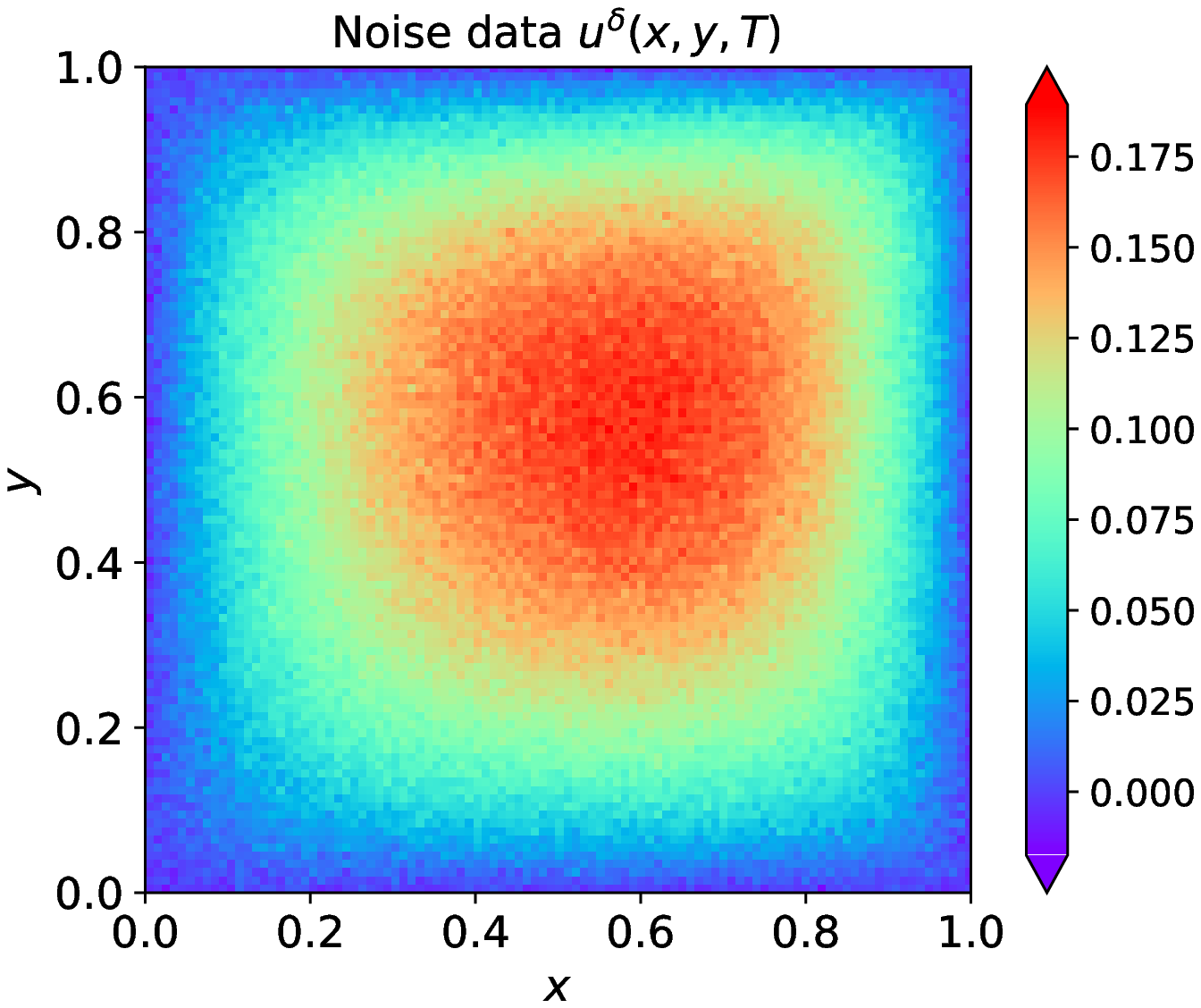}
		\end{minipage}
	}%
	\centering
	\caption{(a) The exact diffusion coefficient to be identified. (b) The exact data with free noise.
		(c) The terminal measurement data was polluted by noise with 0-mean and standard deviation $\sigma=0.001$. (d)
		The terminal measurement data was polluted by noise with 0-mean and standard deviation $\sigma=0.005$.}
	\label{fig4.1_data}
\end{figure}

\begin{figure} 
	\centering
	\subfigure[Inverse $a(x,y)$ in $\delta=0.001$]{
		\begin{minipage}[t]{0.3\linewidth}
			\centering
			\includegraphics[width=2in]{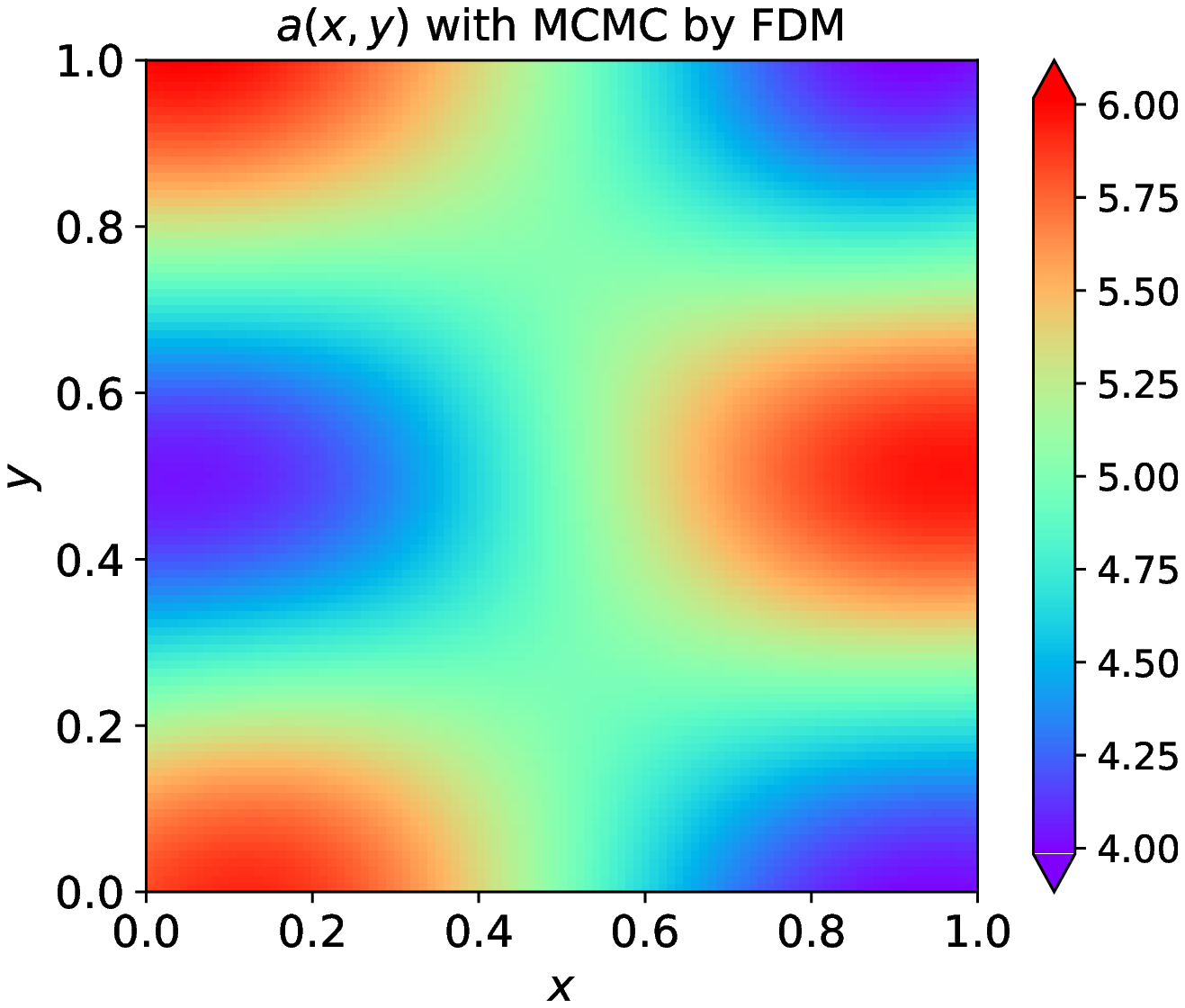}
		\end{minipage}
	}%
	\subfigure[Point-wise errors]{
		\begin{minipage}[t]{0.3\linewidth}
			\centering
			\includegraphics[width=2in]{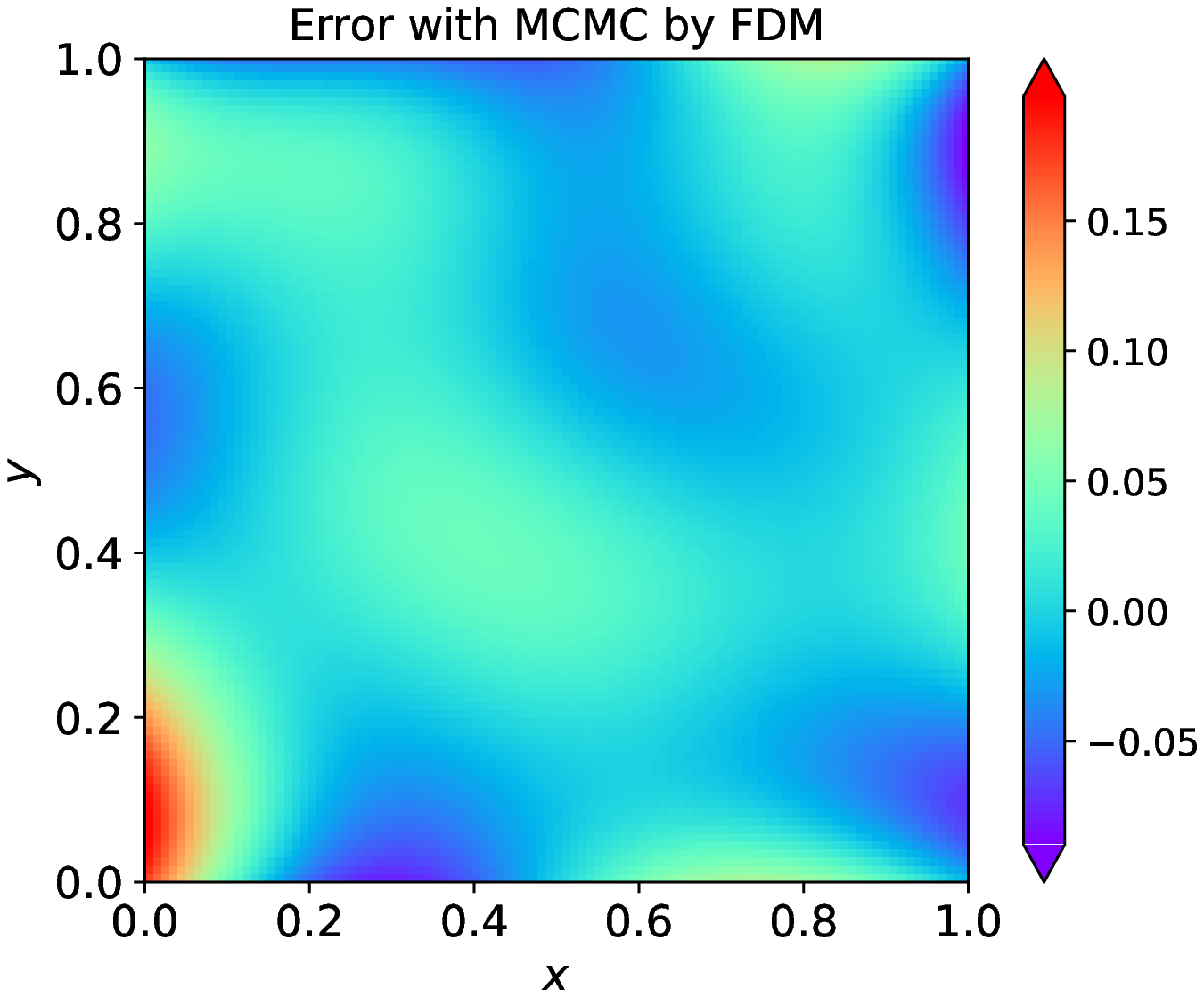}
		\end{minipage}
	}%

	\subfigure[Inverse $a(x,y)$ in $\delta=0.005$]{
		\begin{minipage}[t]{0.3\linewidth}
			\centering
			\includegraphics[width=2in]{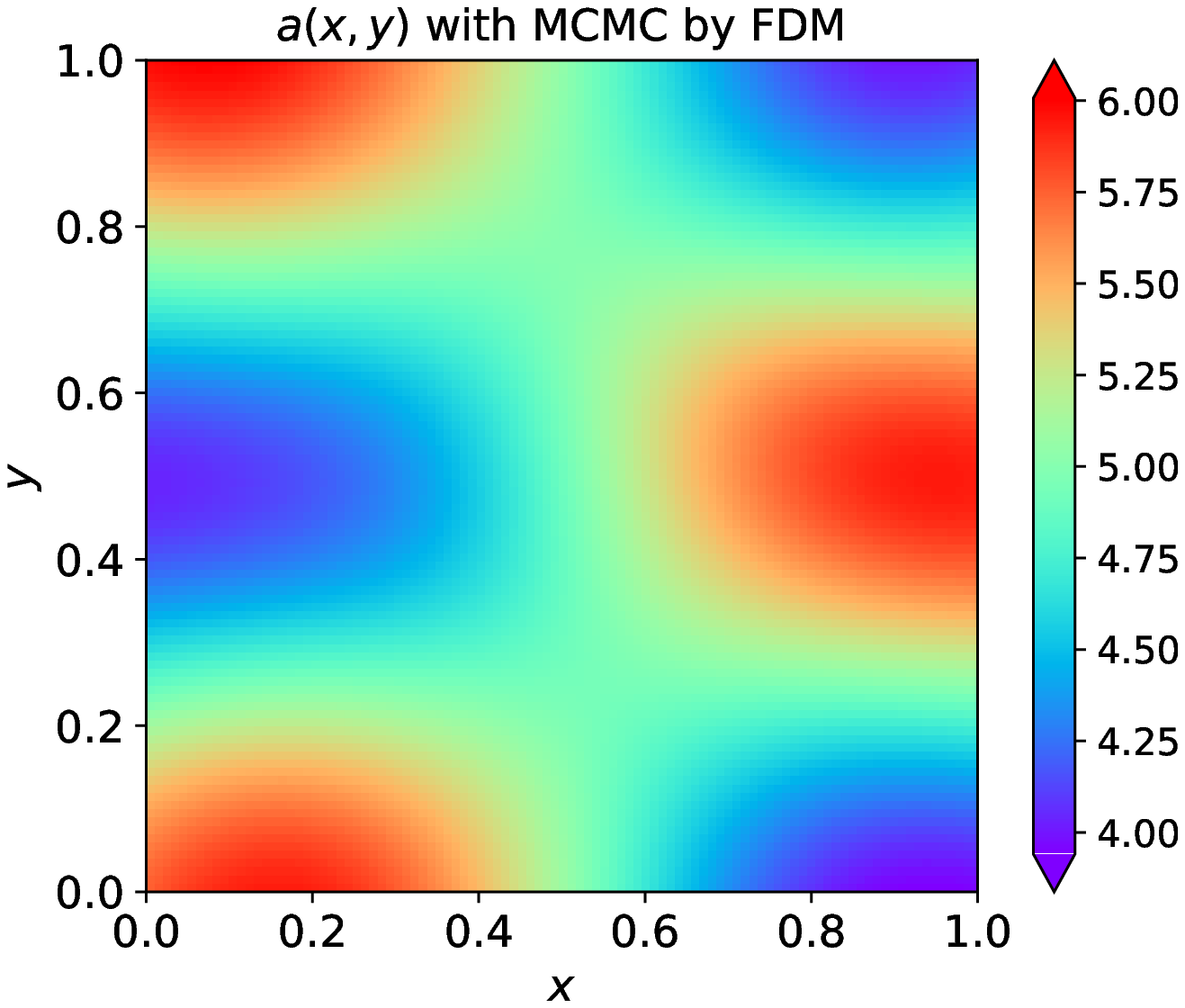}
		\end{minipage}
	}%
	\subfigure[Point-wise errors]{
		\begin{minipage}[t]{0.3\linewidth}
			\centering
			\includegraphics[width=2in]{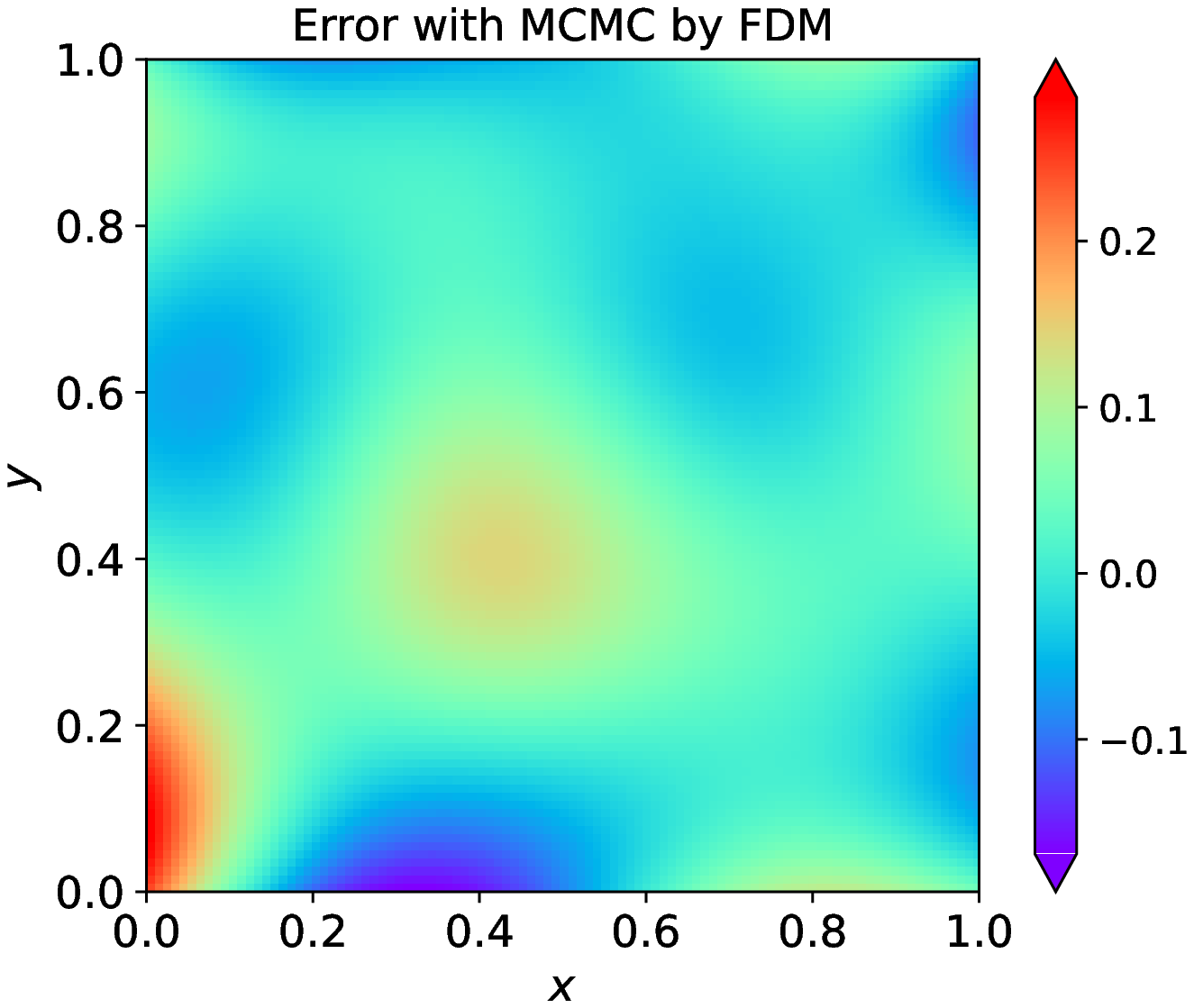}
		\end{minipage}
	}%

	\centering
	\caption{(a), (c) The reconstructed diffusion coefficient by MCMC+FDM from the terminal measurement data that was polluted by noise with 0-mean, and the standard deviation $\sigma$ is $0.001$ and $0.005$, respectively. (b), (d) The difference between the exact diffusion coefficient and reconstructed diffusion coefficient by MCMC+FDM from measurement data that was polluted by noise with 0-mean, and the standard deviation $\sigma$ is $0.001$ and $0.005$, respectively.}
	\label{fig4.1_FDM}
\end{figure}

\begin{figure}[H]
	\centering
	\subfigure[Inverse $a(x,y)$ in $\delta=0.001$]{
		\begin{minipage}[t]{0.3\linewidth}
			\centering
			\includegraphics[width=2in]{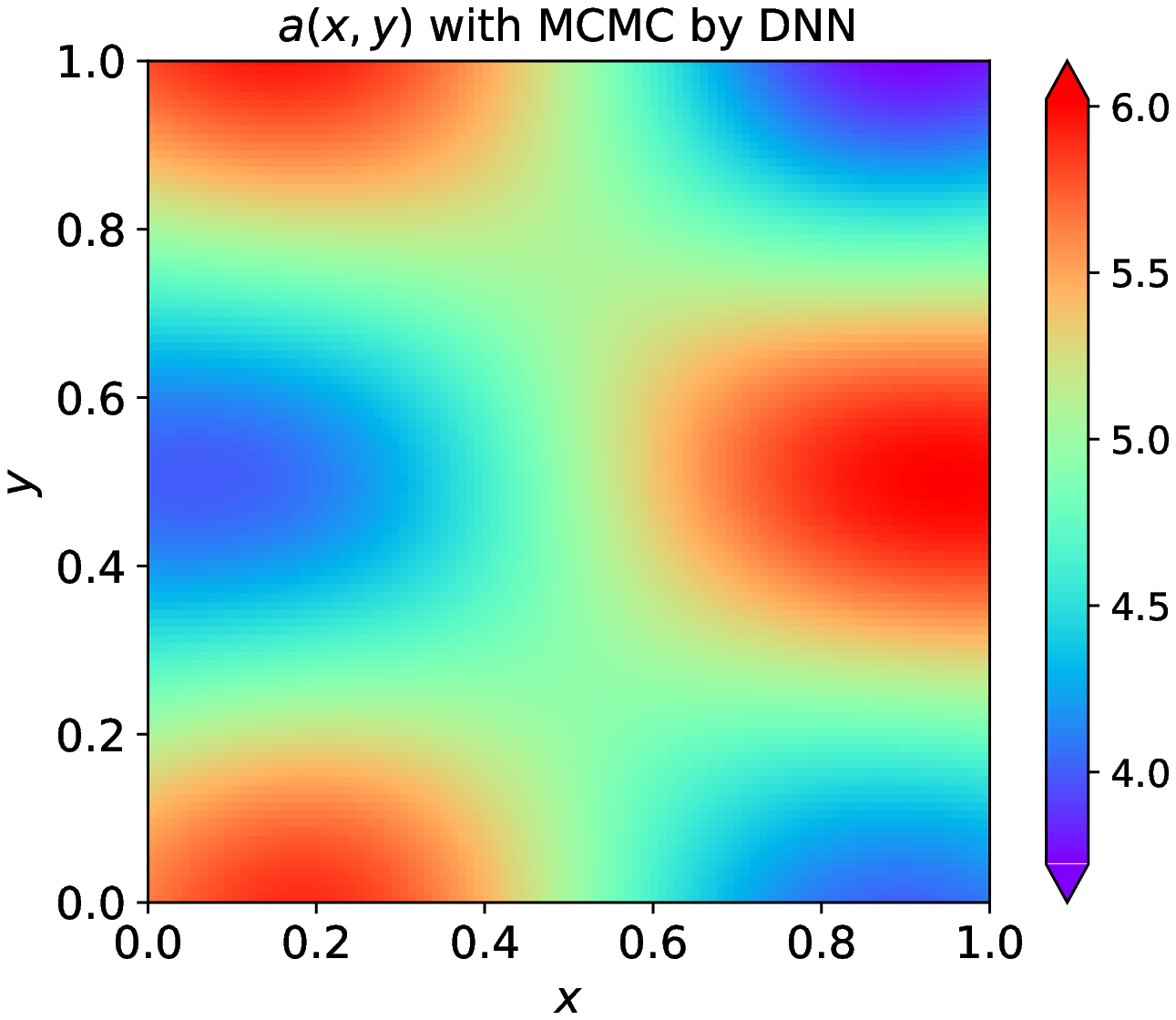}
		\end{minipage}
	}%
	\subfigure[Point-wise errors]{
		\begin{minipage}[t]{0.3\linewidth}
			\centering
			\includegraphics[width=2in]{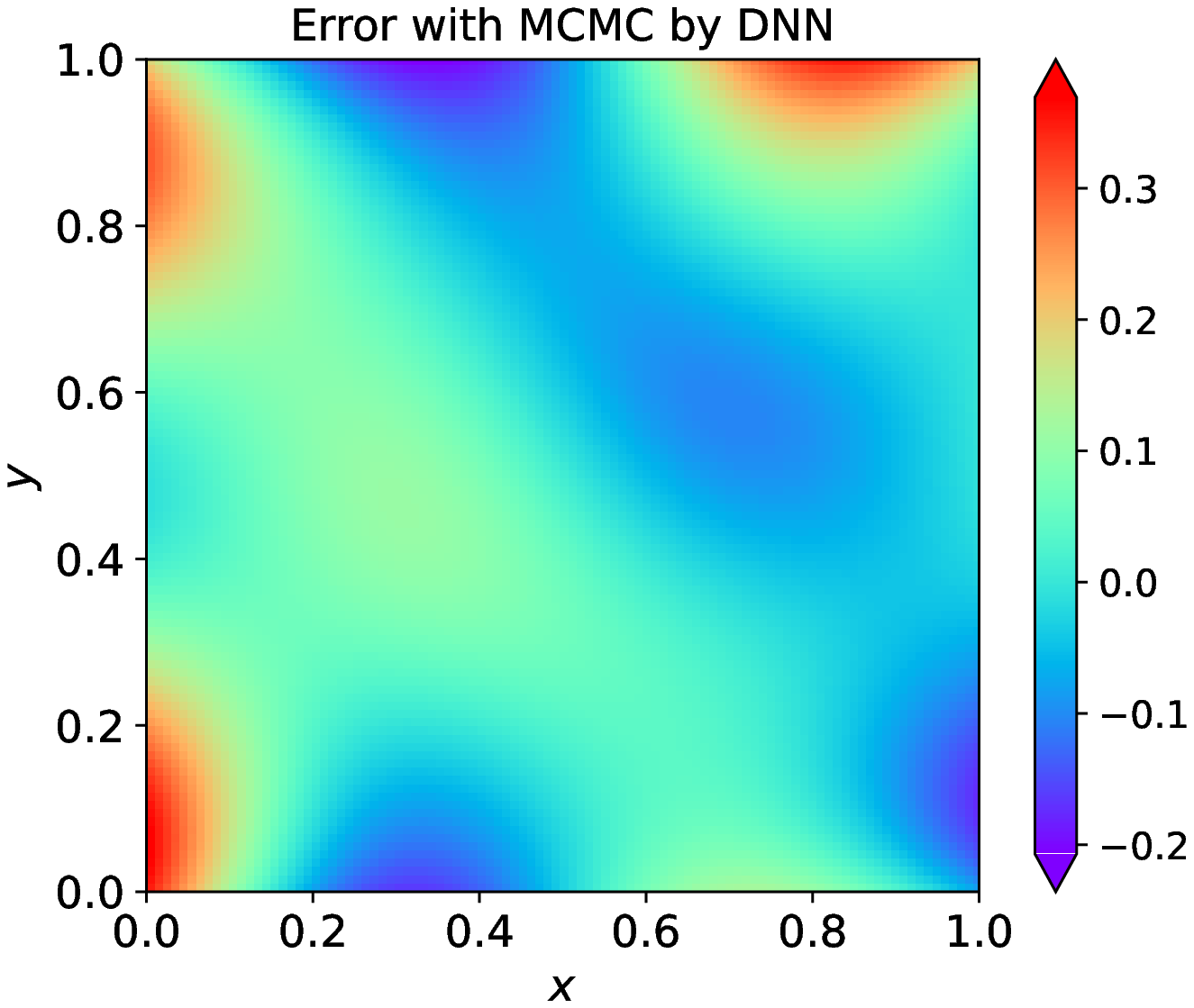}
		\end{minipage}
	}%

	\subfigure[Inverse $a(x,y)$ in $\delta=0.005$]{
		\begin{minipage}[t]{0.3\linewidth}
			\centering
			\includegraphics[width=2in]{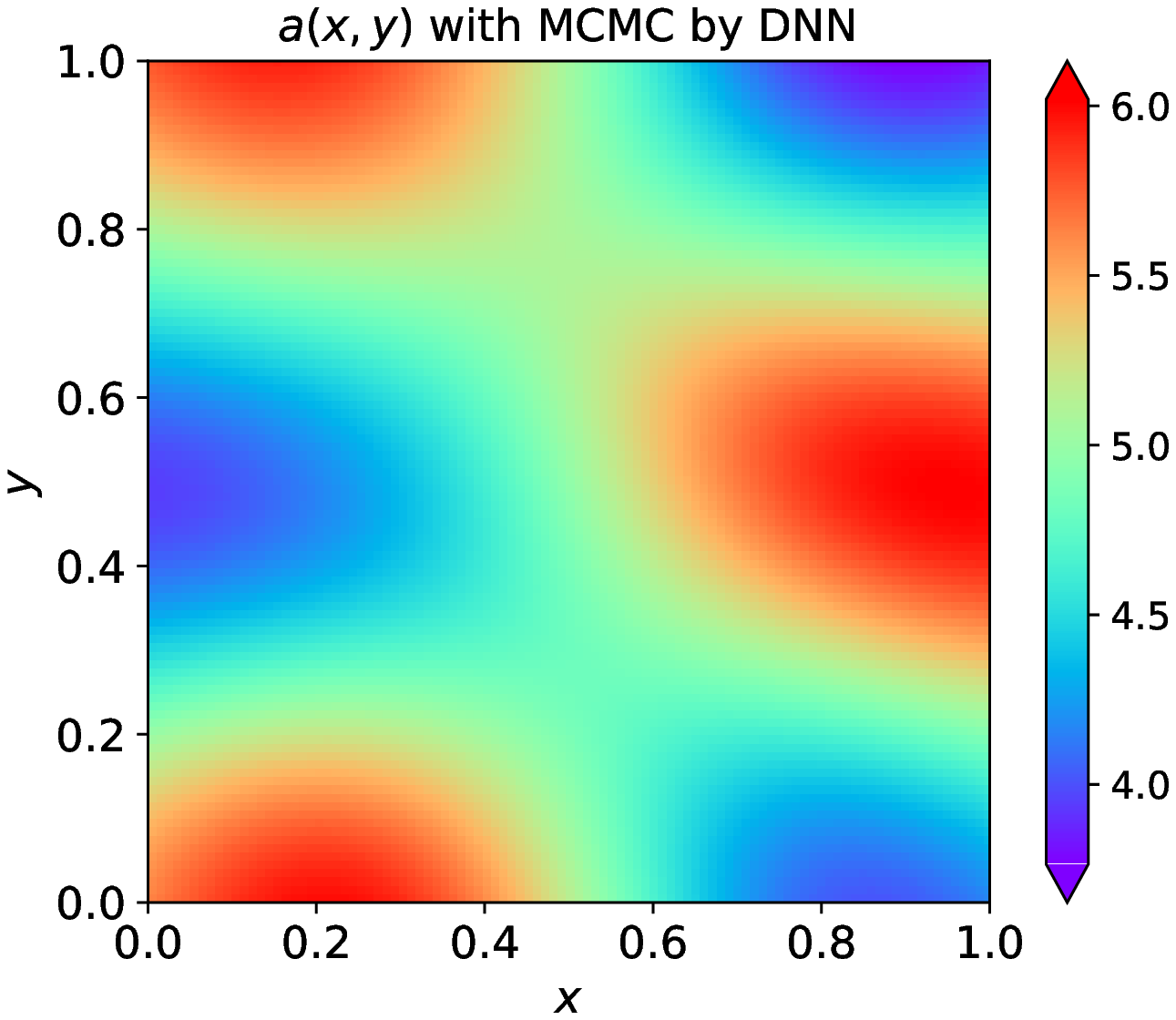}
		\end{minipage}
	}%
	\subfigure[Point-wise errors]{
		\begin{minipage}[t]{0.3\linewidth}
			\centering
			\includegraphics[width=2in]{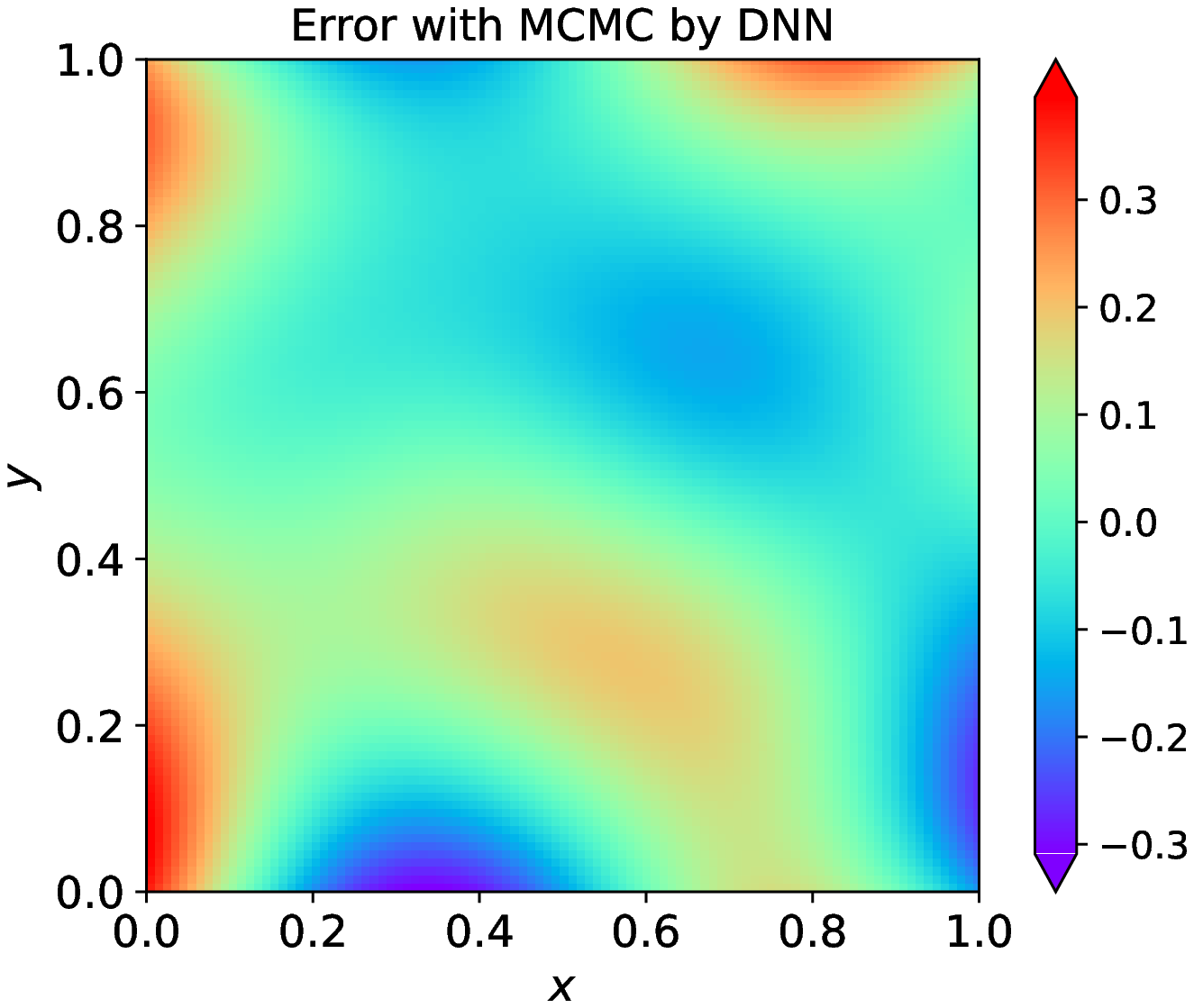}
		\end{minipage}
	}%

	\centering
	\caption{(a), (c) The reconstructed diffusion coefficient by MCMC+$G_{NN}^{(\alpha,a)}$ from the terminal measurement data that was polluted by noise with 0-mean, and the standard deviation $\sigma$ is $0.001$ and $0.005$, respectively. (b), (d) The difference between the exact diffusion coefficient and the reconstructed diffusion coefficient by MCMC+$G_{NN}^{(\alpha,a)}$ from measurement data that was polluted by noise with 0-mean, and the standard deviation $\sigma$ is $0.001$ and $0.005$, respectively.}
	\label{fig4.1_DNN}
\end{figure}

We fix fractional order $\alpha=0.5$, $c(x,y)=-(xy+4)$, $u_{0}(x,y)=6\sin(2\pi x)\sin(3\pi y)$, and source $f(x,y)=\sin(3\pi x)\sin(\pi y)+6\exp(x^2+y^2)$ in the experiment. In Algorithm 2, we take $C=k_{l}(x^{[1]},x^{[2]})$, $l=0.3$, $\beta=0.005$, where $x^{[1]}=(x_1,y_1)$, $x^{[2]}=(x_2, y_2)$. In addition, we run the MCMC sampling in Algorithm 2 for 10000 iterations, and the last 8000 realizations are used to compute the mean of the samples. The exact diffusion coefficient, exact measurement data, and different noises data are shown in Figure \ref{fig4.1_data}. The numerical inversion results under various noises produced by combining the MCMC method with the finite difference method and the MCMC method with the deep operator learning method are shown in Figure \ref{fig4.1_FDM} and Figure \ref{fig4.1_DNN}, respectively. 
We list the cost time of the MCMC sampling carried out by the deep operator learning approach and the finite difference method in Table \ref{table_1}. In comparison to the time needed by the finite difference method, the MCMC sampling with the operator network $G_{NN}^{(\alpha,a)}$ only takes 67 seconds, which is quite short. Even though accounting for the data generation and neural network training time which cost 24484 seconds, the time required for the MCMC sampling with the deep operator learning method is far less than the time 66303 seconds spent for the MCMC sampling with the finite difference method.

To solve the inverse problem more effectively, we retrain a new deep operator network that maps any diffusion coefficients $a(x,y)$ to solution $u(x,y,T)$, i.e.,
\begin{equation}\label{4.1}
	G: a(x,y) \to u(x,y,T), \quad (x,y)\in \bar{\Omega},
\end{equation}
where we denote the approximate deep operator network as $G_{NN}^{a}$. Here, we use four hidden layers with size 256 for the branch net and trunk net. The optimizer and learning rate during neural network training is Adam and $1\times 10^{-4}$, respectively. 
Similarly, in order to generate the training dataset, we fix $\alpha=0.5$, $c(x,y)=-(xy+4)$, $u_{0}(x,y)=6\sin(2\pi x)\sin(3\pi y)$, and source $f(x,y)=\sin(3\pi x)\sin(\pi y)+6\exp(x^2+y^2)$ in the experiment. In addition, we take 1000 samples for aims of training. The settings of the sampling diffusion coefficient $a(x,y)$ and the numerical solutions $u(x,y,T)$ are identical to those in subsection \ref{sub_sec3.1}.

\begin{figure} 
	\centering
	\subfigure[Inverse $a(x,y)$ in $\delta=0.001$]{
		\begin{minipage}[t]{0.3\linewidth}
			\centering
			\includegraphics[width=2in]{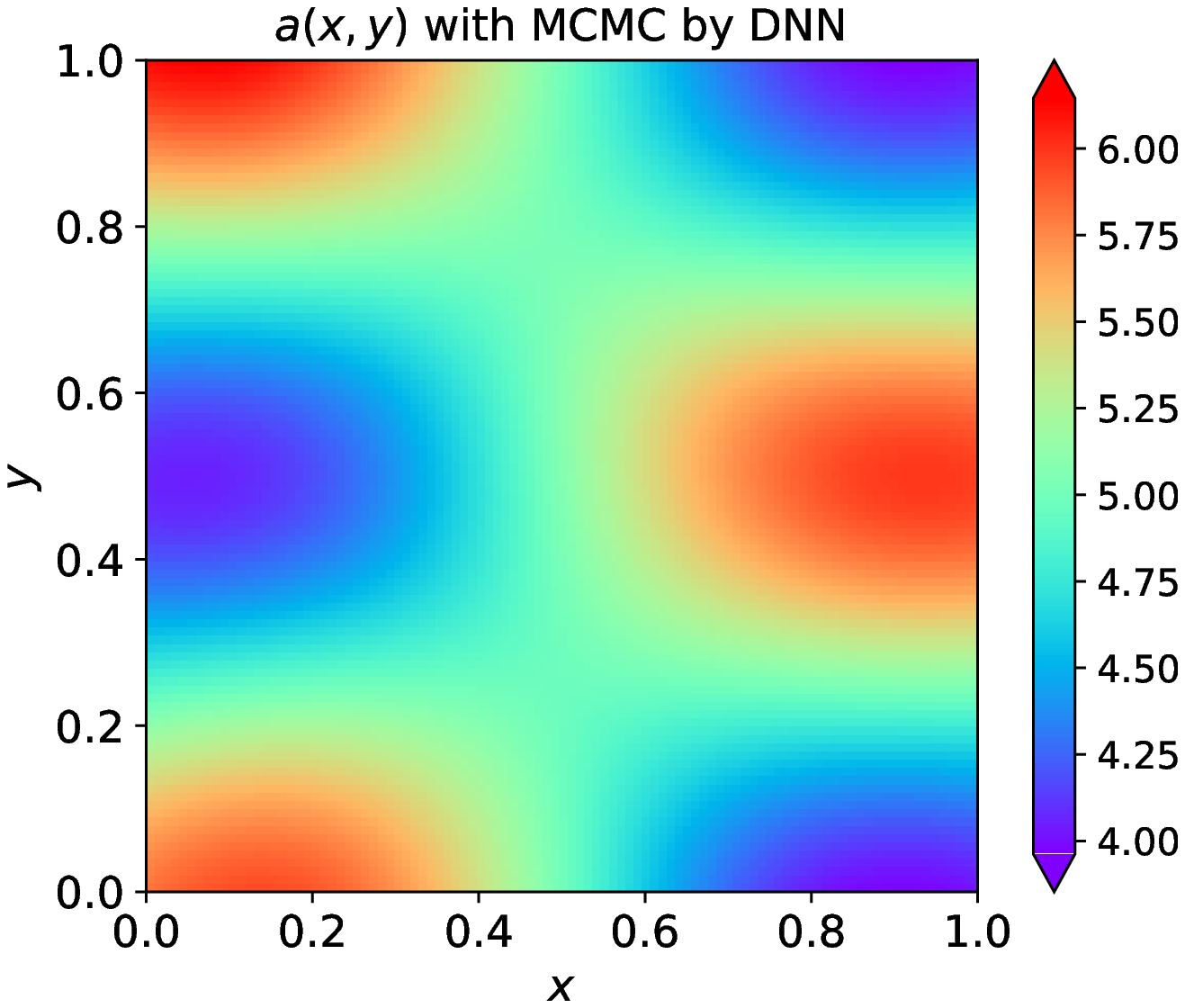}
		\end{minipage}
	}%
	\subfigure[Point-wise errors]{
		\begin{minipage}[t]{0.3\linewidth}
			\centering
			\includegraphics[width=2in]{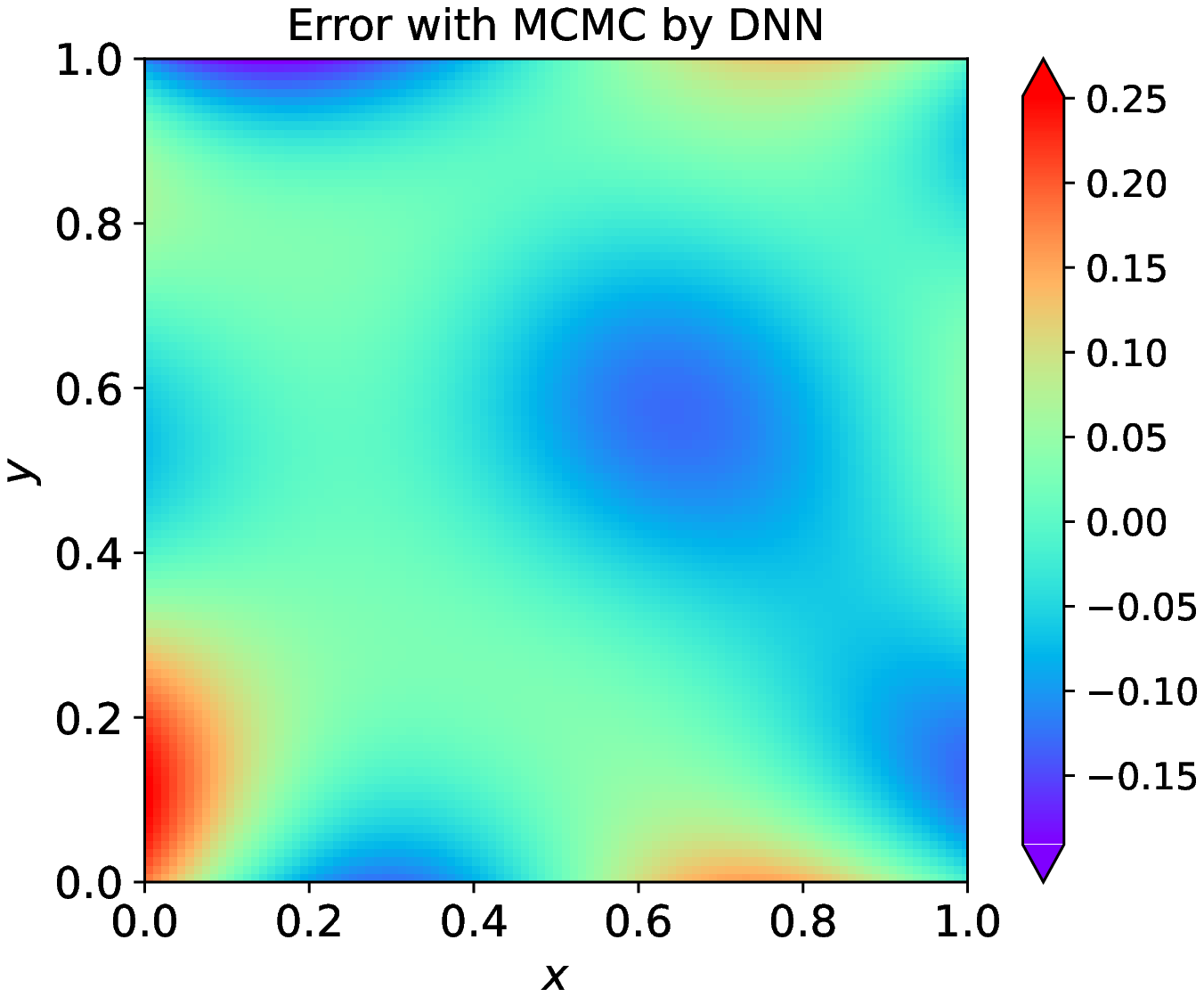}
		\end{minipage}
	}%

	\subfigure[Inverse $a(x,y)$ in $\delta=0.005$]{
		\begin{minipage}[t]{0.3\linewidth}
			\centering
			\includegraphics[width=2in]{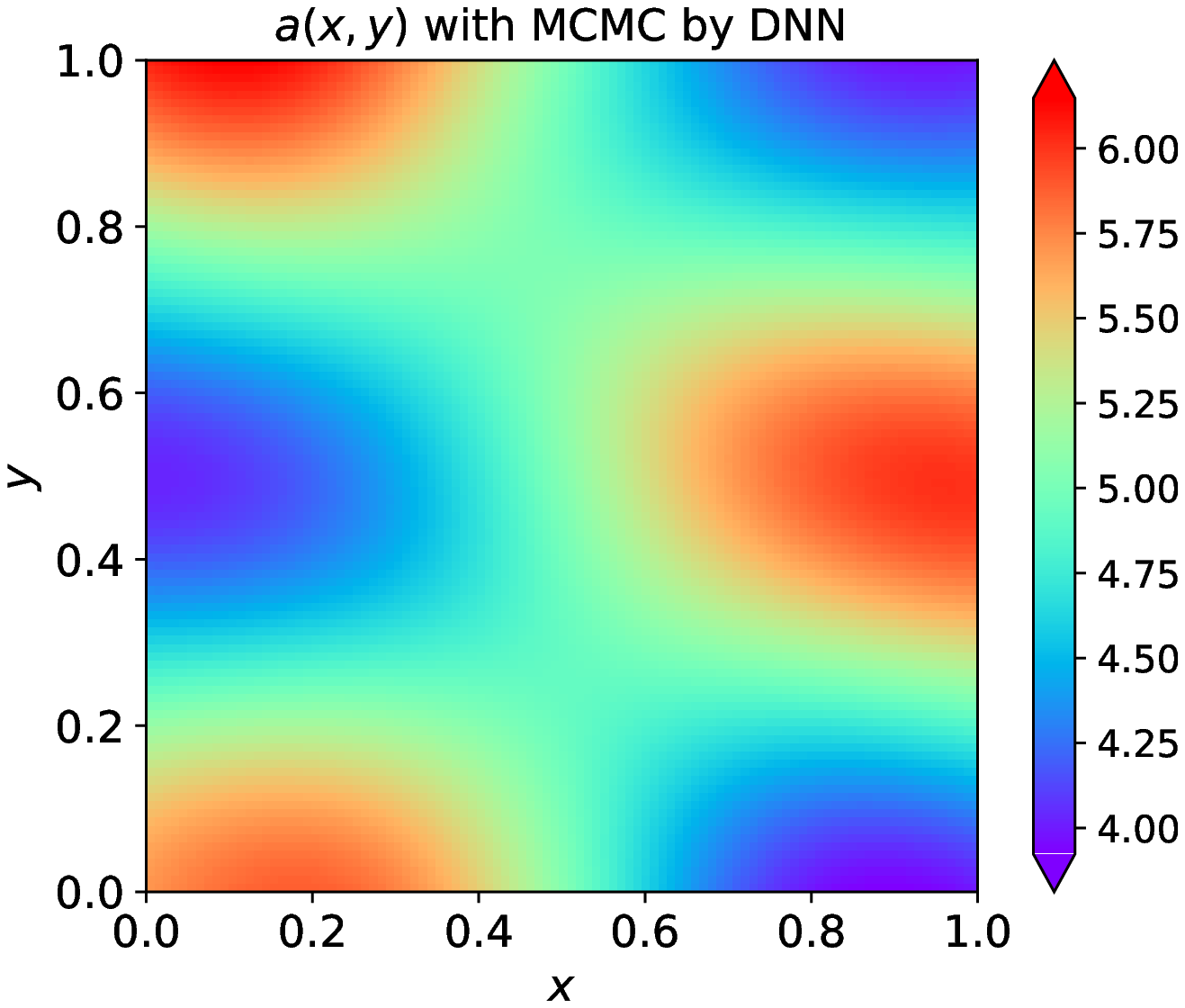}
		\end{minipage}
	}%
	\subfigure[Point-wise errors]{
		\begin{minipage}[t]{0.3\linewidth}
			\centering
			\includegraphics[width=2in]{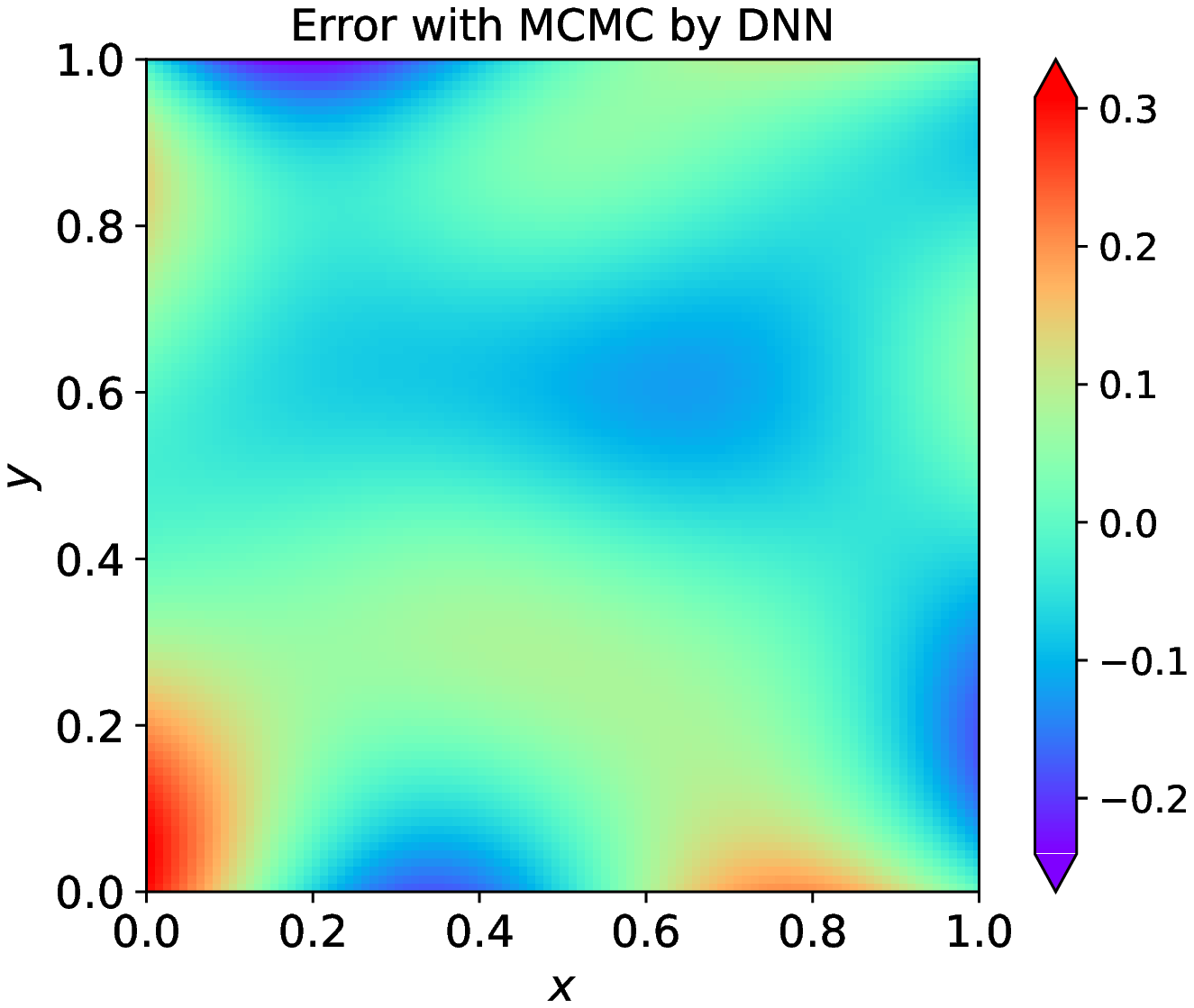}
		\end{minipage}
	}%

	\centering
	\caption{(a), (c) The reconstructed diffusion coefficient by MCMC+$G_{NN}^{a}$ from the terminal measurement data that was polluted by noise with 0-mean, and the standard deviation $\sigma$ is $0.001$ and $0.005$, respectively. (b), (d) The difference between the exact diffusion coefficient and the reconstructed diffusion coefficient by MCMC+$G_{NN}^{a}$ from the terminal measurement data that was polluted by noise with 0-mean, and the standard deviation $\sigma$ is $0.001$ and $0.005$, respectively.}
	\label{fig4.1_DNN_T}
\end{figure}

Figure \ref{fig4.1_DNN_T} shows the numerical inversion results obtained by the MCMC method with the deep operator network $G_{NN}^a$. In addition, in Table \ref{tab01}, we list the relative $l_2$ errors of numerical inversion results generated by above methods.
In comparison to the MCMC method with the deep operator network $G_{NN}^{(\alpha,a)}$, the MCMC sampling with the deep operator network $G_{NN}^{a}$ produces more accurate numerical inversion results. 
Moreover, the cost time of data generation and network training is 6881 seconds, and the inference time is 42 seconds, which is less than the time 66303 seconds required by the MCMC sampling with the finite difference method.

\begin{table}[H]
	\centering
	\setlength{\abovecaptionskip}{0.2cm}
	\setlength{\belowcaptionskip}{0.2cm}
	\caption{Computational times, in seconds, given by three different methods. MCMC+FDM, MCMC+$G_{NN}^{(\alpha,a)}$, MCMC+$G_{NN}^{a}$ indicate that MCMC sampling is carried out using the finite difference method, the operator network $G_{NN}^{(\alpha,a)}$,  and the operator network $G_{NN}^{a}$ to simulate the forward problems.}
	\setlength{\tabcolsep}{3mm}{
		\begin{tabular}{c c c c c}
			\toprule
			Methods                    & Inference & Speed                 & \\
			\midrule
			MCMC+FDM                   & 66303     & -                       \\
			MCMC+$G_{NN}^{(\alpha,a)}$ & 67        & $\approx$1000X faster   \\
			MCMC+$G_{NN}^{a}$          & 42        & $\approx$1500X faster   \\
			\midrule
		\end{tabular}}%
	\label{table_1}%
\end{table}%

\begin{table}[H]
	\centering
	\setlength{\abovecaptionskip}{0.2cm}
	\setlength{\belowcaptionskip}{0.2cm}
	\caption{The relative $l_2$ error between the reconstructed diffusion coefficient and the exact diffusion coefficient for the different inversion methods. MCMC+FDM, MCMC+$G_{NN}^{(\alpha,a)}$, MCMC+$G_{NN}^{a}$ indicate that MCMC sampling is carried out using the finite difference method, the operator network $G_{NN}^{(\alpha,a)}$,  and the operator network $G_{NN}^{a}$ to simulate the forward problems.}
	\begin{tabular}{l|ccc}
		\hline
		\diagbox{Noise}{Relative $l_2$ error}{Methods} & MCMC+FDM & MCMC+$G_{NN}^{(\alpha,a)}$ & MCMC+$G_{NN}^{a}$ \\
		\hline
		$\sigma=0.001$                                 & 0.006282 & 0.018877                   & 0.011988          \\
		$\sigma=0.005$                                 & 0.011683 & 0.022514                   & 0.014754          \\
		\hline
	\end{tabular}\label{tab01}
\end{table}
\subsection{Inverse diffusion coefficient from interior measurements}
\begin{figure} 
	\centering
	{\epsfxsize 0.3\hsize \epsfbox{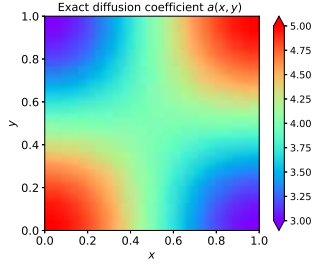}}
	\caption{The exact diffusion coefficient to be identified.}\label{fig4.2_DNN_inDomain_diffu}
\end{figure}

\begin{figure} 
	\centering
	\subfigure[Inverse $a(x,y)$ in $\delta=0.001$]{
		\begin{minipage}[t]{0.3\linewidth}
			\centering
			\includegraphics[width=2in]{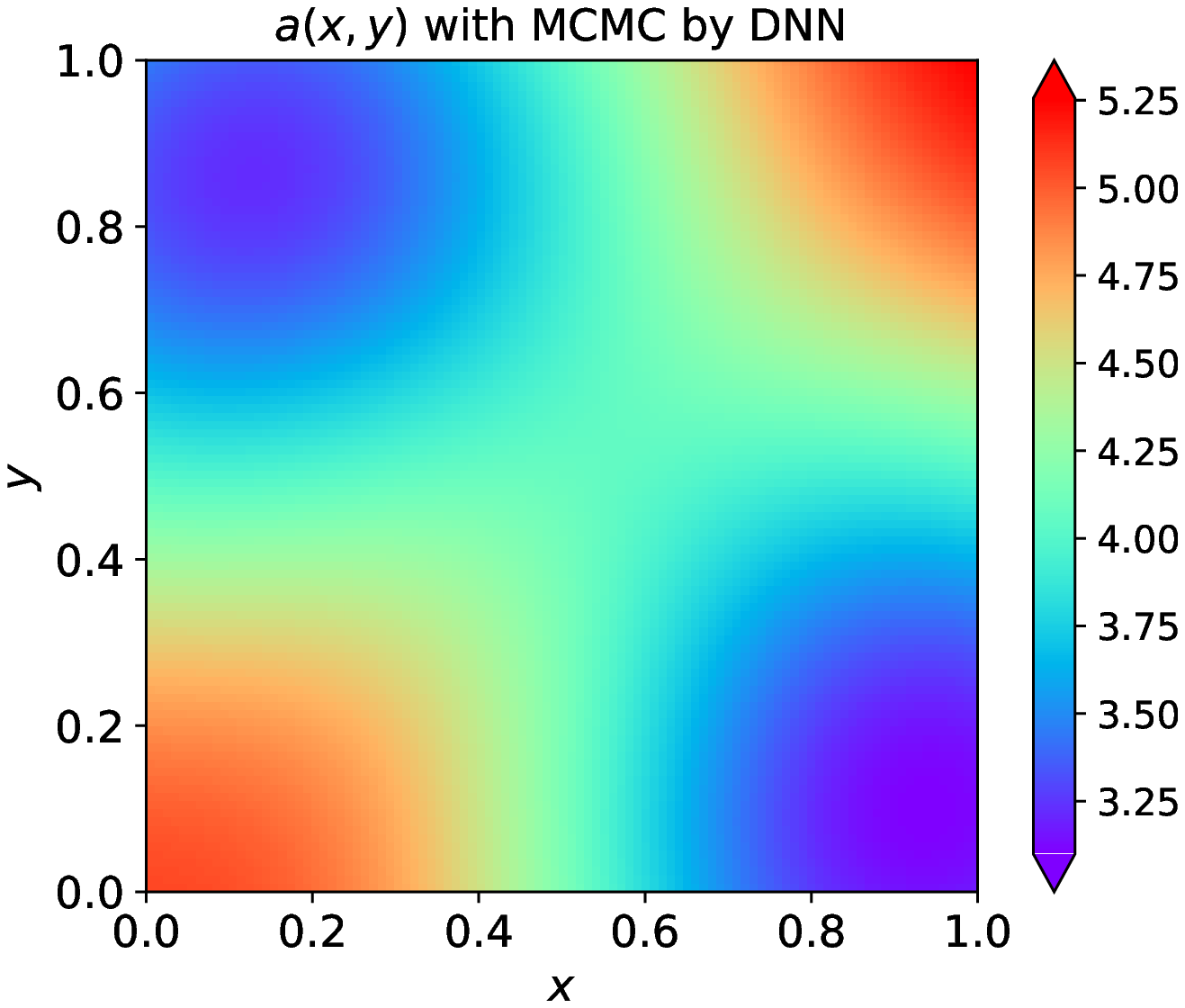}
		\end{minipage}
	}%
	\subfigure[Point-wise errors]{
		\begin{minipage}[t]{0.3\linewidth}
			\centering
			\includegraphics[width=2in]{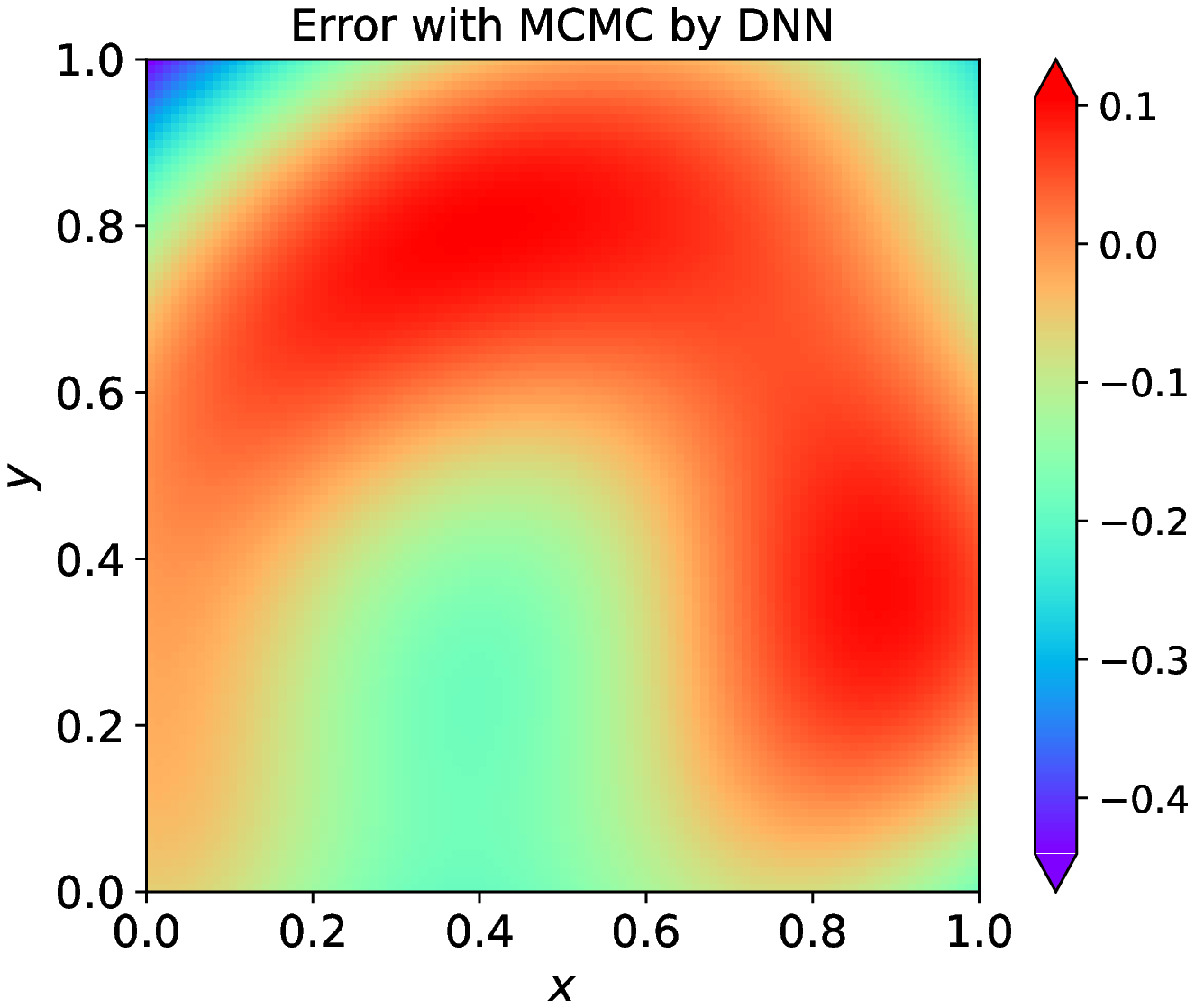}
		\end{minipage}
	}%

	\subfigure[Inverse $a(x,y)$ in $\delta=0.005$]{
		\begin{minipage}[t]{0.3\linewidth}
			\centering
			\includegraphics[width=2in]{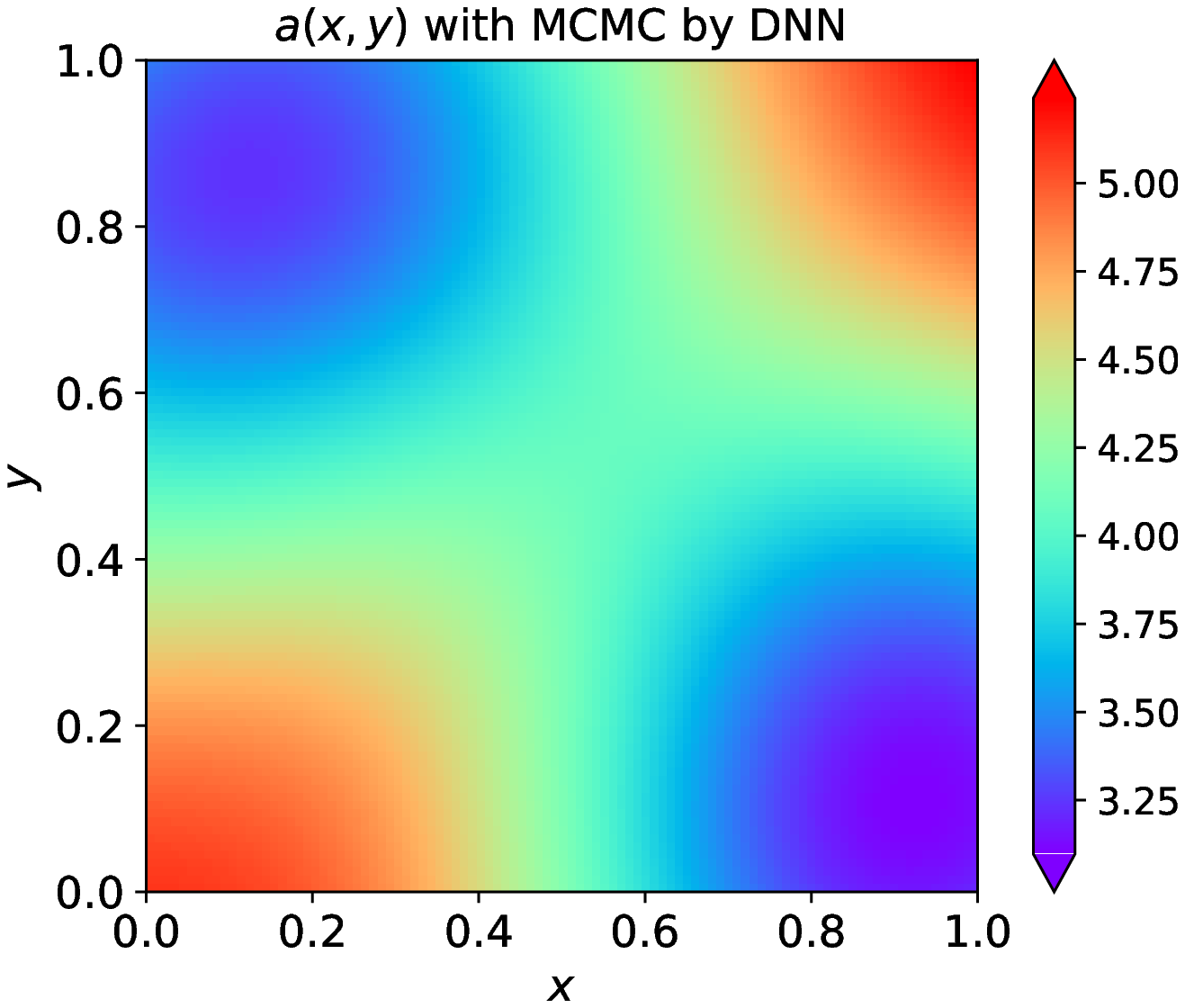}
		\end{minipage}
	}%
	\subfigure[Point-wise errors]{
		\begin{minipage}[t]{0.3\linewidth}
			\centering
			\includegraphics[width=2in]{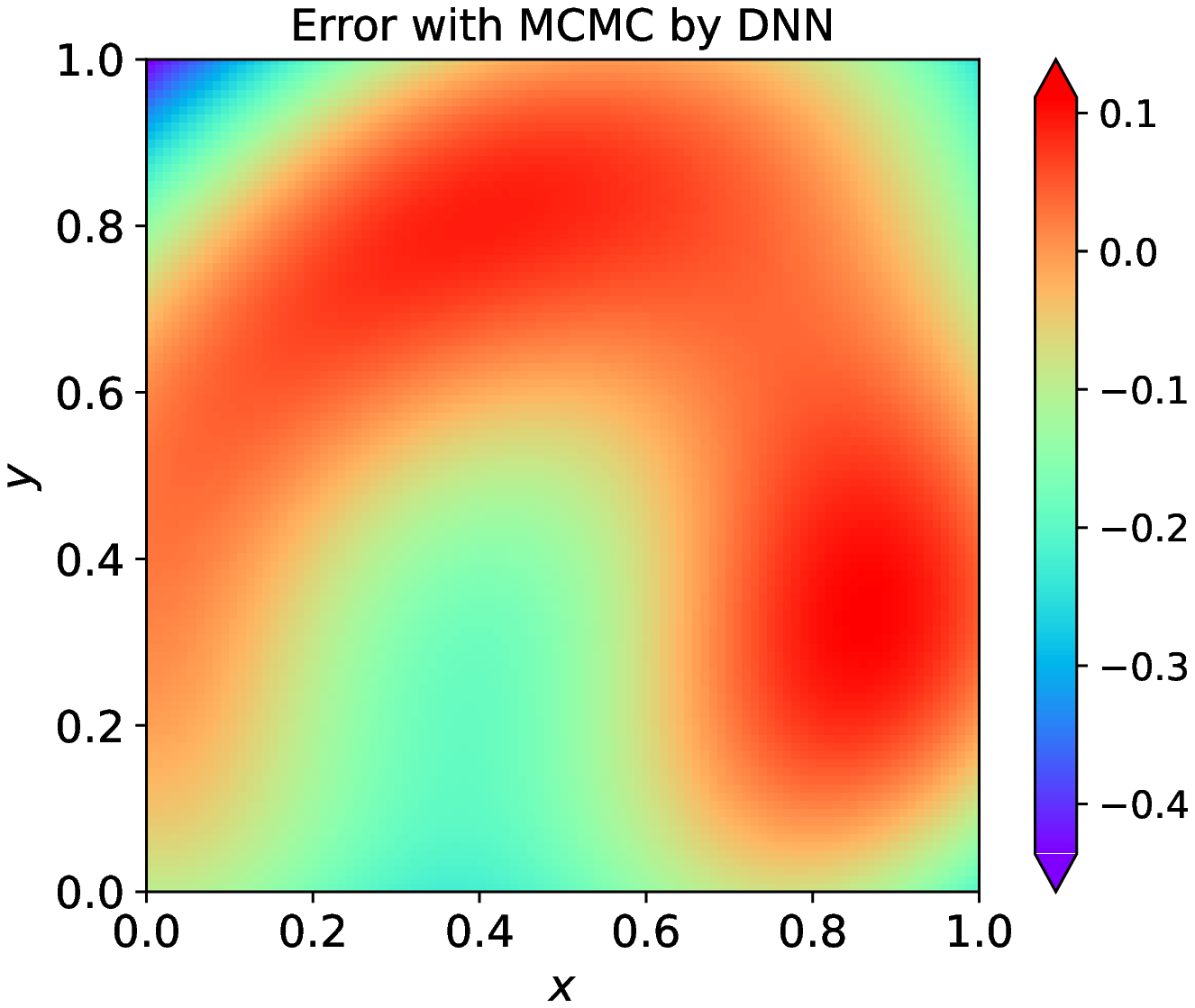}
		\end{minipage}
	}%

	\centering
	\caption{(a), (c) The reconstructed diffusion coefficient by MCMC+$G_{NN}^{(\alpha,a)}$ from the interior measurement data that was polluted by noise with 0-mean, and the standard deviation $\sigma$ is $0.001$ and $0.005$, respectively. (b), (d) The difference between exact diffusion coefficient and reconstruction diffusion coefficient by MCMC+$G_{NN}^{(\alpha,a)}$ from the interior measurement data that was polluted by noise with 0-mean, and the standard deviation $\sigma$ is $0.001$ and $0.005$, respectively. }
	\label{fig4.2_DNN_inDomain}
\end{figure}

In this part, we consider the inverse diffusion coefficient $a(x,y)$ in the subdiffusion problem (\ref{1.1}) based on the interior measurement data $u(x,y,t)$, where $(x,y)\in\omega=[0.25,0.75]^2$, $t\in [0,T]$ in (\ref{IP_4.0}).
Similarly, in the experiment, we fix $\alpha=0.7$, reaction coefficient $c(x,y)=-(xy+4)$, 
initial value $u_0(x,y)=6\sin(2\pi x)\sin(3\pi y)$, source $f(x,y)=\sin(3\pi x)\sin(\pi y)+6\exp(x^2+y^2)$.
Instead of retraining a new deep operator network for this inverse problem, we use the trained deep operator network $G_{NN}^{(\alpha,a)}$ which is used to approximate the solution operator (\ref{3.1.1}) to carry out the MCMC sampling.
Figure \ref{fig4.2_DNN_inDomain_diffu} shows the exact diffusion coefficient. Figure \ref{fig4.2_DNN_inDomain} shows the numerical inversion results obtained by MCMC sampling with the operator network $G_{NN}^{\alpha,a}$. According to the figure, using the MCMC sampling with the deep operator network $G_{NN}^{(\alpha,a)}$ can produce relatively accurate numerical inversion results. This numerical experiment demonstrates that, once trained, the deep operator network can be used to solve the different inverse diffusion coefficient problems of the subdiffusion problem (\ref{1.1}).

\section{Conclusion}\label{sec5}
A deep learning method for solving the forward problems and inverse problems in subdiffusion is presented. By using the deep operator network, we study two kinds of operator learning problems of the subdiffusion problem (\ref{1.1}) through various numerical experiments. Furthermore, we apply a Bayesian inversion method to solve several inverse problems in the subdiffusion problem. To overcome the problem of time-consuming Bayesian inference with conventional numerical methods, we propose a BINO approach for modelling a subdiffusion problem and solving a inverse diffusion coefficient problem in the paper. Several numerical experiments show that the operator learning method we provide can effectively solve the forward problems and the Bayesian inverse problems in the subdiffusion model.

\bigskip
\noindent{\bf Acknowledgments}
This work is sponsored by the National Key R\&D Program of China  Grant No. 2019YFA0709503 (Z. X.) and No. 2020YFA0712000 (Z. M.), the Shanghai Sailing Program (Z. X.), the Natural Science Foundation of Shanghai Grant No. 20ZR1429000  (Z. X.), the National Natural Science Foundation of China Grant No. 62002221 (Z. X.), the National Natural Science Foundation of China Grant No. 12101401 (Z. M.), the National Natural Science Foundation of China Grant No. 12031013 (Z. M.), Shanghai Municipal of Science and Technology Major Project No. 2021SHZDZX0102, and the HPC of School of Mathematical Sciences and the Student Innovation Center at Shanghai Jiao Tong University.
\medskip

\bibliographystyle{unsrt}
\bibliography{mybibfilesource}

\end{document}